\documentclass[12pt,twoside]{article}
\usepackage[hmargin=0.8in,vmargin=0.9in]{geometry}
\usepackage{fancyhdr}
\usepackage{graphicx}
\usepackage{subfigure}
\usepackage{amssymb}
\usepackage{amsmath}
\usepackage{amsfonts}
\usepackage{theorem}
\usepackage{mathrsfs}
\usepackage{mathtools}
\usepackage{bm}
\usepackage{color}
\usepackage{setspace}
\usepackage{exscale}
\usepackage{relsize}
\usepackage{epstopdf}
\usepackage{float}
\usepackage{cite}

\usepackage{nicefrac}
\usepackage{setspace}

\usepackage{booktabs,multirow} 
\usepackage{array} 
\usepackage{paralist} 
\usepackage{verbatim} 
\usepackage{subfigure} 

\theoremstyle{plain}			
\newtheorem{thm}{Theorem}[section]

\newtheorem{rmk}[thm]{Remark}
{\theorembodyfont{\rmfamily}}

\usepackage{tikz}
\usetikzlibrary{positioning}
\usepackage{xcolor}

\allowdisplaybreaks[1]

\numberwithin{equation}{section}
\numberwithin{figure}{section}
\numberwithin{table}{section}

\newcommand\eref[1]{(\ref{#1})}

\newcommand*\xbar[1]{%
  \hbox{%
    \vbox{%
      \hrule height 0.5pt 
      \kern0.4ex
      \hbox{%
        \kern-0.05em
        \ensuremath{#1}%
        \kern-0.00em
      }%
    }%
  }%
}

\newcommand{\mF}{\bm{F}}
\newcommand{\mK}{\bm{K}}

\newcommand{\mH}{\bm{H}}
\newcommand{\mU}{\bm{U}}
\newcommand{\mE}{\bm{E}}

\newcommand{\mW}{\bm{W}}
\newcommand{\mN}{\bm{N}}
\newcommand{\mo}{\bm{0}}
\newcommand{\mI}{\bm{I}}

\newcommand{\dt}{\Delta t}
\newcommand{\dx}{\Delta x}
\newcommand{\dy}{\Delta y}

\newcommand{\hf}{{\frac{1}{2}}}

\newcommand{\jph}{{j+\frac{1}{2}}}
\newcommand{\jmh}{{j-\frac{1}{2}}}
\newcommand{\kph}{{k+\frac{1}{2}}}
\newcommand{\kmh}{{k-\frac{1}{2}}}

\def\softd{{\leavevmode\setbox1=\hbox{d}%
          \hbox to 1.05\wd1{d\kern-0.4ex{\char039}\hss}}}

\title{Local Characteristic Decomposition of Equilibrium Variables for Hyperbolic Systems of Balance Laws}
\author{Shaoshuai Chu\thanks{Department of Mathematics, RWTH Aachen University, 52056 Aachen, Germany; {\tt chu@igpm.rwth-aachen.de}},
Alexander Kurganov\thanks{Department of Mathematics and Shenzhen International Center for Mathematics, Southern University of Science and
Technology, Shenzhen, 518055, China; {\tt alexander@sustech.edu.cn}}, Mingye Na\thanks{Department of Mathematics, Southern University of
Science and Technology, Shenzhen, 518055, China; {\tt 12131231@mail.sustech.edu.cn}},\\
Bao-Shan Wang\thanks{School of Mathematical Sciences \& Laboratory of Marine Mathematics, Ocean University of China, Qingdao, 266100, China;
{\tt wbs@ouc.edu.cn}}, and Ruixiao Xin\thanks{Department of Mathematics, Southern University of Science and Technology, Shenzhen, 518055,
China; {\tt 12331009@mail.sustech.edu.cn}}}

\begin{document}
\date{}
\maketitle
\begin{abstract}
This paper is concerned with high-order numerical methods for hyperbolic systems of balance laws. Such methods are typically based on
high-order piecewise polynomial reconstructions (interpolations) of the computed discrete quantities. However, such reconstructions
(interpolations) may be oscillatory unless the reconstruction (interpolation) procedure is applied to the local characteristic variables via
the local characteristic decomposition (LCD). Another challenge in designing accurate and stable high-order schemes is related to enforcing
a delicate balance between the fluxes, sources, and nonconservative product terms: a good scheme should be well-balanced (WB) in the sense
that it should be capable of exactly preserving certain (physically relevant) steady states. One of the ways to ensure that the
reconstruction (interpolation) preserves these steady states is to apply the reconstruction (interpolation) to the equilibrium variables,
which are supposed to be constant at the steady states. To achieve this goal and to keep the reconstruction (interpolation) non-oscillatory,
we introduce a new LCD of equilibrium variables. We apply the developed technique to the fifth-order Ai-WENO-Z interpolation implemented
within the WB A-WENO framework recently introduced in [{\sc S. Chu, A. Kurganov, and R. Xin}, Beijing J. of Pure and Appl. Math., 2 (2025),
pp. 87--113], and illustrate its performance on a variety of numerical examples.
\end{abstract}

\noindent
{\bf Key words:} High-order reconstructions (interpolations); A-WENO schemes; well-balanced schemes; equilibrium variables; local
characteristic decomposition.

\smallskip
\noindent
{\bf AMS subject classification:} 76M20, 65M06, 35L65, 35L67.

\section{Introduction}
This paper is focused on the development of a new local characteristic decomposition (LCD) of equilibrium variables for hyperbolic systems
of balance laws, which, in the two-dimensional (2-D) case, read as
\begin{equation}
\bm U_t+\bm F(\bm U)_x+\bm G(\bm U)_y=B^x(\bm U)\bm U_x+B^y(\bm U)\bm U_y+\bm S(\bm U).
\label{1.1}
\end{equation}
Here, $x$ and $y$ are the spatial variables, $t$ is time, $\mU\in\mathbb R^d$ is a vector of unknowns, $\mF$ and $\bm G$ are the fluxes,
$\bm S$ is the source term, and $B^x(\mU)\mU_x$ and $B^y(\mU)\mU_y$ are nonconservative product terms.

Development of high-order numerical methods for \eref{1.1} is a challenging task for the following two main reasons. First, solutions of
\eref{1.1} may develop discontinuities (even for infinitely smooth initial data) and thus high-order schemes should rely on non-oscillatory
high-order reconstructions (interpolations) of the solutions out of the computed discrete quantities. To make these reconstructions
(interpolations) non-oscillatory, one typically needs to use nonlinear limiting techniques often applied to the local characteristic
variables via the LCD; see, e.g., \cite{DLGW,JSZ,Joh,Liu17,Nonomura20,Qiu02,Shu20,WLGD_18}. The LCD is applied by computing local
linearizations of the matrices $\frac{\partial{\bm F}}{\partial{\bm U}}-B^x(\bm U)$ and
$\frac{\partial{\bm G}}{\partial{\bm U}}-B^y(\bm U)$, locally switching to the corresponding characteristic variables, performing the
reconstruction (interpolation) to these local variables, and then switching back to the original variables $\bm U$.

The use of characteristic-wise reconstructions has played an important role in the development of high-order WENO-type schemes for
hyperbolic systems. In particular, LCD is often essential for reducing spurious oscillations near shocks, contact discontinuities, and sharp
transition regions; see, e.g., \cite{JSZ,Joh,Nonomura20,Qiu02,Shu20} and references therein. Recently, other transformed-variable
approaches, such as WENO reconstructions based on Riemann invariants, have also been investigated as alternatives to the standard LCD in
certain systems; see, e.g., \cite{XuShu2024}. These works demonstrate that a suitable choice of variables for the nonlinear reconstruction
is crucial for the robustness of high-order schemes.

Second, many (physically relevant) solutions of \eref{1.1} are, in fact, small perturbations of certain steady states, which are supposed to
be exactly preserved by good high-order schemes---this is a so-called well-balanced (WB) property. A large body of work has been devoted to
the design of WB schemes for balance laws, including a variety of shallow water models, the Euler equations with gravitation, and other
systems with sources and nonconservative product terms. For the Euler equations with gravitation, for example, high-order WB finite-volume
(FV) and finite-difference (FD) WENO-type schemes have been developed; see, e.g.,
\cite{KaeppeliMishra2014,LiXing2016,LiXing2018,KlingenbergPuppoSemplice2019,GrosheintzLavalKaeppeli2019,LiWangDon2021} and references
therein. These methods are often based on local hydrostatic reconstructions, equilibrium-perturbation reconstructions, or carefully balanced
discretizations of the flux and source terms. In particular, characteristic-wise WB WENO schemes have been developed for the Euler equations
with gravitation in \cite{LiWangDon2021}.

One of the ways to ensure that the high-order reconstruction (interpolation) preserves these steady states is to reconstruct (interpolate)
the equilibrium variables instead of the conservative ones. We note that in many cases, the equilibrium variables can be obtained by
rewriting the system \eref{1.1} in the following form (see, e.g., \cite{CKL23}):
\begin{equation}
\bm U_t+M^x(\bm U)\bm E^x(\bm U)_x+M^y(\bm U)\bm E^y(\bm U)_y=\bm0,
\label{1.4}
\end{equation}
where
\begin{equation}
\begin{aligned}
&M^x(\bm U)\bm E^x(\bm U)_x=\bm F(\bm U)_x-B^x(\bm U)\bm U_x-\bm S^x(\bm U),\\[1.ex]
&M^y(\bm U)\bm E^y(\bm U)_y=\bm G(\bm U)_y-B^y(\bm U)\bm U_y-\bm S^y(\bm U),
\end{aligned}
\label{1.4a}
\end{equation}
and $\bm S^x(\bm U)+\bm S^y(\bm U)=\bm S(\bm U)$. In \eref{1.4a}, $M^x,M^y\in\mathbb R^{d\times d}$ and $\bm E^x(\bm U(x,y))$,
$\bm E^y(\bm U(x,y))$ are equilibrium variables, since
\begin{equation}
\bm E^x(\bm U(x,y))=\bm E^x(y)\quad\mbox{and}\quad\bm E^y(\bm U(x,y))=\bm E^y(x)
\label{1.3a}
\end{equation}
at steady states satisfying $\bm E^x(\bm U)_x=\bm E^y(\bm U)_y\equiv0$. Therefore, one may prefer to reconstruct $\bm E^x$ in the
$x$-direction and $\bm E^y$ in the $y$-direction and then to recalculate the corresponding values of $\mU$ to ensure that \eref{1.3a} is
satisfied at the discrete level.

The reconstruction of equilibrium variables and the use of LCD can be, in principle, viewed as two separate components of the numerical
solution algorithm. However, when the reconstruction is performed for the equilibrium variables rather than for the conservative variables,
the LCD based on the Jacobians of the conservative-variable formulation is no longer a natural characteristic decomposition for the
quantities being interpolated. To combine the WB reconstruction of equilibrium variables with the oscillation-suppressing effect of the LCD,
one needs to construct the characteristic decomposition associated with the equilibrium-variable formulation itself.

In this paper, we introduce a new LCD of equilibrium variables. To this end, we first rewrite the system \eref{1.4} in the following
two equivalent (for smooth solution) formulations:
\begin{equation}
\begin{aligned}
&\bm E^x(\bm U)_t+C^x(\bm U)\bm E^x(\bm U)_x+D^x(\bm U)\bm U_y=\tilde\mI^x(\mU),\\[1.ex]
&\bm E^y(\bm U)_t+D^y(\bm U)\bm U_x+C^y(\bm U)\bm E^y(\bm U)_y=\tilde\mI^y(\mU),
\end{aligned}
\label{1.5}
\end{equation}
where the matrices $C^x$ and $C^y$ and the source terms $\tilde{\bm I}^x$ and $\tilde{\bm I}^y$ are specified in \S\ref{sec3}. The proposed
LCD is based on the matrices $C^x$ and $C^y$, which are associated with the equilibrium-variable formulations \eref{1.5}. We then compute
the matrices $C^x$ and $C^y$ at the grid points and use them to compute the local characteristic equilibrium variables, which are
reconstructed (interpolated) to obtain high-order values of $\bm E^x$ and $\bm E^y$, which, in turn, give us the corresponding high-order
values of $\bm U$ (solving nonlinear systems of equations may be required).

We implement the new LCD technique in the framework of flux globalization based WB alternative weighted essentially non-oscillatory
(A-WENO) FD schemes recently introduced in \cite{CKX_24WB}. The local characteristic equilibrium variables are interpolated using the
fifth-order affine-invariant WENO-Z (Ai-WENO-Z) interpolations \cite{DLWW22,WD22,LLWDW23}. The developed A-WENO scheme is applied to
five systems of balance laws including the nozzle flow system, the one- and two-layer shallow water equations, the compressible Euler
equations with gravitation, and the 2-D Ripa system. We conduct several numerical experiments to demonstrate that the proposed LCD of
equilibrium variables reduces spurious oscillations while preserving the WB property of the underlying flux-globalization based scheme.

\section{Flux Globalization Based WB A-WENO Schemes: An Overview}
In this section, we give an overview of the flux globalization based WB A-WENO schemes introduced in\cite{CKX_24WB} for general nonconservative systems.

\subsection{1-D Scheme}
The one-dimensional (1-D) hyperbolic systems of balance laws
\begin{equation}
\bm U_t+\bm F(\bm U)_x=B(\bm U)\bm U_x+\bm S(\bm U)
\label{1.4b}
\end{equation}
can be written in an equivalent quasi-conservative form:
\begin{equation*}
\bm U_t+\bm K(\bm U)_x=\bm0,
\end{equation*}
where $\bm K(\bm U)$ is a global flux
\begin{equation*}
\bm K(\bm U)=\bm F(\bm U)-\bm R(\bm U),\quad\bm R(\bm U)=
\int\limits_{\hat x}^x\Big[B(\bm U(\xi,t))\bm U_\xi(\xi,t)+\bm S(\bm U(\xi,t))\Big]\,{\rm d}\xi,
\end{equation*}
and ${\hat x}$ is an arbitrary number.

We first introduce a uniform mesh consisting of the cells $[x_\jmh,x_\jph]$ of size $x_\jph-x_\jmh\equiv\dx$ centered at
$x_j=(x_\jmh+x_\jph)/2$, $j=1,\ldots,N$. We assume that at a certain time level $t$, the approximate solution, realized in terms of its cell
centered values $\bm U_j\approx\bm U(x_j,t)$, is available (in the rest of the paper, we will suppress the time-dependence of all of the
indexed quantities for the sake of brevity). The solution is then evolved in time by solving the following system of ODEs:
\begin{equation}
\frac{{\rm d}\bm U_j}{{\rm d}t}=-\frac{\bm{{\cal K}}_\jph-\bm{{\cal K}}_\jmh}{\dx},
\label{2.3}
\end{equation}
where $\bm{{\cal K}}_\jph$ are the fifth-order A-WENO numerical fluxes (see \cite{CKX23,CKX_24WB}):
\begin{equation*}
\bm{{\cal K}}_\jph=\bm{{\cal K}}^{\rm FV}_\jph-\frac{(\dx)^2}{24}(\bm K_{xx})_\jph+\frac{7(\dx)^4}{5760}(\bm K_{xxxx})_\jph.
\end{equation*}
Here, $\bm{{\cal K}}^{\rm FV}_\jph$ is a FV numerical flux (in the numerical experiments reported in \S\ref{sec4}, we have used
the second-order WB path-conservative central-upwind numerical flux introduced in \cite{KLX_21}), and $(\bm K_{xx})_\jph$ and
$({\bm K_{xxxx}})_\jph$ are the high-order correction terms. The numerical fluxes
$\bm{{\cal K}}^{\rm FV}_\jph=\bm{{\cal K}}^{\rm FV}_\jph(\mU_\jph^\pm,\widehat\mU_\jph^\pm)$ are computed using the one-sided interpolated
values of $\mU$, and to enforce the WB evolution, one needs to use two copies of those values denoted by $\mU_\jph^\pm$ and
$\widehat\mU_\jph^\pm$; see \cite{CKX_24WB,KLX_21} for details. The correction terms $(\bm K_{xx})_\jph$ and $({\bm K_{xxxx}})_\jph$ are
computed using the numerical fluxes $\bm{{\cal K}}^{\rm FV}_\jph$, which have been already obtained:
\begin{equation*}
\begin{aligned}
&(\mK_{xx})_\jph=\frac{1}{12(\dx)^2}\Big[-\bm{{\cal K}}^{\rm FV}_{j-\frac{3}{2}}+16\bm{{\cal K}}^{\rm FV}_\jmh-
30\bm{{\cal K}}^{\rm FV}_\jph+16\bm{{\cal K}}^{\rm FV}_{j+\frac{3}{2}}-\bm{{\cal K}}^{\rm FV}_{j+\frac{5}{2}}\Big],\\
&(\mK_{xxxx})_\jph=\frac{1}{(\dx)^4}\Big[\bm{{\cal K}}^{\rm FV}_{j-\frac{3}{2}}-4\bm{{\cal K}}^{\rm FV}_\jmh+6\bm{{\cal K}}^{\rm FV}_\jph-
4\bm{{\cal K}}^{\rm FV}_{j+\frac{3}{2}}+\bm{{\cal K}}^{\rm FV}_{j+\frac{5}{2}}\Big];
\end{aligned}
\end{equation*}
see \cite{CKX23} for details. We would like to stress that these correction terms are needed to increase the order of the resulting scheme
to the fifth order and that adding these terms typically does not cause oscillations as long as the reconstruction is performed in the local
characteristic variables; see, e.g., \cite{JSZ,CCK23_Adaptive,Liu17,WDGK_20,WLGD_18,WDKL}.

\subsection{2-D Scheme}
The 2-D hyperbolic systems of balance laws \eref{1.1} can be similarly written in an equivalent quasi-conservative form:
\begin{equation*}
\bm U_t+\bm K(\bm U)_x+\bm L(\bm U)_y=\bm0,
\end{equation*}
where $\bm K(\bm U)$ and $\bm L(\bm U)$ are the global fluxes
\begin{equation*}
\begin{aligned}
&\bm K(\bm U)=\bm F(\bm U)-\bm R^x(\bm U),&&\bm R^x(\bm U)=
\int\limits_{\hat x}^x\Big[B^x(\bm U(\xi,y,t))\bm U_\xi(\xi,y,t)+\bm S^x(\bm U(\xi,y,t))\Big]\,{\rm d}\xi,\\
&\bm L(\bm U)=\bm G(\bm U)-\bm R^y(\bm U),&&\bm R^y(\bm U)=
\int\limits_{\hat y}^y\Big[B^y(\bm U(x,\eta,t))\bm U_\eta(x,\eta,t)+\bm S^y(\bm U(x,\eta,t))\Big]\,{\rm d}\eta,
\end{aligned}
\end{equation*}
and ${\hat x}$ and ${\hat y}$ are arbitrary numbers.

Let $[x_\jmh,x_\jph]\times[y_\kmh,y_\kph]$, $j=1,\cdots,N_x$, $k=1,\cdots,N_y$ be the uniform 2-D cells centered at $(x_j,y_k)$ with
$x_\jph-x_\jmh\equiv\dx$, $y_\kph-y_\kmh\equiv\dy$, $x_j=(x_\jmh+x_\jph)/2$, and $y_k=(y_\kmh+y_\kph)/2$. We denote by
$\mU_{j,k}\approx\mU(x_j,y_k,t)$ the computed cell centered values, which are assumed to be available at a certain time level $t$. The point
values $\mU_{j,k}$ are evolved in time by solving the following system of ODEs:
\begin{equation}
\frac{{\rm d}\mU_{j,k}}{{\rm d}t}=-\frac{\bm{{\cal K}}_{\jph,k}-\bm{{\cal K}}_{\jmh,k}}{\dx}-
\frac{\bm{{\cal L}}_{j,\kph}-\bm{{\cal L}}_{j,\kmh}}{\dy},
\label{2.8a}
\end{equation}
where $\bm{{\cal K}}_{\jph,k}$ and $\bm{{\cal L}}_{j,\kph}$ are the fifth-order A-WENO numerical fluxes (see \cite{CKX23,CKX_24WB}):
\begin{equation*}
\begin{aligned}
&\bm{{\cal K}}_{\jph,k}=
\bm{{\cal K}}^{\rm FV}_{\jph,k}-\frac{(\dx)^2}{24}(\bm K_{xx})_{\jph,k}+\frac{7(\dx)^4}{5760}(\bm K_{xxxx})_{\jph,k},\\
&\bm{{\cal L}}_{j,\kph}=
\bm{{\cal L}}^{\rm FV}_{j,\kph}-\frac{(\dy)^2}{24}(\bm L_{yy})_{j,\kph}+\frac{7(\dy)^4}{5760}(\bm L_{yyyy})_{j,\kph}.
\end{aligned}
\end{equation*}
Here, $\bm{{\cal K}}^{\rm FV}_{\jph,k}$ and $\bm{{\cal L}}^{\rm FV}_{j,\kph}$ are finite-volume numerical fluxes and
$(\bm K_{xx})_{\jph,k}$, $({\bm K_{xxxx}})_{\jph,k}$, $(\bm L_{yy})_{j,\kph}$, and $({\bm L_{yyyy}})_{j,\kph}$ are the high-order correction
terms. The numerical fluxes $\bm{{\cal K}}^{\rm FV}_{\jph,k}=\bm{{\cal K}}^{\rm FV}_{\jph,k}\big(\mU_{\jph,k}^\pm\big)$ and
$\bm{{\cal L}}^{\rm FV}_{j,\kph}=\bm{{\cal L}}^{\rm FV}_{j,\kph}\big(\mU_{j,\kph}^\pm\big)$ are computed using the one-sided interpolated
values of $\mU$, and to enforce the WB evolution, one may need to use the numerical diffusion switch functions; see \cite{Kurganov25Na} for
an example of such switch function used for the compressible Euler equations with gravitation. The correction terms $(\bm K_{xx})_{\jph,k}$,
$({\bm K_{xxxx}})_{\jph,k}$, $(\bm L_{yy})_{j,\kph}$, and $({\bm L_{yyyy}})_{j,\kph}$ are computed using the numerical fluxes
$\bm{{\cal K}}^{\rm FV}_{\jph,k}$ and $\bm{{\cal L}}^{\rm FV}_{\jph,k}$, which have been already obtained:
\begin{equation*}
\begin{aligned}
&(\mK_{xx})_{\jph,k}=\frac{1}{12(\dx)^2}\Big[-\bm{{\cal K}}^{\rm FV}_{j-\frac{3}{2},k}+16\bm{{\cal K}}^{\rm FV}_{\jmh,k}-
30\bm{{\cal K}}^{\rm FV}_{\jph,k}+16\bm{{\cal K}}^{\rm FV}_{j+\frac{3}{2},k}-\bm{{\cal K}}^{\rm FV}_{j+\frac{5}{2},k}\Big],\\
&(\mK_{xxxx})_{\jph,k}=\frac{1}{(\dx)^4}\Big[\bm{{\cal K}}^{\rm FV}_{j-\frac{3}{2},k}-4\bm{{\cal K}}^{\rm FV}_{\jmh,k}+
6\bm{{\cal K}}^{\rm FV}_{\jph,k}-4\bm{{\cal K}}^{\rm FV}_{j+\frac{3}{2},k}+\bm{{\cal K}}^{\rm FV}_{j+\frac{5}{2},k}\Big];\\
&(\bm L_{yy})_{j,\kph}=\frac{1}{12(\dy)^2}\Big[-\bm{{\cal L}}^{\rm FV}_{j,k-\frac{3}{2}}+16\bm{{\cal L}}^{\rm FV}_{j,\kmh}-
30\bm{{\cal L}}^{\rm FV}_{j,\kph}+16\bm{{\cal L}}^{\rm FV}_{j,k+\frac{3}{2}}-\bm{{\cal L}}^{\rm FV}_{j,k+\frac{5}{2}}\Big],\\
&(\bm L_{yyyy})_{j,\kph}=\frac{1}{(\dy)^4}\Big[\bm{{\cal L}}^{\rm FV}_{j,k-\frac{3}{2}}-4\bm{{\cal L}}^{\rm FV}_{j,\kmh}+
6\bm{{\cal L}}^{\rm FV}_{j,\kph}- 4\bm{{\cal L}}^{\rm FV}_{j,k+\frac{3}{2}}+\bm{{\cal L}}^{\rm FV}_{j,k+\frac{5}{2}}\Big];
\end{aligned}
\end{equation*}
see \cite{CKX23} for details.

\section{Local Characteristic Decomposition of Equilibrium Variables}\label{sec3}
In this section, we apply the LCD approach introduced in \cite{DLGW,JSZ,Joh,Liu17,Nonomura20,Qiu02,Shu20,WLGD_18} to the Ai-WENO-Z interpolation of the equilibrium variables.

Throughout this section, we use the following notation for the equilibrium variables:
\begin{equation}\label{3.1a}
\mE^\nu(\mU)=\mW^\nu(\mU)+\mI^\nu(\mU),
\end{equation}
where $\nu$ stands for either $x$ or $y$, $\mW^\nu$ denotes the local part of the equilibrium variable $\mE^\nu$, and $\mI^\nu$ is its
global (integral) part. When the equilibrium variables are local, $\mI^\nu\equiv0$, and thus $\mE^\nu=\mW^\nu$. In the 1-D case, the
directional superscript $\nu$ is omitted so that we write
\begin{equation}
\mE(\mU)=\mW(\mU)+\mI(\mU).
\label{3.1b}
\end{equation}

\subsection{1-D Case}\label{sec31}
We first rewrite the system \eref{1.4b} in the following equivalent (for smooth solution) form:
\begin{equation}
\bm U_t+M(\bm U)\bm E(\bm U)_x=\bm0,
\label{1.5b}
\end{equation}
where
\begin{equation}
M(\bm U)\bm E(\bm U)_x=\bm F(\bm U)_x-B(\bm U)\bm U_x-\bm S(\bm U),
\label{1.4c}
\end{equation}
which vanishes at steady states. In \eref{1.4c}, $M\in\mathbb R^{d\times d}$ and $\bm E$ is the vector of equilibrium variables, which are
constant at steady states. We then rewrite \eref{1.5b} again as
\begin{equation}
\bm E(\bm U)_t+C(\bm U)\bm E(\bm U)_x=\tilde\mI(\mU).
\label{1.5a}
\end{equation}

We emphasize that \eref{1.5a} is not used as a time-evolution equation in the numerical method; it is introduced only to identify the matrix
$C(\mU)$ used in the LCD of equilibrium variables.

We now need to specify $C$ and $\tilde{\bm I}$ in \eref{1.5a}. To this end, we consider two possible cases.

\medskip
\noindent
$\bullet$ If $\mE(\mU)$ is a local equilibrium variable, that is $\mI(\mU)\equiv\mo$ in \eref{3.1b}, then we simply multiply \eref{1.4c} by
$\frac{\partial{\bm E}}{\partial{\bm U}}$ and obtain \eref{1.5a} with $C(\bm U)=\frac{\partial{\bm E}}{\partial{\bm U}}M(\bm U)$ and
$\tilde\mI(\bm U)\equiv\bm0$.

\smallskip
\noindent
$\bullet$ If $\mE(\mU)$ is a global equilibrium variable, that is, $\mI(\mU)\ne\mo$ and
\begin{equation}
\mI(\mU)=\int\limits^x\big[\bm H(\bm U)+N(\mU)\mU_x\big]\,{\rm d}x,
\label{3.1c}
\end{equation}
where $\bm H(\bm U)$ is a given vector function and $N(\mU)$ is a given matrix, then we proceed in a different way. We first differentiate
\eref{3.1b} and \eref{3.1c} with respect to $t$ to obtain
\begin{equation}
\mE(\mU)_t+\frac{\partial\mW}{\partial\mU}(\mU)\,\mU_t+\int\limits^x\bigg[\frac{\partial\mH}{\partial\mU}(\mU)\,\mU_t+N(\mU)_t\,\mU_x+
N(\mU)\mU_{xt}\bigg]\,{\rm d}x.
\label{3.1d}
\end{equation}
We then use \eref{1.5b} and the identity
\begin{equation*}
N(\mU)\mU_{xt}=(N(\mU)\mU_t)_x-N(\mU)_x\,\mU_t
\end{equation*}
to rewrite \eref{3.1d} as
\begin{equation}
\begin{aligned}
\mE(\mU)_t&+\bigg[\frac{\partial\mW}{\partial\mU}(\mU)+N(\mU)\bigg]M(\mU)\mE(\mU)_x\\
=&-\int\limits^x\bigg[\frac{\partial\mH}{\partial\mU}(\mU)M(\mU)\mE(\mU)_x-\bigg(\frac{\partial\mN_1}{\partial\mU}(\mU)\,\mU_x\dots
\frac{\partial\mN_d}{\partial\mU}(\mU)\,\mU_x\bigg)M(\mU)\mE(\mU)_x\\[0.5ex]
&\hspace*{1.0cm}+\bigg(\frac{\partial\mN_1}{\partial\mU}(\mU)M(\mU)\mE(\mU)_x\dots\frac{\partial\mN_d}{\partial\mU}(\mU)M(\mU)\mE(\mU)_x
\bigg)\bigg]\,{\rm d}x.
\end{aligned}
\label{3.1e}
\end{equation}
Here, we have used the following notation: $N=(\mN_1\dots\mN_d)$, where $\mN_i$ is the $i^{\rm th}$ column of the matrix $N$. Notice that
\eref{3.1e} corresponds to \eref{1.5a} with
$$
C(\mU)=\bigg[\frac{\partial \mW (\mU) }{\partial \mU}(\mU)+N(\mU) \bigg]M(\mU)
$$
and $\tilde\mI(\mU)$ given by the right-hand side of \eref{3.1e}.

In fact,  the details of $\tilde\mI$ are not important as this source term does not influence the LCD of equilibrium variables, which is
based on the matrix $C$ only.

\medskip
We then evaluate the matrix $C$ at the grid points to obtain the constant matrices $C_j:=C(\mU_j)$, which can be diagonalized using the
matrices $Q_j$ and $Q_j^{-1}$ to obtain $\Lambda_j=Q_j^{-1}C_jQ_j$, where $\Lambda_j$ is a diagonal matrix containing the eigenvalues of
$C_j$.

Next, we introduce the local characteristic equilibrium variables in the neighborhood of $x=x_j$:
\begin{equation}
\bm\Gamma_\ell=Q_j^{-1}\mE_\ell,\quad\ell=j\pm2,j\pm1,j,
\label{3.3}
\end{equation}
apply the fifth-order Ai-WENO-Z (or any other fifth-order WENO-type) interpolation to evaluate the values $\bm\Gamma^+_\jmh$ and
$\bm\Gamma^-_\jph$, and finally obtain
\begin{equation}
\bm E^\pm_{j\mp\hf}=Q_j\bm\Gamma^\pm_{j\mp\hf}.
\label{3.4}
\end{equation}
Equipped with these values, we proceed as in \cite{CKX_24WB,KLX_21} and solve the nonlinear equations (see
\cite[equations (2.8) and (2.16)]{CKLX_22}) to recover the values $\mU^\pm_{j\mp\hf}$ and $\widehat\mU^\pm_{j\mp\hf}$ needed to evaluate the
numerical fluxes $\bm{{\cal K}}^{\rm FV}_\jph$.
\begin{rmk}
It should be pointed out that the presented LCD-based reconstruction algorithm is different from the one used in, e.g., \cite{Chu21} as the
LCD process is now performed at the cell centers $x=x_j$ not at the cell interfaces $x_\jph$. The current approach has three advantages.
First, no averaged values of any quantities from cells $j$ and $j+1$ are required. Second, only five---not six---values of $\bm\Gamma$
should be computed for every $j$. Third, we can use only one Ai-WENO-Z interpolant in every cell $C_j$ to evaluate $\bm\Gamma^+_\jmh$ and
$\bm\Gamma^-_\jph$, while in the LCD algorithm in \cite{Chu21}, one had to use two Ai-WENO-Z interpolants to compute the one-sided values of
$\bm\Gamma^\pm_\jph$.
\end{rmk}

\subsection{2-D Case}
In the 2-D case, the LCD of the equilibrium variables is based on the formulation of the 2-D system \eref{1.1} given in \eref{1.5}. In fact,
to perform the LCD, we will only need to specify $C^x$ and $C^y$ since the terms with $D^x$, $D^y$, $\tilde\mI^x$, and $\tilde\mI^y$ do not
influence the LCD process. The computation of $C^x$ and $C^y$ is similar to the computation of the matrix $C$ in \S\ref{sec31}, and they are
given by
\begin{equation}
C^\nu(\mU)=\bigg[\frac{\partial\mW^\nu}{\partial\mU}(\mU)+N^\nu(\mU)\bigg]M^\nu(\mU),\quad\nu\in\{x,y\},
\label{3.4aa}
\end{equation}
where the matrices $N^\nu(\mU)$ are from the expressions of the equilibrium variables \eref{3.1a} with
$$
\mI^x=\int\limits^x\bigg[\mH^x(\mU)+N^x(\mU)\,\mU_x \bigg]\,{\rm d}x,\quad
\mI^y=\int\limits^y\bigg[\mH^y(\mU)+N^y(\mU)\,\mU_y \bigg]\,{\rm d}y.
$$
Here, $\mH^x(\mU)$ and $\mH^y(\mU)$ are given vector functions, and $N^x(\mU)$ and $N^y(\mU)$ are given matrices. In particular, if the
equilibrium variables are local, then \eref{3.4aa} reduces to
$$
C^\nu(\mU)=\frac{\partial\mE^\nu}{\partial\mU}(\mU)M^\nu(\mU),\quad\nu\in\{x,y\}.
$$

We then evaluate the matrices $C^x$ and $C^y$ at the grid points to obtain the constant matrices $C^x_{j,k}:=C^x(\mU_{j,k})$ and
$C^y_{j,k}:=C^y(\mU_{j,k})$, which can be diagonalized using the matrices $Q^x_{j,k}$ and $Q^y_{j,k}$ to obtain
$\Lambda^x_{j,k}=(Q^x_{j,k})^{-1}C^x_{j,k}Q^x_{j,k}$ and $\Lambda^y_{j,k}=(Q^y_{j,k})^{-1}C^y_{j,k}Q^y_{j,k}$, where $\Lambda^x_{j,k}$ and
$\Lambda^y_{j,k}$ are the diagonal matrices containing the eigenvalues of $C^x_{j,k}$ and $C^y_{j,k}$, respectively.

Next, we introduce the local characteristic equilibrium variables in the neighborhood of $(x,y)=(x_j,y_k)$:
\begin{equation}
\bm\Gamma^x_{\ell,k}=(Q^x_{j,k})^{-1}\mE^x_{\ell,k},\quad\ell=j\pm2,j\pm1,j,\quad\mbox{and}\quad
\bm\Gamma^y_{j,\ell}=(Q^y_{j,k})^{-1}\mE^y_{j,\ell},\quad\ell=k\pm2,k\pm1,k,
\label{3.3a}
\end{equation}
apply the fifth-order Ai-WENO-Z (or any other fifth-order WENO-type) interpolation in the $x$- and $y$-directions to evaluate the values
$(\bm\Gamma^x_{j\mp\hf,k})^\pm$ and $(\bm\Gamma^y_{j,k\mp\hf})^\pm$, respectively, and finally obtain
\begin{equation}
\big(\bm E^x_{j\mp\hf,k}\big)^\pm=Q^x_{j,k}\big(\bm\Gamma^x_{j\mp\hf,k}\big)^\pm,\quad
\big(\bm E^y_{j,k\mp\hf}\big)^\pm=Q^y_{j,k}\big(\bm\Gamma^y_{j,k\mp\hf}\big)^\pm.
\label{3.4a}
\end{equation}

Next, we show the application of the introduced LCD of equilibrium variables to five particular systems of balance laws.
\begin{rmk}
We emphasize that the proposed LCD does not change the class of steady states preserved by the underlying flux-globalization based WB
scheme. Its role is to perform the nonlinear reconstruction in local characteristic equilibrium variables, thereby reducing spurious
oscillations while maintaining the WB property. In the 1-D case, the scheme preserves all steady states that can be represented by
$M(\mU)\mE(\mU)_x=\mo$, that is, all equilibria for which the equilibrium variables $\mE(\mU)$ are constant. In the 2-D case, the
corresponding preserved equilibria are those satisfying $\mE^x(\mU)_x=\mE^y(\mU)_y\equiv\mo$. Equivalently, at such steady states,
$\mE^x(\mU(x,y))$ is independent of $x$ and $\mE^y(\mU(x,y))$ is independent of $y$. We would like to stress that the studied method is not
designed to preserve a-priori prescribed equilibria, but it is capable of maintaining a whole family of equilibria associated with the
equilibrium variables. For the particular systems studied below, these equilibrium variables are explicitly specified in
\S\ref{sec3.3}--\S\ref{sec3.8}.
\end{rmk}

\subsection{Application to the Nozzle Flow System}\label{sec3.3}
In this section, we consider the 1-D nozzle flow system, which reads as \eref{1.4b} with
\begin{equation*}
\bm U=(\sigma\rho,\sigma\rho u)^\top,\quad\bm F(\bm U)=(\sigma\rho u,\sigma\rho u^2+\sigma p)^\top,\quad
B(\bm U)=0,\quad\bm S(\bm U)=(0,p\sigma_x)^\top,
\end{equation*}
where $\rho$ is the density, $u$ is the velocity, $p(\rho)=\kappa\rho^\gamma$ is the pressure, $\kappa>0$ and $1<\gamma<\frac{5}{3}$ are
constants, and $\sigma=\sigma(x)$ denotes the cross-section of the nozzle. The studied nozzle flow system admits steady-state solutions
satisfying $M(\bm U)\bm E(\bm U)_x=\bm0$ with
\begin{equation*}
\begin{aligned}
&M(\bm U)=\begin{pmatrix}1&0\\u&\sigma\rho\end{pmatrix},\quad\bm E(\bm U)=\begin{pmatrix}q\\{\cal E}\end{pmatrix},\quad
q=\sigma\rho u,\quad{\cal E}=\frac{u^2}{2}+\frac{\kappa\gamma}{\gamma-1}\rho^{\gamma-1}.
\end{aligned}
\end{equation*}

In order to apply the fifth-order Ai-WENO-Z interpolation to the equilibrium variables, we first compute
\begin{equation*}
{\cal E}_j=\frac{u_j^2}{2}+\frac{\kappa\gamma}{\gamma-1}(\rho_j)^{\gamma-1},
\end{equation*}
where $u_j=q_j/(\sigma\rho)_j$, $\rho_j=(\sigma\rho)_j/\sigma_j$, and $\sigma_j=\sigma(x_j)$, and then evaluate the matrices
\begin{equation*}
\begin{aligned}
&C_j=\begin{pmatrix}u_j&(\sigma\rho)_j\\\dfrac{\kappa\gamma(\rho_j)^{\gamma-1}}{(\sigma\rho)_j}&u_j\end{pmatrix},\quad
Q_j=\begin{pmatrix}(\sigma\rho)_j&(\sigma\rho)_j\\
-\sqrt{\kappa\gamma}(\rho_j)^{\frac{\gamma-1}{2}}&\sqrt{\kappa\gamma}(\rho_j)^{\frac{\gamma-1}{2}}\end{pmatrix},\\
&Q_j^{-1}=\frac{1}{2\sqrt{\kappa\gamma}(\sigma\rho)_j(\rho_j)^{\frac{\gamma-1}{2}}}
\begin{pmatrix}\sqrt{\kappa\gamma}(\rho_j)^{\frac{\gamma-1}{2}}&-(\sigma\rho)_j\\
\sqrt{\kappa\gamma}(\rho_j)^{\frac{\gamma-1}{2}}&(\sigma\rho)_j\end{pmatrix}.
\end{aligned}
\end{equation*}
We now implement the LCD of the equilibrium variables followed by the fifth-order Ai-WENO-Z interpolation giving $\bm\Gamma^\pm_{j\mp\hf}$
and then $q^\pm_{j\mp\hf}$ and ${\cal E}^\pm_{j\mp\hf}$. After that, we apply the same fifth-order Ai-WENO-Z interpolation to obtain the
one-sided values of the cross-section of the nozzle $\sigma^\pm_{j\mp\hf}$, and then solve the nonlinear equations as described in
\cite[Equations (3.5) and (3.6)]{CKX_24WB} to obtain $(\sigma\rho)^\pm_{j\mp\hf}$ and $(\widehat{\sigma\rho})^\pm_{j\mp\hf}$.

\subsection{Application to the Saint-Venant System with Manning Friction}\label{sec32}
In this section, we consider the 1-D Saint-Venant system of shallow water equations with Manning friction, which reads as \eref{1.4b} with
\begin{equation*}
\bm U=(h,q)^\top,\quad\bm F(\bm U)=\Big(q,hu^2+\hf gh^2\Big)^\top,\quad B(\bm U)=0,\quad\bm S(\bm U)=(0,-ghZ_x-ghS_f)^\top,
\end{equation*}
where $h$ is the water depth, $u$ is the velocity, $q=hu$ represents the discharge, $Z(x)$ is a function describing the bottom topography,
which can be discontinuous, $g$ is the constant acceleration due to gravity, $S_f$ is the Manning friction term (see, e.g., \cite{Man91})
given by $S_f=n^2q|q|h^{-\frac{10}{3}}$. The studied Saint-Venant system admits steady-state solutions satisfying
$M(\bm U)\bm E(\bm U)_x=\bm0$ with
\begin{equation*}
\begin{aligned}
M(\bm U)=\begin{pmatrix}1&0\\u&h\end{pmatrix},\quad\bm E(\bm U)=\begin{pmatrix}q\\{\cal E}\end{pmatrix},\quad
{\cal E}=\frac{u^2}{2}+g(h+Z)+\int\limits_{\hat x}^x{gS_f}\,{\rm d}\xi.
\end{aligned}
\end{equation*}

In order to apply the fifth-order Ai-WENO-Z interpolation to the equilibrium variables, we first compute
\begin{equation*}
{\cal E}_j=\frac{u_j^2}{2}+g(h_j+Z_j)+{\cal I}_j,
\end{equation*}
where $u_j=q_j/h_j$, $Z_j=Z(x_j)$, and ${\cal I}_j$ is a fifth-order approximation of the integral
$\int_{x_{-\frac52}}^{x_j}gS_f\,{\rm d}x$, in which we have set $\hat x=x_{-\frac{5}{2}}$. The values ${\cal I}_j$ are computed as follows.
First, we evaluate $f_j:=g(S_f)_j={gn^2q_j|q_j|}{h_j^{-\frac{10}{3}}}$ for $j=1,\ldots,N$ and extend these values to $j=-5,\ldots,0$ and
$j=N+1,\ldots,N+5$ using the prescribed boundary conditions (for $h$ and $q$) implemented within the ghost cell framework. We then construct
the interpolating polynomial for $f$ using the five points $(x_{-\frac{5}{2}},f_{-\frac{5}{2}}^+)$, $(x_{-\frac{9}{4}},f_{-\frac{9}{4}})$,
$(x_{-2},f_{-2})$, $(x_{-\frac{7}{4}},f_{-\frac{7}{4}})$, and $(x_{-\frac{3}{2}},f_{-\frac{3}{2}}^-)$ and integrate it over the interval
$[x_{-\frac{5}{2}},x_{-2}]$ to obtain
\begin{equation}
{\cal I}_{-2}=\frac{\dx}{360}\left[29f_{-\frac{5}{2}}^++124f_{-\frac{9}{4}}+24f_{-2}+4f_{-\frac{7}{4}}-f_{-\frac{3}{2}}^-\right],
\label{3.5a}
\end{equation}
where $f_{j\pm\frac{1}{4}}$ are evaluated as in \cite[\S4]{Chu21}, but now we use the Ai-WENO-Z interpolant instead of the WENO-Z one
implemented in \cite{Chu21}.

We then proceed recursively: construct the interpolating polynomials for $f$ using $(x_{j-1},f_{j-1})$,
$(x_{j-\frac{3}{4}},f_{j-\frac{3}{4}})$, $(x_\jmh,f_\jmh)$, $(x_{j-\frac{1}{4}},f_{j-\frac{1}{4}})$, and $(x_j,f_j)$, integrate them over
the corresponding interval $[x_{j-1},x_j]$, and end up with
\begin{equation}
{\cal I}_j={\cal I}_{j-1}+\frac{\dx}{90}\left[7f_{j-1}+32f_{j-\frac{3}{4}}+12f_\jmh+32f_{j-\frac{1}{4}}+7f_{j}\right],\quad j=1,\ldots,N+3,
\label{3.5b}
\end{equation}
where the point values $f_{j-\frac{3}{4}}$ and $f_{j-\frac{1}{4}}$ are computed as in \cite[\S4]{Chu21}, but using the Ai-WENO-Z
interpolant and $f_{\jmh}:=\big(f_\jmh^-+f_\jmh^+\big)/2$.

We then evaluate the matrices
\begin{equation*}
\begin{aligned}
C_j=\begin{pmatrix}u_j&h_j\\g&u_j\end{pmatrix},\quad Q_j=\begin{pmatrix}\sqrt{h_j}&\sqrt{h_j}\\-\sqrt{g}&\sqrt{g}\end{pmatrix},\quad
Q_j^{-1}=\frac{1}{2\sqrt{gh_j}}\begin{pmatrix}\sqrt{g}&-\sqrt{h_j}\\\sqrt{g}&\sqrt{h_j}\end{pmatrix},
\end{aligned}
\end{equation*}
and implement the LCD of the equilibrium variables followed by the fifth-order Ai-WENO-Z interpolation giving $\bm\Gamma^\pm_{j\mp\hf}$ and
then $q^\pm_{j\mp\hf}$ and ${\cal E}^\pm_{j\mp\hf}$. After that, we apply the same fifth-order Ai-WENO-Z interpolation to obtain the
one-sided values of the bottom topography $Z^\pm_{j\mp\hf}$, and then solve the nonlinear equations
\begin{equation}
{\cal E}^\pm_\jph=\frac{\big(q^\pm_\jph\big)^2}{2\big(h^\pm_\jph\big )^2}+g\left(h^\pm_\jph+Z^\pm_\jph\right)+{\cal I}_\jph,
\label{3.1}
\end{equation}
and
\begin{equation}
{\cal E}^\pm_\jph=\frac{\big(q^\pm_\jph\big)^2}{2\big(\,\widehat h^\pm_\jph\big)^2}+g\left(\widehat h^\pm_\jph+Z_\jph\right)+
{\cal I}_\jph,\quad Z_\jph=\hf\left(Z^+_\jph+Z^-_\jph\right)
\label{3.2}
\end{equation}
to obtain $h^\pm_\jph$ and $\,\widehat h^\pm_\jph$, respectively. Here, the integrals ${\cal I}_\jph$ are evaluated recursively by the
fifth-order quadrature: we first set ${\cal I}_{-\frac52}=0$ and then compute
\begin{equation}
{\cal I}_\jph={\cal I}_\jmh+\frac{\dx}{90}\left[7f_\jmh^++32f_{j-\frac{1}{4}}+12f_j+32f_{j+\frac{1}{4}}+7f_\jph^-\right]
\label{3.5c}
\end{equation}
for $j=-2,\ldots,N+2$, where we have used the quadrature, which is obtained by constructing the interpolating polynomials for $f$ using the
five points $(x_\jmh,f_\jmh^+)$, $(x_{j-\frac{1}{4}},f_{j-\frac{1}{4}})$, $(x_j,f_j)$, $(x_{j+\frac{1}{4}},f_{j+\frac{1}{4}})$, and
$(x_\jph,f_\jph^-)$ and integrating them over the corresponding intervals $[x_\jmh,x_\jph]$. As before, the values $f_{j\pm\frac{1}{4}}$ are
computed as in \cite[\S4]{Chu21}, but using the Ai-WENO-Z interpolant. Notice that equations \eref{3.1} and \eref{3.2} are cubic and we
solve them exactly as described in \cite{CKLX_22}.

\subsection{Application to the Two-Layer Shallow Water System}
In this section, we consider the 1-D two-layer shallow water system, which reads as \eref{1.4b} with
$$
\begin{aligned}
&\bm U=(h_1,q_1,h_2,q_2)^\top,\\
&\bm F(\bm U)=\big(q_1,h_1u_1^2+\dfrac{g}{2}h_1^2,q_2,h_2u_2^2+\dfrac{g}{2}h_2^2\big)^\top,\\
&\bm S(\bm U)=(0,-gh_1Z_x,0,-gh_2Z_x)^\top,
\end{aligned}
\qquad B(\bm U)=\begin{pmatrix}0&0&0&0\\0&0&-gh_1&0\\0&0&0&0\\-rgh_2&0&0&0\end{pmatrix}.
$$
Here, $h_1$ and $h_2$ are the water depths in the upper and lower layers, respectively, $u_1$ and $u_2$ are the corresponding velocities,
$q_1=h_1u_1$ and $q_2=h_2u_2$ represent the corresponding discharges, $Z(x)$ and $g$ are the same as in \S\ref{sec32}, and
$r=\frac{\rho_1}{\rho_2}<1$ is the ratio of the constant densities $\rho_1$ (upper layer) and $\rho_2$ (lower layer). The studied two-layer
shallow water system admits steady-state solutions satisfying $M(\bm U)\bm E(\bm U)_x=\bm0$ with
$$
\begin{aligned}
M(\bm U)=\begin{pmatrix}1&0&0&0\\u_1&h_1&0&0\\0&0&1&0\\0&0&u_2&h_2\end{pmatrix},\quad
\bm E(\bm U)=\begin{pmatrix}q_1\\{\cal E}_1\\q_2\\{\cal E}_2\end{pmatrix},\qquad
\begin{aligned}
&{\cal E}_1:=\frac{q_1^2}{2h_1^2}+g(h_1+h_2+Z),\\
&{\cal E}_2:=\frac{q_2^2}{2h_2^2}+g(rh_1+h_2+Z).
\end{aligned}
\end{aligned}
$$

In order to apply the fifth-order Ai-WENO-Z interpolation to the equilibrium variables, we first compute
\begin{equation*}
\begin{aligned}
({\cal E}_1)_j&=\frac{(q_1)_j^2}{2(h_1)_j^2}+g\big[(h_1)_j+(h_2)_j+Z_j\big],\\
({\cal E}_2)_j&=\frac{(q_2)_j^2}{2(h_2)_j^2}+g\big[r(h_1)_j+(h_2)_j+Z_j\big],
\end{aligned}
\end{equation*}
and evaluate the matrices
\begin{equation*}
\begin{aligned}
C_j=\begin{pmatrix}(u_1)_j&(h_1)_j&0&0\\g&(u_1)_j&g&0\\0&0&(u_2)_j&(h_2)_j\\rg&0&g&(u_2)_j\end{pmatrix},
\end{aligned}
\end{equation*}
where $(u_1)_j=(q_1)_j/(h_1)_j$ and $(u_2)_j=(q_2)_j/(h_2)_j$. We then compute the matrices $Q_j$ and $Q_j^{-1}$ numerically and implement
the LCD of the equilibrium variables followed by the fifth-order Ai-WENO-Z interpolation giving $\bm\Gamma^\pm_{j\mp\hf}$ and then
$\bm E^\pm_{j\mp\hf}
=\big((q_1)^\pm_{j\mp\hf},({\cal E}_1)^\pm_{j\mp\hf},(q_2)^\pm_{j\mp\hf},({\cal E}_2)^\pm_{j\mp\hf}\big)^\top$. After that, we apply the
same fifth-order Ai-WENO-Z interpolation to obtain the one-sided values of the bottom topography $Z^\pm_{j\mp\hf}$, and then solve the
nonlinear equations as described in \cite[Equations (3.11)--(3.14)]{CKX_24WB} to obtain $(h_i)^\pm_\jph$ and $(\widehat h_i)^\pm_\jph$,
$i=1,2$.

\subsection{Application to the 1-D Euler Equations with Gravitation}
In this section, we consider the 1-D compressible Euler equations with gravitation, which can be written as \eref{1.4b} with
\begin{equation*}
\bm U=(\rho,m,{\cal E})^\top,\quad\bm F(\bm U)=\left(m,\rho u^2+p,u({\cal E}+p)\right)^\top,\quad B(\bm U)\equiv0,\quad
\bm S(\bm U)=(0,-\rho \phi_x,0)^\top,
\end{equation*}
where $\rho$ is the density, $u$ is the velocity, $m:=\rho u$ is the momentum, ${\cal E}:=E+\rho\phi$, $E$ is the total energy, $p$ is the
pressure, and $\phi(x)$ is the time-independent gravitational potential. The system is completed with an equation of state (EOS), which, in
the case of ideal gas, reads as $E=\frac{p}{\gamma-1}+\hf\rho u^2,$ where $\gamma$ is a specific heat ratio. The studied 1-D Euler equations
with gravitation admit steady-state solutions satisfying $M(\bm U)\bm E(\bm U)_x=\bm0$ with
\begin{equation*}
M(\bm U)=\begin{pmatrix}1&0&0\\0&1&0\\L&0&m\end{pmatrix},\quad\bm E(\bm U)=\begin{pmatrix}m\\K\\L\end{pmatrix},\quad
K:=\rho u^2+p+\int\limits_{\hat x}^x{\rho\phi_x}\,{\rm d}\xi,\quad L:=\frac{{\cal E}+p}{\rho}.
\end{equation*}

In order to apply the fifth-order Ai-WENO-Z interpolation to the equilibrium variables, we first compute
\begin{equation*}
K_j=\rho_ju_j^2+p_j+{\cal I}_j,\quad L_j=\frac{{\cal E}_j+p_j}{\rho_j}
\end{equation*}
where ${\cal I}_j$ is a fifth-order approximation of the integral $\int_{x_{-\frac52}}^{x_j}\rho\phi_x\,{\rm d}x$, in which we have set
$\hat x=x_{-\frac{5}{2}}$. The values ${\cal I}_j$ are computed as follows. First, we evaluate $f_j:=\rho_j\phi_x(x_j)$ for $j=1,\ldots,N$
and extend these values to $j=-5,\ldots,0$ and $j=N+1,\ldots,N+5$ using the prescribed boundary conditions implemented within the ghost cell
framework. We then compute the values ${\cal I}_j$, $j=-2,\dots,N+3$ using \eref{3.5a}--\eref{3.5b}.

We then evaluate the matrices
\allowdisplaybreaks
\begin{align*}
&C_j=\begin{pmatrix}0&1&0\\(\gamma-2)u_j^2+c_j^2&(3-\gamma)u_j&(\gamma-1)m_j\\[1.2ex]
\dfrac{(\gamma-1)u_j^2+c_j^2}{\rho_j}&\dfrac{(1-\gamma)u_j}{\rho_j}&\gamma u_j\end{pmatrix},\quad
Q_j=\begin{pmatrix}\dfrac{(1-\gamma)m_j}{c_j^2}&-\dfrac{\rho_j}{c_j}&\dfrac{\rho_j}{c_j}\\[1.5ex]
\dfrac{(1-\gamma)m_j^2}{c_j^2\rho_j}&\rho_j-\dfrac{m_j}{c_j}&\rho_j+\dfrac{m_j}{c_j}\\[1.5ex]1&1&1\end{pmatrix},\\
&Q_j^{-1}=\begin{pmatrix}\dfrac{u_j}{\rho_j}&-\dfrac{1}{\rho_j}&1\\[1.5ex]
\dfrac{(1-\gamma)u_j^2}{2c_j\rho_j}-\dfrac{1}{2\rho_j}(u_j+c_j)&\dfrac{(\gamma-1)u_j}{2c_j\rho_j}+\dfrac{1}{2\rho_j}&
\dfrac{(1-\gamma)u_j}{2 c_j}\\[1.8ex]
\dfrac{(\gamma-1)u_j^2}{2c_j\rho_j}-\dfrac{1}{2\rho_j}(u_j-c_j)&\dfrac{(1-\gamma)u_j}{2c_j\rho_j}+\dfrac{1}{2\rho_j}&
\dfrac{(\gamma-1)u_j}{2c_j}\end{pmatrix},
\end{align*}
where $c_j:=\sqrt{\gamma p_j/\rho_j}$, and implement the LCD of the equilibrium variables followed by the fifth-order Ai-WENO-Z
interpolation giving $\bm\Gamma^\pm_{j\mp\hf}$ and then $m^\pm_{j\mp\hf}$, $K^\pm_{j\mp\hf}$, and $L^\pm_{j\mp\hf}$. After that, we solve
the quadratic equations
\begin{equation}
(\gamma-1)\big(L_\jph^\pm-\phi(x_\jph)\big)\big(\rho_\jph^\pm\big)^2-\gamma\big(K_\jph^\pm-{\cal I}_\jph\big)\rho_\jph^\pm+
\frac{(\gamma+1)\big(m_\jph^\pm\big)^2}{2}=0
\label{3.13}
\end{equation}
using the method described in \cite[Appendix B.1]{Kurganov25Na} to obtain $\rho^\pm_{j\mp\hf}$ and then ${\cal E}^\pm_{j\mp\hf}$. In
\eref{3.13}, the integrals ${\cal I}_\jph$ are evaluated recursively using \eref{3.5c}. For the numerical flux
$\bm{{\cal K}}^{\rm FV}_\jph$, we use the one given in \cite[(2.2)]{Kurganov25Na}.

\subsection{Application to the 2-D Euler Equations with Gravitation}\label{sec3.7}
In this section, we consider the 2-D Euler equations with gravitation, which read as \eref{1.1} with
\begin{equation*}
\begin{aligned}
&\bm U=(\rho,m,n,{\cal E})^\top\!,~\bm F(\bm U)=\Big(m,\rho u^2+p,\frac{mn}{\rho},u({\cal E}+p)\Big)^\top\!,~
\bm G(\bm U)=\Big(n,\frac{mn}{\rho},\rho v^2+p,v({\cal E}+p)\Big)^\top\!,\\
&B^x(\bm U)=B^y(\bm U)\equiv0,\quad\bm S(\bm U)=(0,-\rho\phi_x,-\rho\phi_y,0)^\top,
\end{aligned}
\end{equation*}
where $v$ is the $y$-velocity, $n:=\rho v$ is the $y$-momentum, and the rest of the notation is the same as in the 1-D case. The system is
completed with the EOS for an ideal gas: $E=\frac{p}{\gamma-1}+\hf\rho(u^2+v^2)$.

The studied 2-D Euler equations with gravitation admit steady-state solutions satisfying \\
$M^x(\bm U)\bm E^x(\bm U)_x=\bm0$ and $M^y(\bm U)\bm E^y(\bm U)_y=\bm0$ with
$$
\begin{aligned}
&M^x(\bm U)=
\begin{pmatrix}1&0&0&0\\0&0&1&0\\v+(\gamma+1)u\psi^x&u+(1-\gamma)v\psi^x&-\gamma\psi^x&(\gamma-1)\rho\psi^x\\L&0&0&m\end{pmatrix},\\[0.5ex]
&M^y(\bm U)=
\begin{pmatrix}0&1&0&0\\v+(1-\gamma)u\psi^y&u+(\gamma+1)v\psi^y&-\gamma\psi^y&(\gamma-1)\rho\psi^y\\0&0&1&0\\0&L&0&n\end{pmatrix},\\[0.5ex]
&\bm E^x(\bm U)=(m,n,K^x,L)^\top,\quad\bm E^y(\bm U)=(m,n,K^y,L)^\top,
\end{aligned}
$$
where
$$
\begin{aligned}
&K^x:=\rho u^2+p+\int\limits_{\hat x}^x{\rho\phi_x}\,{\rm d}\xi,\quad K^y:=\rho v^2+p+\int\limits_{\hat y}^y{\rho\phi_y}\,{\rm d}\eta,\quad
L:=\frac{{\cal E}+p}{\rho},\\
&\psi^x:=\frac{2mn}{(\gamma-1)(2\rho^2L+n^2)-(\gamma+1)m^2},\quad\psi^y:=\frac{2mn}{(\gamma-1)(2\rho^2L+m^2)-(\gamma+1)n^2}.
\end{aligned}
$$

In order to apply the fifth-order Ai-WENO-Z interpolation to the equilibrium variables, we first compute
\begin{equation*}
K_{j,k}^x=\rho_{j,k}u_{j,k}^2+p_{j,k}+{\cal I}_{j,k}^x,\quad K_{j,k}^y=\rho_{j,k}v_{j,k}^2+p_{j,k}+{\cal I}_{j,k}^y,\quad
L_{j,k}=\frac{{\cal E}_{j,k}+p_{j,k}}{\rho_{j,k}},
\end{equation*}
where $u_{j,k}=m_{j,k}/\rho_{j,k}$, $v_{j,k}=n_{j,k}/\rho_{j,k}$, and ${\cal I}_{j,k}^x$ and ${\cal I}_{j,k}^y$ are fifth-order
approximations of the integral $\int_{x_{-\frac52}}^{x_j}\rho\phi_x\,{\rm d}x$ and $\int_{y_{-\frac52}}^{y_k}\rho\phi_y\,{\rm d}y$,
respectively. The values ${\cal I}_{j,k}^x$ are computed as follows. For each $k=1,\dots,N_y$, we evaluate
$f_{j,k}:=\rho_{j,k}\phi_x(x_j,y_k)$ for $j=1,\ldots,N_x$, and extend these values to $j=-5,\ldots,0$ and $j=N+1,\ldots,N+5$ using the
prescribed boundary conditions implemented within the ghost cell framework. We then compute ${\cal I}_{j,k}^x$, $j=-2,\dots,N_x+3$ using
\eref{3.5a}--\eref{3.5b}. The values of ${\cal I}_{j,k}^y$, $j=1,\ldots,N_x$, $k=-2,\dots,N_y+3$ can be computed similarly.

We then take the local parts of $\mE^x(\mU)$ and $\mE^y(\mU)$, which are $W^x(\mU)=(m,n,\rho u^2+p,L)$ and $W^y(\mU)=(m,n,\rho v^2+p,L)$,
and evaluate their Jacobians
\begin{align*}
&\frac{\partial W^x}{\partial U}(\bm U)=\begin{pmatrix}0&1&0&0\\0&0&1&0\\[0.5ex]
(1-\gamma)\phi+\dfrac{\gamma-3}{2}u^2+\dfrac{\gamma-1}{2}v^2 &(3-\gamma)u&(1-\gamma)v&\gamma-1\\[1.5ex]
(\gamma-1)\dfrac{u^2+v^2}{\rho}-\dfrac{\gamma{\cal E}}{\rho^2}&(1-\gamma)\dfrac{u}{\rho}&(1-\gamma)\frac{v}{\rho}&\dfrac{\gamma}{\rho}
\end{pmatrix},\\[1.5ex]
&\frac{\partial W^y}{\partial U}(\bm U)=\begin{pmatrix}0&1&0&0\\0&0&1&0\\[0.5ex]
(1-\gamma)\phi+\dfrac{\gamma-1}{2}u^2+\dfrac{\gamma-3}{2}v^2&(1-\gamma)u&(3-\gamma)v&\gamma-1\\[1.5ex]
(\gamma-1)\dfrac{u^2+v^2}{\rho}-\dfrac{\gamma{\cal E}}{\rho^2}&(1-\gamma)\dfrac{u}{\rho}&(1-\gamma)\frac{v}{\rho}&\dfrac{\gamma}{\rho}
\end{pmatrix}.
\end{align*}

\noindent
Notice that computing the matrices $C^x_{j,k}=\frac{\partial W^x}{\partial U}(\bm U_{j,k})M^x(\bm U_{j,k})$ and
$C^y_{j,k}=\frac{\partial W^y}{\partial U}(\bm U_{j,k})M^y(\bm U_{j,k})$ analytically is quite cumbersome and their evaluation will be
computationally expensive. Moreover, their analytical eigenstructures are unavailable. We therefore compute $C^x_{j,k}$ and $C^y_{j,k}$
together with the matrices $Q^x_{j,k}$, $(Q^x_{j,k})^{-1}$, $Q^y_{j,k}$, and $(Q^y_{j,k})^{-1}$ numerically. We then implement the LCD of
the equilibrium variables followed by the fifth-order Ai-WENO-Z interpolation in the $x$- and $y$-directions separately. In the
$x$-direction, this gives $\bm\Gamma^\pm_{j\mp\hf,k}$ and then $m^\pm_{j\mp\hf,k}$, $n^\pm_{j\mp\hf,k}$, $(K^x)^\pm_{j\mp\hf,k}$, and
$L^\pm_{j\mp\hf,k}$. After that, we solve the quadratic equations
\begin{equation}
\begin{aligned}
(\gamma-1)\big(L^\pm_{\jph,k}-\phi(x_\jph,y)\big)\big(\rho^\pm_{\jph,k}\big)^2-\gamma\big((K^x)^\pm_{\jph,k}-I^x_{\jph,k}\big)
\rho^\pm_{\jph,k}&\\
+\frac{(\gamma+1)\big(m^\pm_{\jph,k}\big)^2}{2}-\frac{(\gamma-1)\big(n^\pm_{\jph,k}\big)^2}{2}&=0
\end{aligned}
\label{3.14}
\end{equation}
using the method described in \cite[Appendix B.2]{Kurganov25Na} to obtain $\rho^\pm_{j\mp\hf,k}$ and then ${\cal E}^\pm_{j\mp\hf,k}$. In
\eref{3.14}, the integrals $I^x_{\jph,k}$ are evaluated recursively using \eref{3.5c}. The point values $m^\pm_{j,k\mp\hf}$,
$n^\pm_{j,k\mp\hf}$, $(K^y)^\pm_{j,k\mp\hf}$, $L^\pm_{j,k\mp\hf}$, $\rho^\pm_{j,k\mp\hf}$, and ${\cal E}^\pm_{j,k\mp\hf}$ can be computed
using a similar interpolation procedure carried out in the $y$-direction. Finally, we use the numerical fluxes described in
\cite[(2.7)]{Kurganov25Na} for $\bm{{\cal K}}_{\jph,k}$ and $\bm{{\cal L}}_{j,\kph}$.

\subsection{Application to the 2-D Ripa System}\label{sec3.8}
In this section, we consider the 2-D Ripa system, which was introduced in \cite{Rip93,Rip95} and reads as \eref{1.1} with
\begin{equation*}
\begin{aligned}
&\bm U=(h,q^x,q^y,h\theta,Z)^\top\!,~\bm F(\bm U)=\Big(q^x,q^xu+P,\frac{q^xq^y}{h},q^x\theta,0\Big)^\top\!,~
\bm G(\bm U)=\Big(q^y,\frac{q^xq^y}{h},q^yv+P,q^y\theta, 0\Big)^\top\!,\\
&B^x(\bm{U})=\begin{pmatrix}0&0&0&0&0\\0&0&0&0&-h\theta\\0&0&0&0&0\\0&0&0&0&0\\0&0&0&0&0\end{pmatrix},\quad
B^y(\bm{U})=\begin{pmatrix}0&0&0&0&0\\0&0&0&0&0\\0&0&0&0&-h\theta\\0&0&0&0&0\\0&0&0&0&0\end{pmatrix},\quad\bm S(\bm U)\equiv\bm0,
\end{aligned}
\end{equation*}
where $h$ is the water depth, $u$ and $v$ are the velocities in $x$- and $y$-directions, respectively, $q^x=hu$, $q^y=hv$, $Z(x,y)$ is the
bottom topography, $\theta$ is the potential temperature, and $P=\hf h^2\theta$ is the pressure.

We restrict our consideration to the quasi 1-D moving-water equilibria: the $x$-directional,
\begin{equation}
\begin{aligned}
&q^x_x=q^y=\theta_x=({\cal E}^x)_x\equiv0,\quad{\cal E}^x=\frac{u^2}{2}+\theta(h+Z)+{\cal I}^x,\\
&{\cal I}^x=-\int\limits_{\widehat x}^x\left[\sqrt{2P(\xi,y,t)}\big(\sqrt{\theta(\xi,y,t)}\big)_\xi+Z(\xi,y)\theta_\xi(\xi,y,t)\right]
{\rm d}\xi,
\end{aligned}
\label{4.7}
\end{equation}
and the $y$-directional,
\begin{equation}
\begin{aligned}
&q^x=q^y_y=\theta_y=({\cal E}^y)_y\equiv0,\quad{\cal E}^y=\frac{v^2}{2}+\theta(h+Z)+{\cal I}^y,\\
&{\cal I}^y=-\int\limits_{\widehat y}^y\left[\sqrt{2P(x,\eta,t)}\big(\sqrt{\theta(x,\eta,t)}\big)_\eta+Z(x,\eta)\theta_\eta(x,\eta,t)
\right]{\rm d}\eta,
\end{aligned}
\label{4.8}
\end{equation}
ones. The steady states \eref{4.7} and \eref{4.8} satisfy $M^x(\bm U)\bm E^x(\bm U)_x=\bm0$ and $M^y(\bm U)\bm E^y(\bm U)_y=\bm0$,
respectively, with
$$
\begin{aligned}
&M^x(\bm U)=\begin{pmatrix}1&0&0&0&0\\u&h&0&0&0\\v&0&q^x&0&0\\\theta&0&0&q^x&0\\0&0&0&0&0\end{pmatrix},\quad
M^y(\bm U)=\begin{pmatrix}1&0&0&0&0\\u&0&q^y&0&0\\v&h&0&0&0\\\theta&0&0&q^y&0\\0&0&0&0&0\end{pmatrix},\quad
\begin{matrix}\bm E^x(\bm U)=\Big(q^x,{\cal E}^x,v,\theta,Z\Big)^\top\!,\\
\bm E^y(\bm U)=\Big(q^y,{\cal E}^y,u,\theta,Z\Big)^\top\!.
\end{matrix}
\end{aligned}
$$
For the 2-D Ripa system, the matrices $C^x(\bm U)$ and $C^y(\bm U)$ can be explicitly computed:
\begin{equation}
C^x(\bm U)=\begin{pmatrix}u&h&0&0&0\\[0.5ex]\theta&u&0&\dfrac{q^x}{2}&0\\[1.5ex]0&0&u&0&0\\0&0&0&u&0\\0&0&0&0&0\end{pmatrix},\quad
C^y(\bm U)=\begin{pmatrix}v&h&0&0&0\\\theta&v&0&\dfrac{q^y}{2}&0\\0&0&v&0&0\\0&0&0&v&0\\0&0&0&0&0\end{pmatrix},
\label{4.14}
\end{equation}
and they have complete eigensystems so that the corresponding matrices $Q^x(\bm U)$ and $Q^y(\bm U)$ along with their inverses are given by
$$
\begin{aligned}
&Q^x(\bm U)=\begin{pmatrix}0&0&-\dfrac{q^x}{2\theta}&\dfrac{c}{\theta}&-\dfrac{c}{\theta}\\[1.5ex]0&0&0&1&1\\0&1&0&0&0\\0&0&1&0&0\\1&0&0&0&0
\end{pmatrix},\quad[Q^x(\bm U)]^{-1}=\begin{pmatrix}0&0&0&0&1\\0&0&1&0&0\\0&0&0&1&0\\[0.5ex]
\dfrac{\theta}{2c}&\dfrac{1}{2}&0&-\dfrac{q^x}{4c}&0\\[2.0ex]-\dfrac{\theta}{2c}&\dfrac{1}{2}&0&-\dfrac{q^x}{4c}&0\end{pmatrix},\\[1.5ex]
&Q^y(\bm U)=\begin{pmatrix}0&0&-\dfrac{q^y}{2\theta}&\dfrac{c}{\theta}&-\dfrac{c}{\theta}\\[1.5ex]0&0&0&1&1\\0&1&0&0&0\\0&0&1&0&0\\1&0&0&0&0
\end{pmatrix},\quad[Q^y(\bm U)]^{-1}=\begin{pmatrix}0&0&0&0&1\\0&0&1&0&0\\0&0&0&1&0\\[0.5ex]
\dfrac{\theta}{2c}&\dfrac{1}{2}&0&-\dfrac{q^y}{4c}&0\\[2.0ex]-\dfrac{\theta}{2c}&\dfrac{1}{2}&0&-\dfrac{q^y}{4c}&0\end{pmatrix}.
\end{aligned}
$$
In order to apply the fifth-order Ai-WENO-Z interpolation to the equilibrium variables, we first compute
\begin{equation*}
{\cal E}^x_{j,k}=\frac{u_{j,k}^2}{2}+\theta_{j,k}(h_{j,k}+Z_{j,k})+{\cal I}^x_{j,k},\quad
{\cal E}^y_{j,k}=\frac{v_{j,k}^2}{2}+\theta_{j,k}(h_{j,k}+Z_{j,k})+{\cal I}^y_{j,k},
\end{equation*}
where ${\cal I}^x_{j,k}$ and ${\cal I}^y_{j,k}$ are computed using the path-conservative integration as in \cite{CKL23}. We then implement
the LCD of the equilibrium variables followed by the fifth-order Ai-WENO-Z interpolation in the $x$- and $y$-directions separately. In the
$x$-direction, this gives $\bm\Gamma^\pm_{j\mp\hf,k}$ and then $({\cal E}^x)_{\jph,k}^\pm$, $(q^x)^\pm_{j\mp\hf,k}$, $v^\pm_{j\mp\hf,k}$,
$\theta^\pm_{j\mp\hf,k}$, and $Z^\pm_{j\mp\hf,k}$. After that, we solve the quadratic equations
\begin{equation}
\hf\left(\frac{(q^x)_{\jph,k}^\pm}{h_{\jph,k}^\pm}\right)^2+\theta_{\jph,k}^\pm\big(h_{\jph,k}^\pm+Z_{\jph,k}^\pm\big)+
({\cal I}^x)_{\jph,k}^\pm-({\cal E}^x)_{\jph,k}^\pm=0
\label{3.24f}
\end{equation}
to obtain $h^\pm_{j\mp\hf,k}$. In \eref{3.24f}, the values $({\cal I}^x)_{\jph,k}^\pm$ are evaluated recursively also using the
path-conservative integration. A similar interpolation procedure is carried out in the $y$-direction. Finally, we use the numerical fluxes
for $\bm{{\cal K}}_{\jph,k}$ and $\bm{{\cal L}}_{j,\kph}$. We refer the reader to \cite[\S5]{QGKWW} for details.

\section{Numerical Examples}\label{sec4}
In this section, we test the proposed fifth-order WB A-WENO scheme based on the LCD of equilibrium variables on several numerical examples
for the nozzle flow system, one- and two-layer shallow water equations, compressible Euler equations with gravitation, and the 2-D Ripa
system. For the sake of brevity, this scheme will be referred to as Scheme 1 and its performance will be compared with the following two
schemes:

\medskip
\noindent
$\bullet$ Scheme 2: The WB A-WENO scheme from \cite{CKX_24WB}, in which the Ai-WENO-Z interpolation is applied to the equilibrium variables
$\bm E$ without any LCD;

\smallskip
\noindent
$\bullet$ Scheme 3: The A-WENO scheme from \cite{Chu21}, which can only preserve the simplest ``lake-at-rest'' steady states, but uses the LCD applied to the conservative variables $\bm U$ to reduce the spurious oscillations.

\medskip
In all of the examples, we solve the ODE systems \eref{2.3} and \eref{2.8a} using the three-stage third-order SSP Runge--Kutta method  (see,
e.g., \cite{Gottlieb11,Gottlieb12}) with the time-step restricted by the CFL number $0.5$.

When the eigendecomposition of matrices $C$, $C^x$, or $C^y$ is required, we perform it using the FORTRAN LAPACK routine GEEV.

\subsection{Nozzle Flow System}
\subsubsection*{Example 1---Flow in Continuous Divergent Nozzle}
In the first example taken from \cite{KLX_21,CKX_24WB}, we consider the divergent nozzle described using the smooth cross-section
\begin{equation*}
\sigma(x)=0.976+0.748\tanh(0.8x-4).
\end{equation*}
We first take the steady states with $q_{\rm eq}(x)\equiv8$, ${\cal E}_{\rm eq}(x)\equiv21.9230562619897$ and compute the discrete values of
$\rho_{\rm eq}(x)$ by solving the corresponding nonlinear equations; see \cite{KLX_21,CKX_24WB}. We then obtain
$u_{\rm eq}(x)=q_{\rm eq}(x)/(\sigma(x)\rho_{\rm eq}(x))$.

Equipped with these steady states, we add a small perturbation to the density field and consider the initial data
\begin{equation*}
\rho(x,0)=\rho_{\rm eq}(x)+\begin{cases}10^{-2},&x\in[0.5,1.5],\\0,&\mbox{otherwise},\end{cases}\qquad
q(x,0)=\sigma(x)\rho(x,0)u_{\rm eq}(x),
\end{equation*}
which are prescribed in the computational domain $[0,10]$ subject to the homogeneous Neumann boundary conditions.

\begin{figure}[ht!]
\centerline{\includegraphics[trim=0.4cm 0.3cm 0.9cm 0.1cm, clip, width=5.7cm]{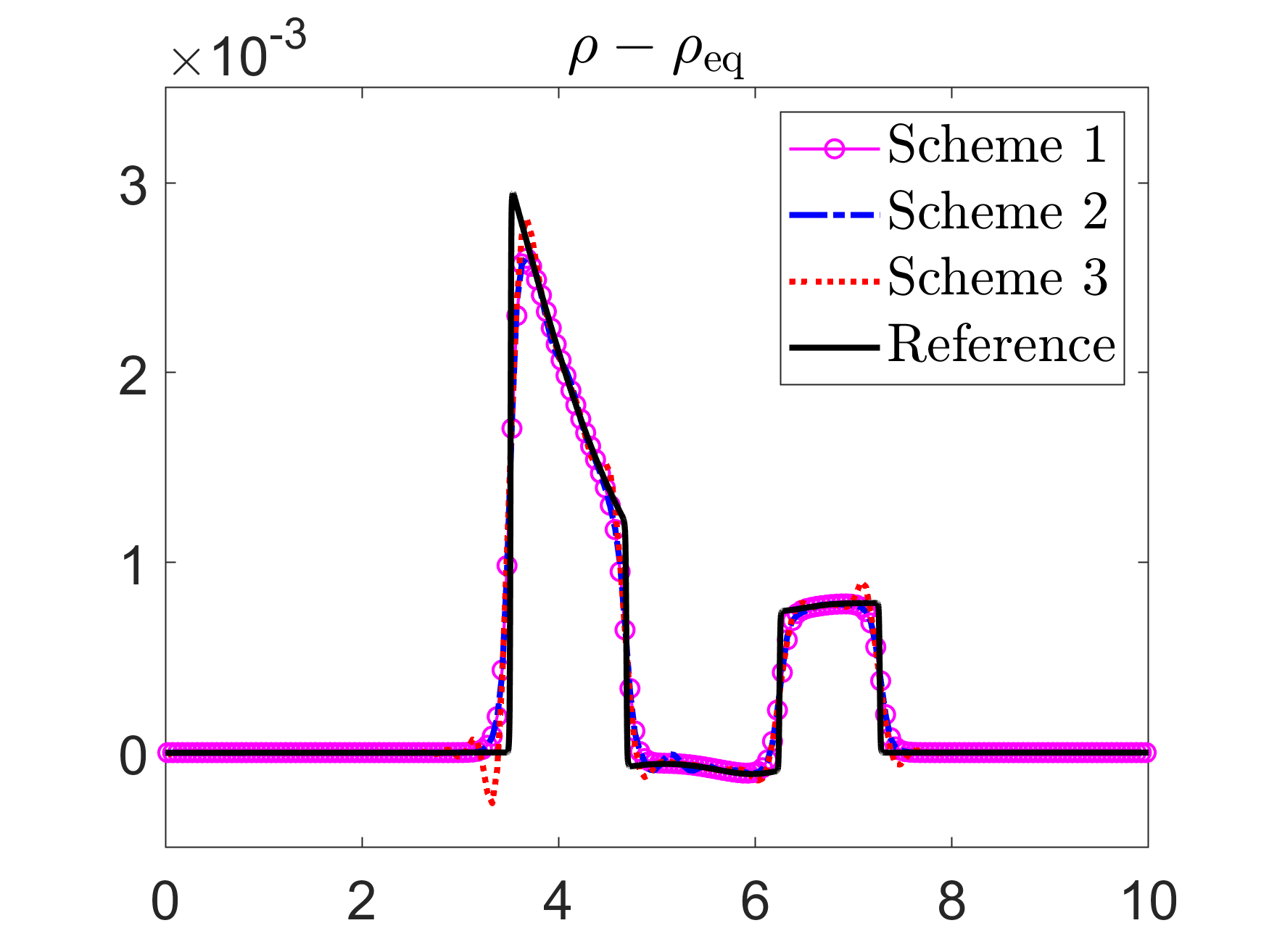}\hspace{0.3cm}
            \includegraphics[trim=0.4cm 0.3cm 0.9cm 0.1cm, clip, width=5.7cm]{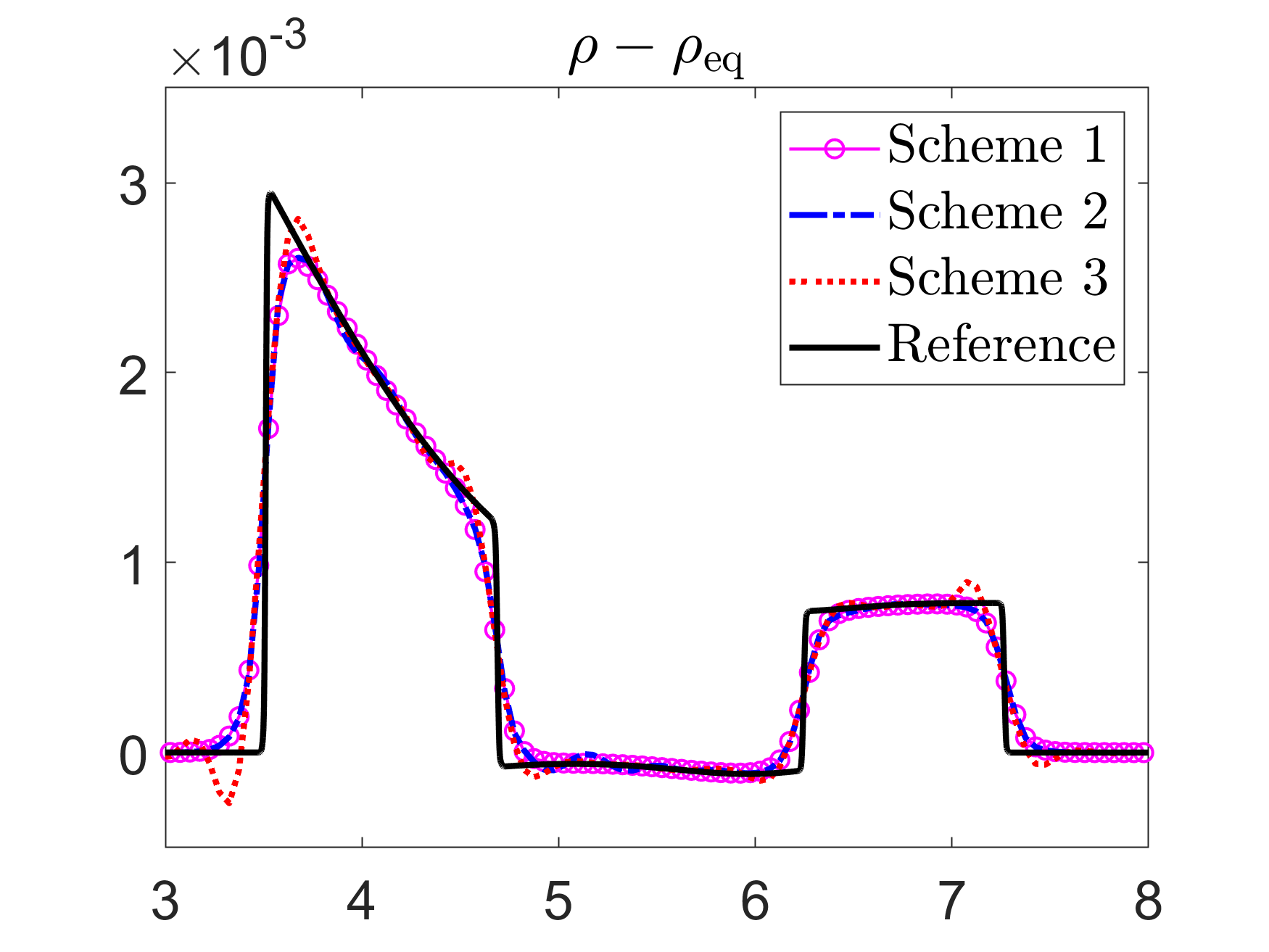}}
\vskip7pt
\centerline{\includegraphics[trim=0.4cm 0.3cm 0.9cm 0.1cm, clip, width=5.7cm]{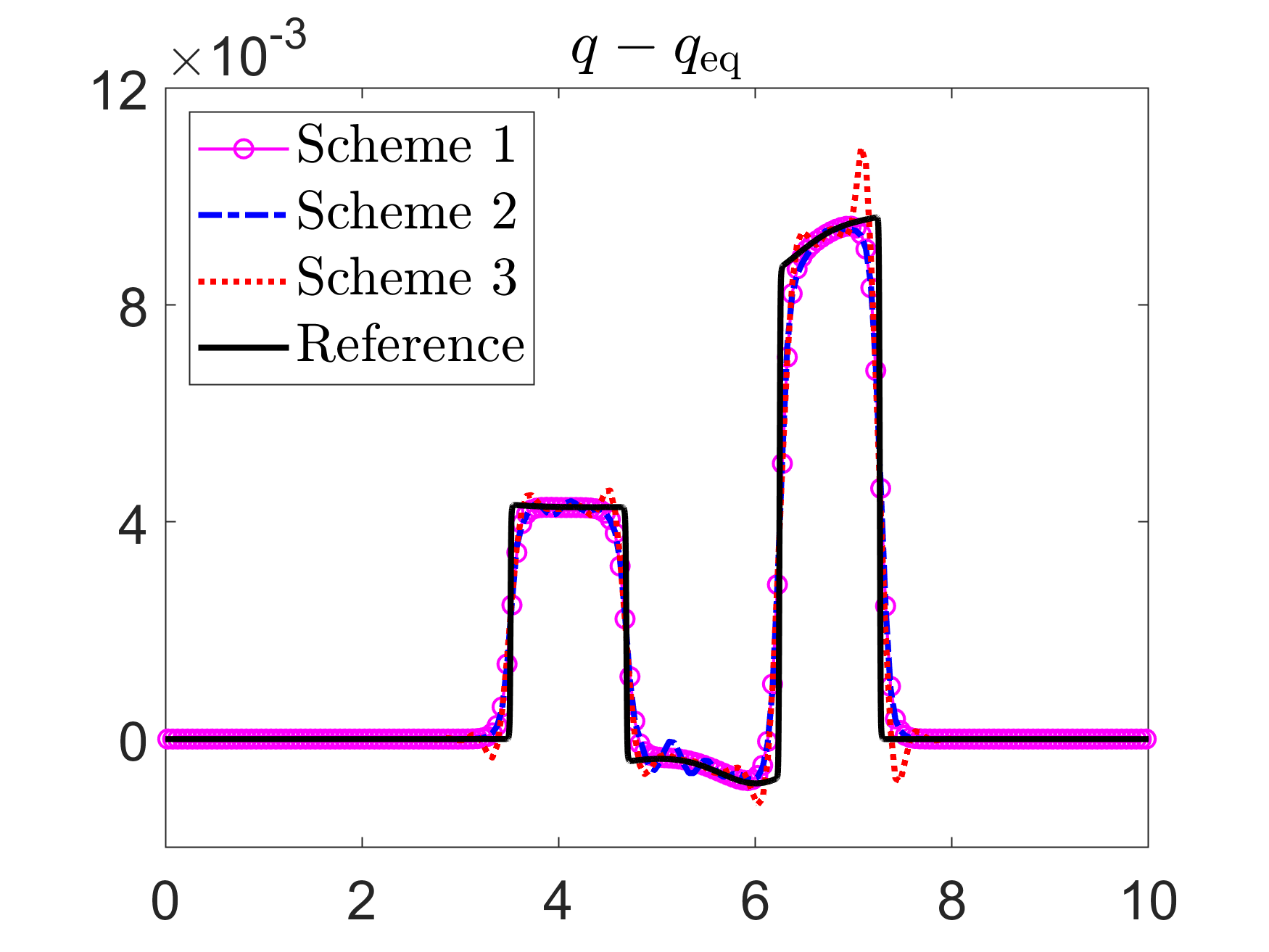}\hspace{0.3cm}
            \includegraphics[trim=0.4cm 0.3cm 0.9cm 0.1cm, clip, width=5.7cm]{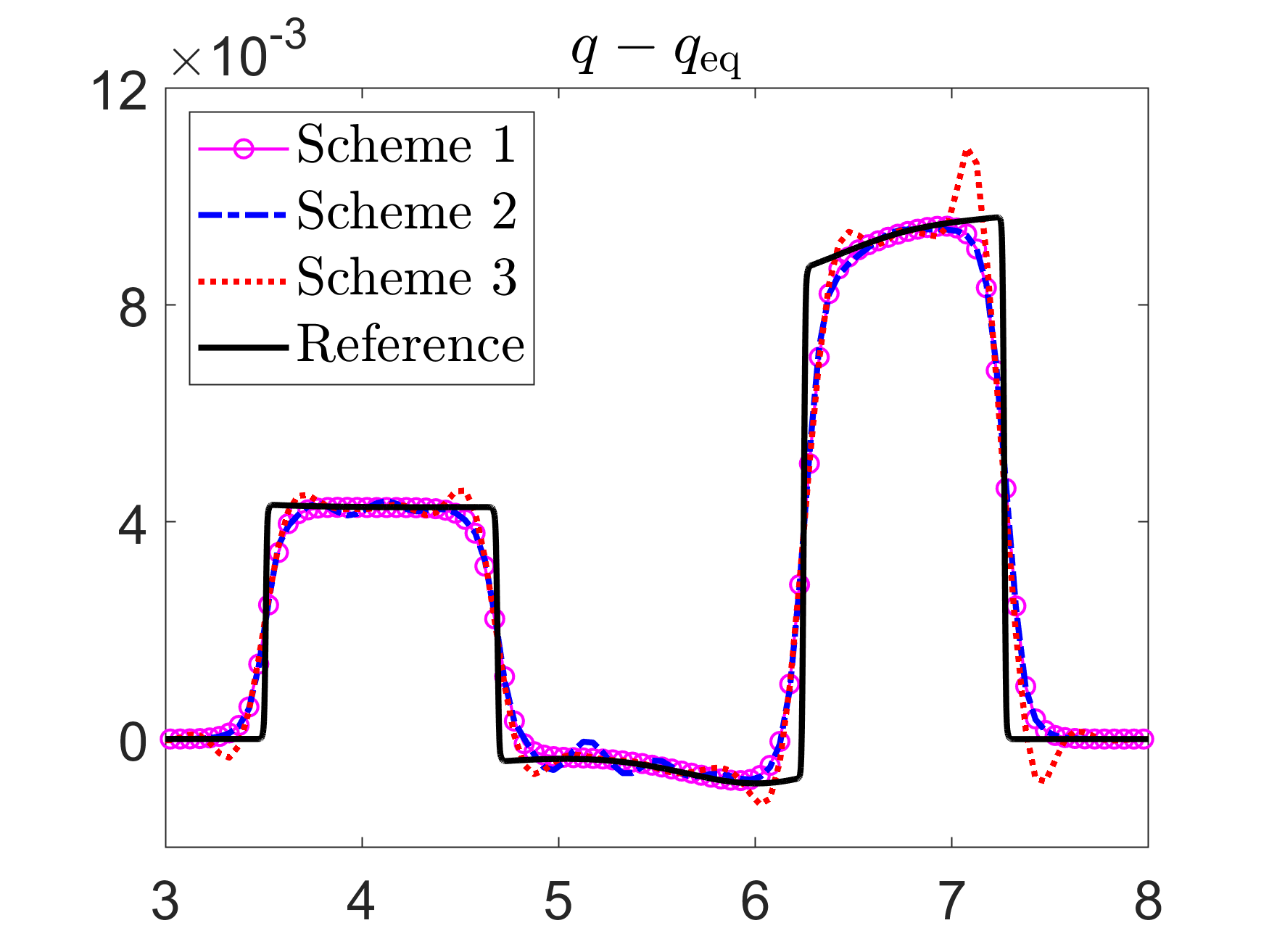}}
\caption{\sf Example 1: The differences $\rho(x,0.8)-\rho_{eq}(x)$ (top row) and  $q(x,0.8)-q_{eq}(x)$ (bottom row) computed by Schemes
1--3, and zoom at $x\in[3,8]$ (right column).\label{fig1}}
\end{figure}

We compute the numerical solutions until the final time $t=0.8$ by Schemes 1--3 on a uniform mesh with $\dx=1/20$ and plot the differences
$\rho(x,0.8)-\rho_{\rm eq}(x)$ and $q(x,0.8)-q_{\rm eq}(x)$ together with the reference solution computed by Scheme 1 on a much finer mesh
with $\dx=1/400$ in Figure \ref{fig1}. As one can see, Scheme 1 clearly outperforms Schemes 2 and 3 as there are no oscillations in the
results computed by Scheme 1. We also stress that in this example, Scheme 3, which is not WB, produces the largest oscillations as we
capture small perturbations of the steady state.

\subsubsection*{Example 2---Flow in Continuous Convergent Nozzle}
In the second example also taken from \cite{KLX_21,CKX_24WB}, we consider the convergent nozzle described using the smooth cross-section
\begin{equation*}
\sigma(x)=0.976-0.748\tanh(0.8x-4).
\end{equation*}
We first take the steady states with $q_{\rm eq}(x)\equiv8$, ${\cal E}_{\rm eq}(x)\equiv58.3367745090349$, compute the discrete values of
$\rho_{\rm eq}(x)$ by solving the corresponding nonlinear equations, and then obtain
$u_{\rm eq}(x)=q_{\rm eq}(x)/(\sigma(x)\rho_{\rm eq}(x))$; see \cite{KLX_21,CKX_24WB}. We then consider the initial data containing a
substantially larger perturbation than the one studied in Example 1:
\begin{equation*}
\rho(x,0)=\rho_{\rm eq}(x)+\begin{cases}0.3,&x\in[0.5,1.5],\\0,&\mbox{otherwise},\end{cases}\qquad
q(x,0)=\sigma(x)\rho(x,0)u_{\rm eq}(x).
\end{equation*}
As in Example 1, the computational domain is $[0,10]$ and the homogeneous Neumann boundary conditions are imposed.

\begin{figure}[ht!]
\centerline{\includegraphics[trim=0.4cm 0.3cm 0.9cm 0.1cm, clip, width=5.7cm]{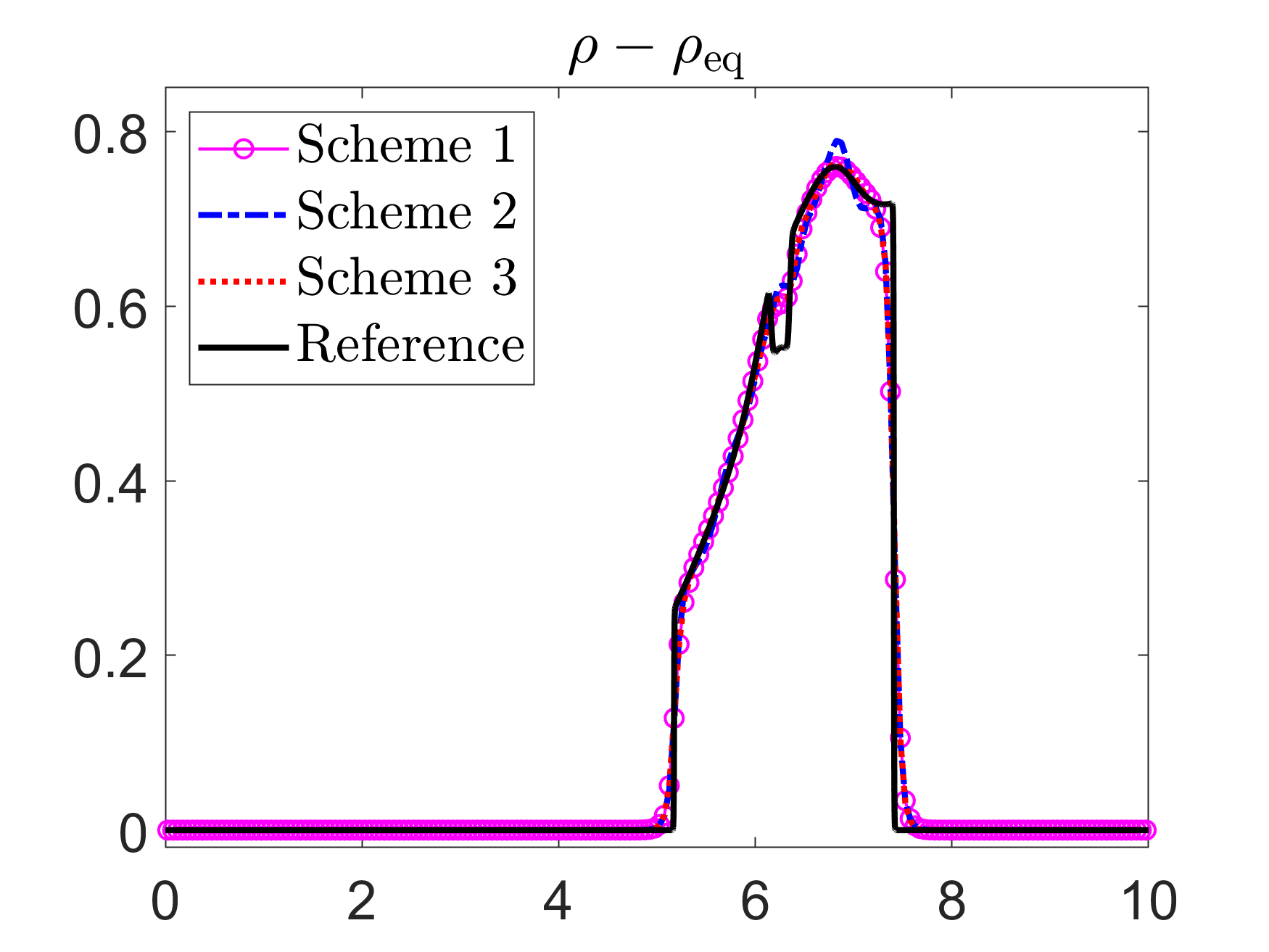}\hspace{0.3cm}
            \includegraphics[trim=0.4cm 0.3cm 0.9cm 0.1cm, clip, width=5.7cm]{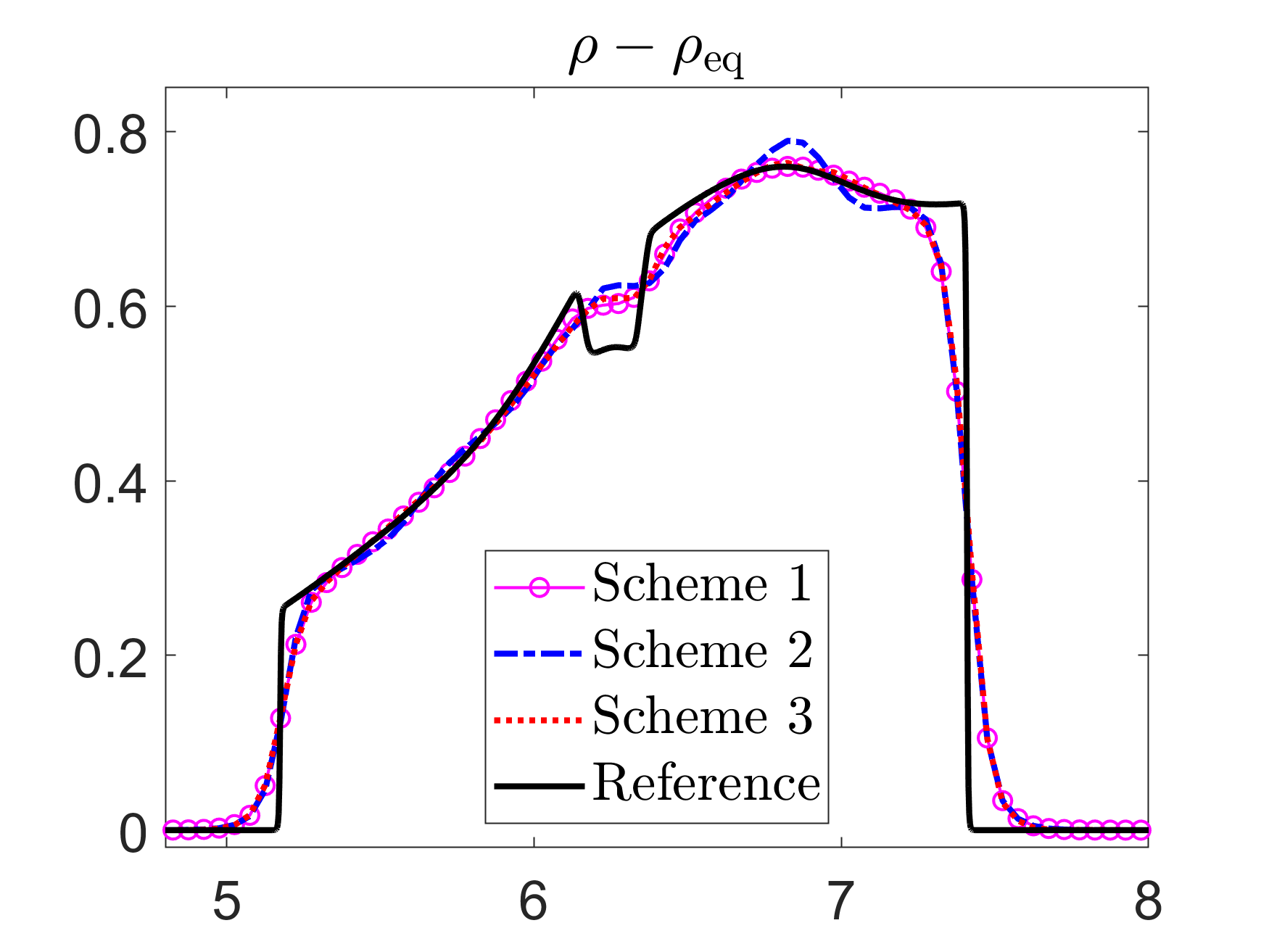}}
\vskip 7pt
\centerline{\includegraphics[trim=0.4cm 0.3cm 0.9cm 0.1cm, clip, width=5.7cm]{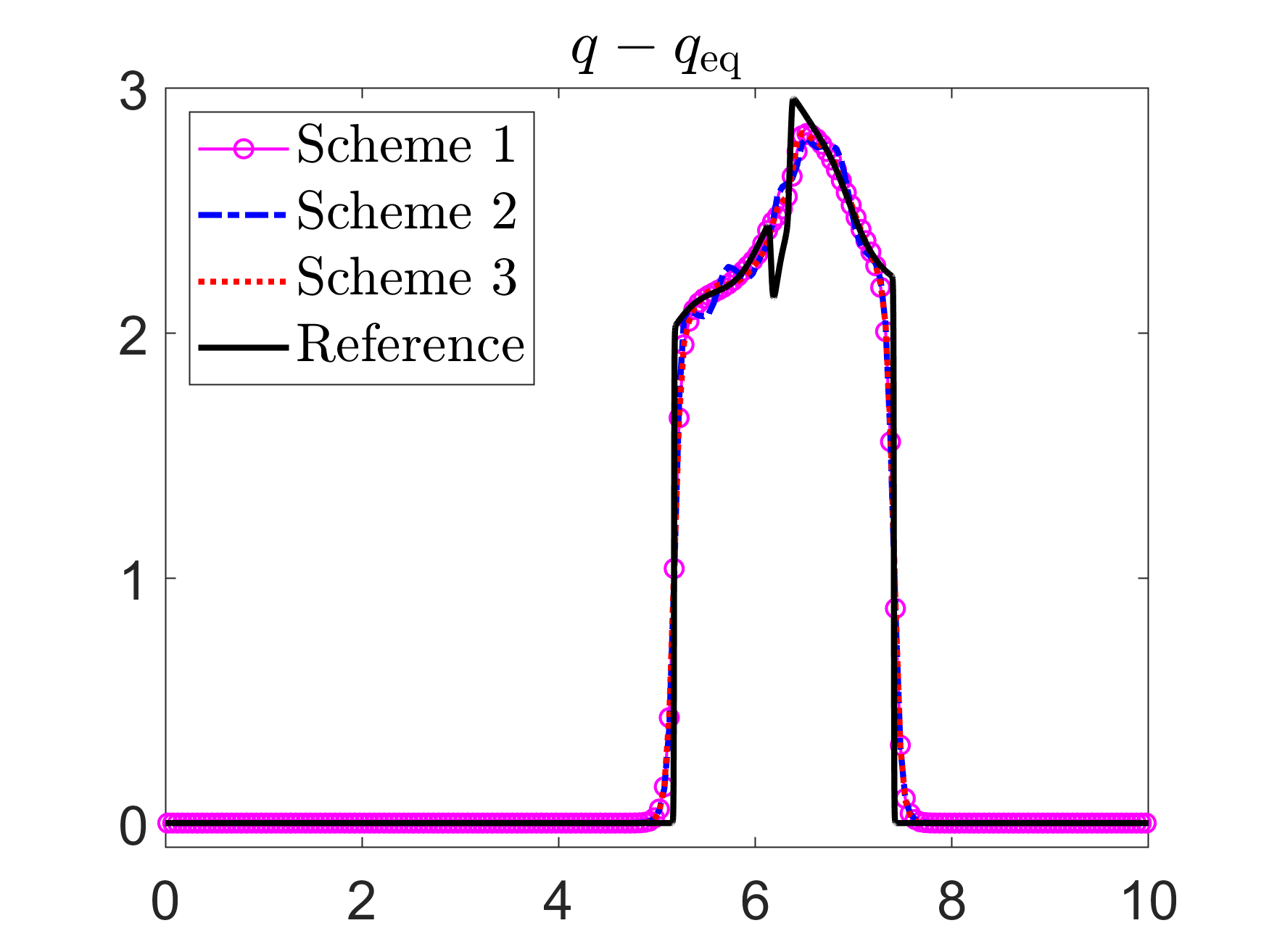}\hspace{0.3cm}
            \includegraphics[trim=0.4cm 0.3cm 0.9cm 0.1cm, clip, width=5.7cm]{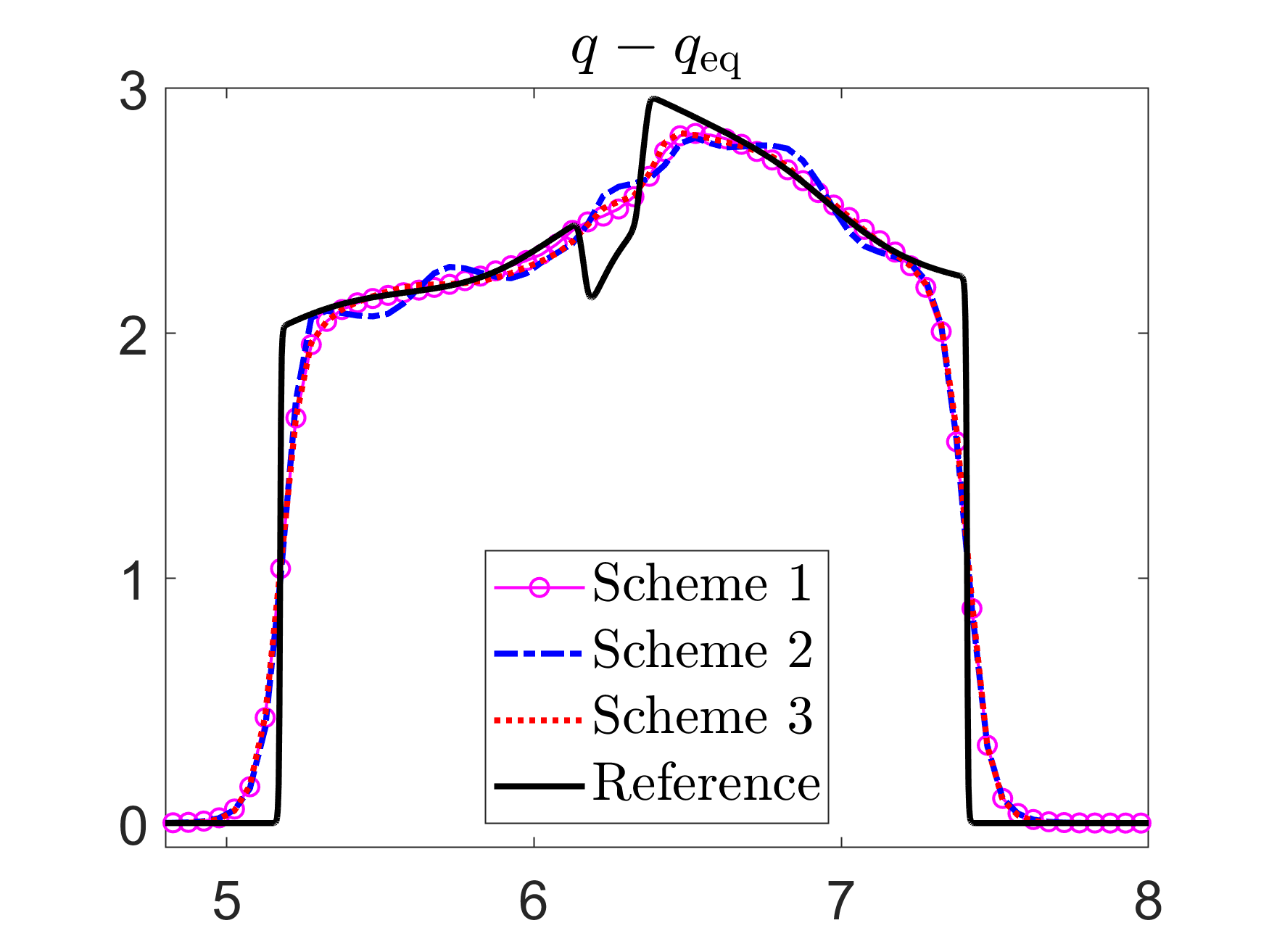}}
\caption{\sf Example 2: The differences $\rho(x,0.5)-\rho_{eq}(x)$ (top row) and  $q(x,0.5)-q_{eq}(x)$ (bottom row) computed by Schemes
1--3, and zoom at $x\in[4.8,8]$ (right column).\label{fig3}}
\end{figure}

We compute the numerical solutions until the final time $t=0.5$ by Schemes 1--3 on a uniform mesh with $\dx=1/20$ together with the
reference solution computed by Scheme 1 on a much finer mesh with $\dx=1/400$ and plot the differences $\rho(x,0.5)-\rho_{\rm eq}(x)$ and
$q(x,0.5)-q_{\rm eq}(x)$ in Figure \ref{fig3}. As one can see, Scheme 3 does not produce visible oscillations as the magnitude of the
perturbation is apparently larger than the size of the truncation errors. On the contrary, Scheme 2, which does not use any LCD, now
produces larger oscillations than in Example 1, where the size of perturbation was much smaller.

\subsection{Saint-Venant System with Manning Friction}
\subsubsection*{Example 3---Riemann Problem ($n=0.4$)}
In this example, we test the performance of Schemes 1--3 on a Riemann problem with the following initial data (prescribed in the
computational domain $[-0.1,0.3]$ subject to the homogeneous Neumann boundary conditions):
\begin{equation*}
h(x,0)=\left\{\begin{aligned}&1,&&x<0,\\&0.8,&&x>0,\end{aligned}\right.\qquad
u(x,0)=\left\{\begin{aligned}&2,&&x<0,\\&4,&&x>0,\end{aligned}\right.
\end{equation*}
and the bottom topography also containing a jump at $x=0$:
\begin{equation*}
Z(x)=\left\{\begin{aligned}&1,&&x<0,\\&1.9,&&x>0.\end{aligned}\right.
\end{equation*}

\begin{figure}[ht!]
\centerline{\includegraphics[trim=0.4cm 0.3cm 0.9cm 0.1cm, clip, width=5.7cm]{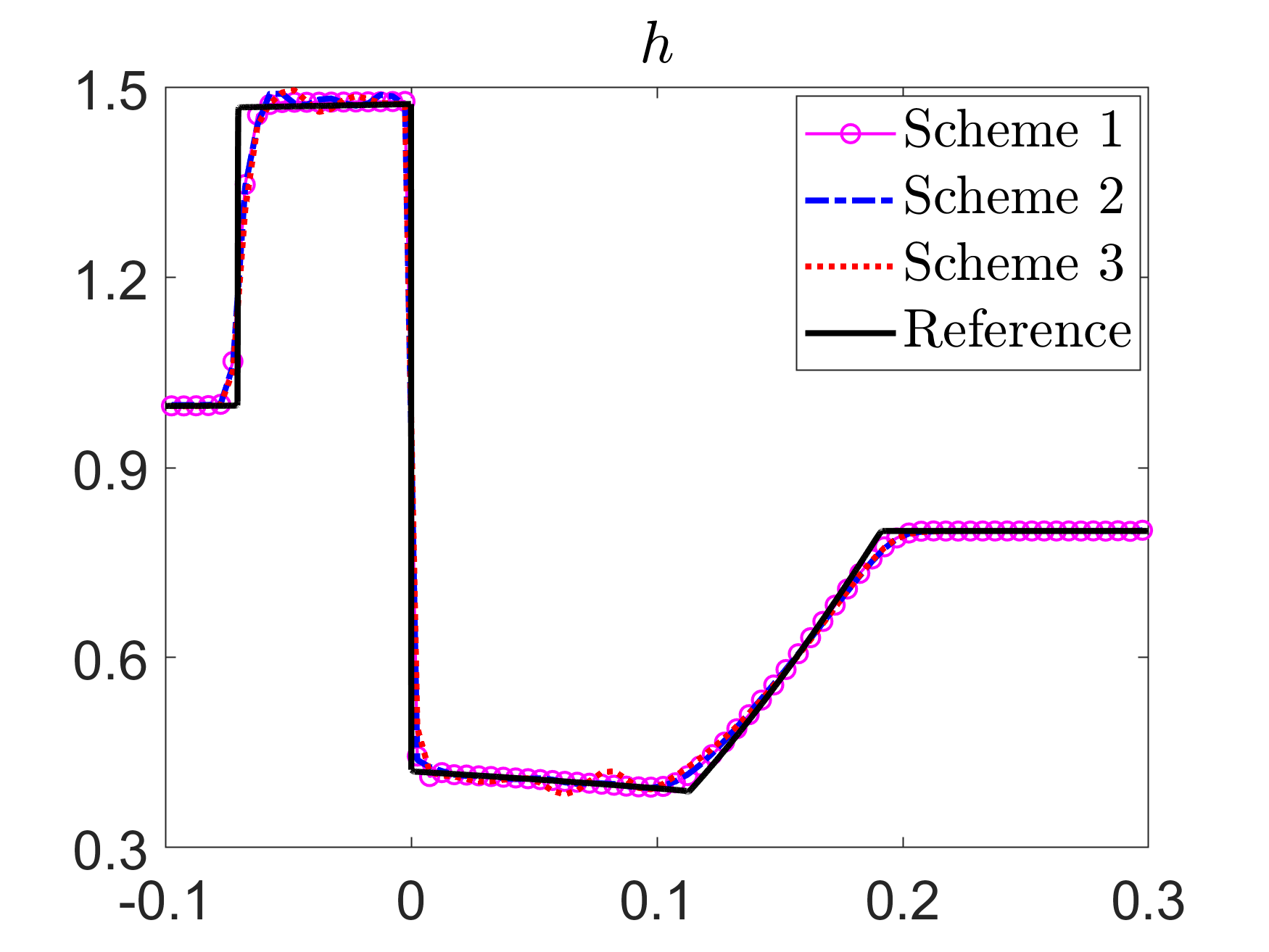}\hspace*{0.0cm}
	    \includegraphics[trim=0.4cm 0.3cm 0.9cm 0.1cm, clip, width=5.7cm]{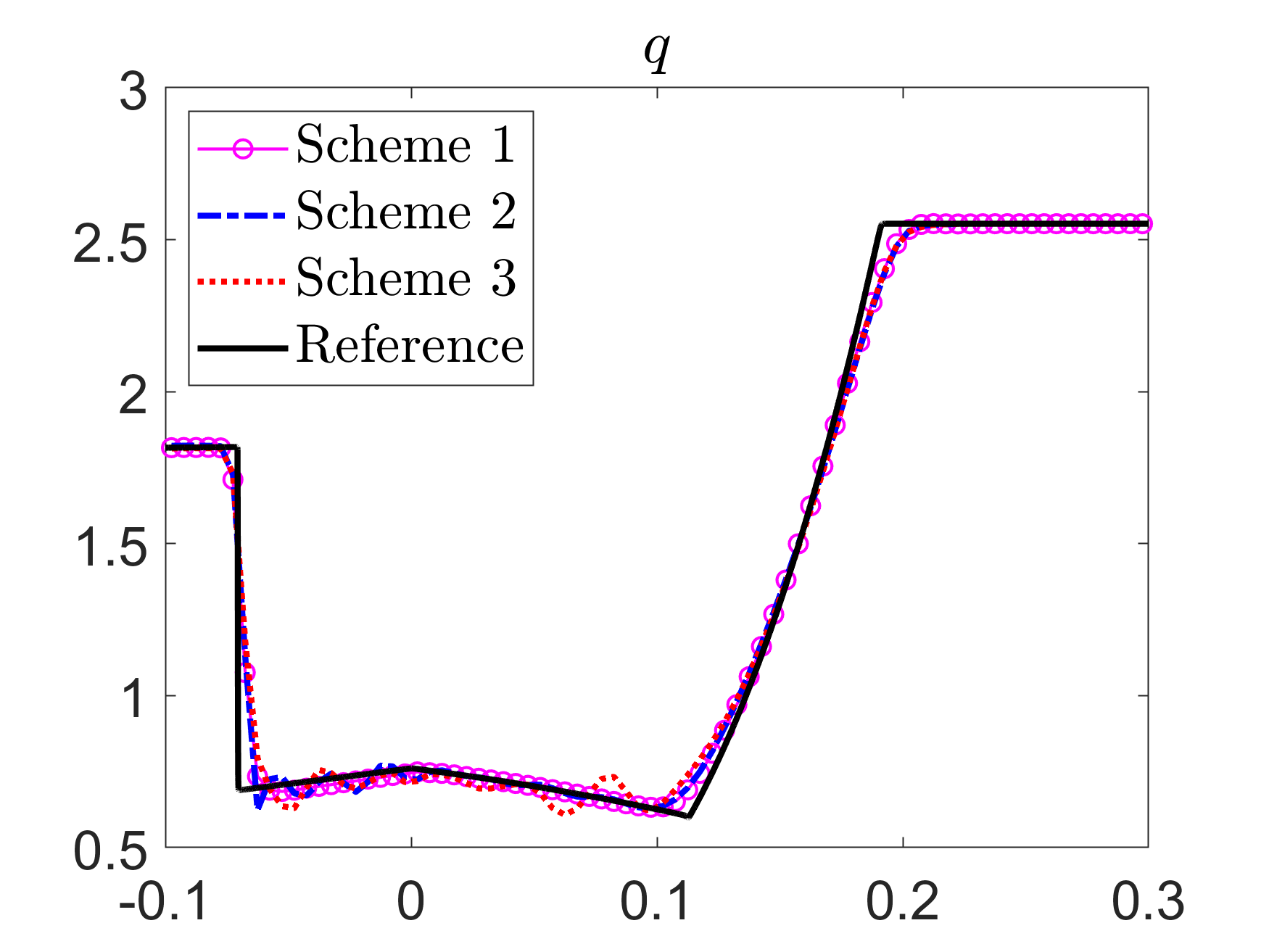}\hspace*{0.0cm}
            \includegraphics[trim=0.4cm 0.3cm 0.9cm 0.1cm, clip, width=5.7cm]{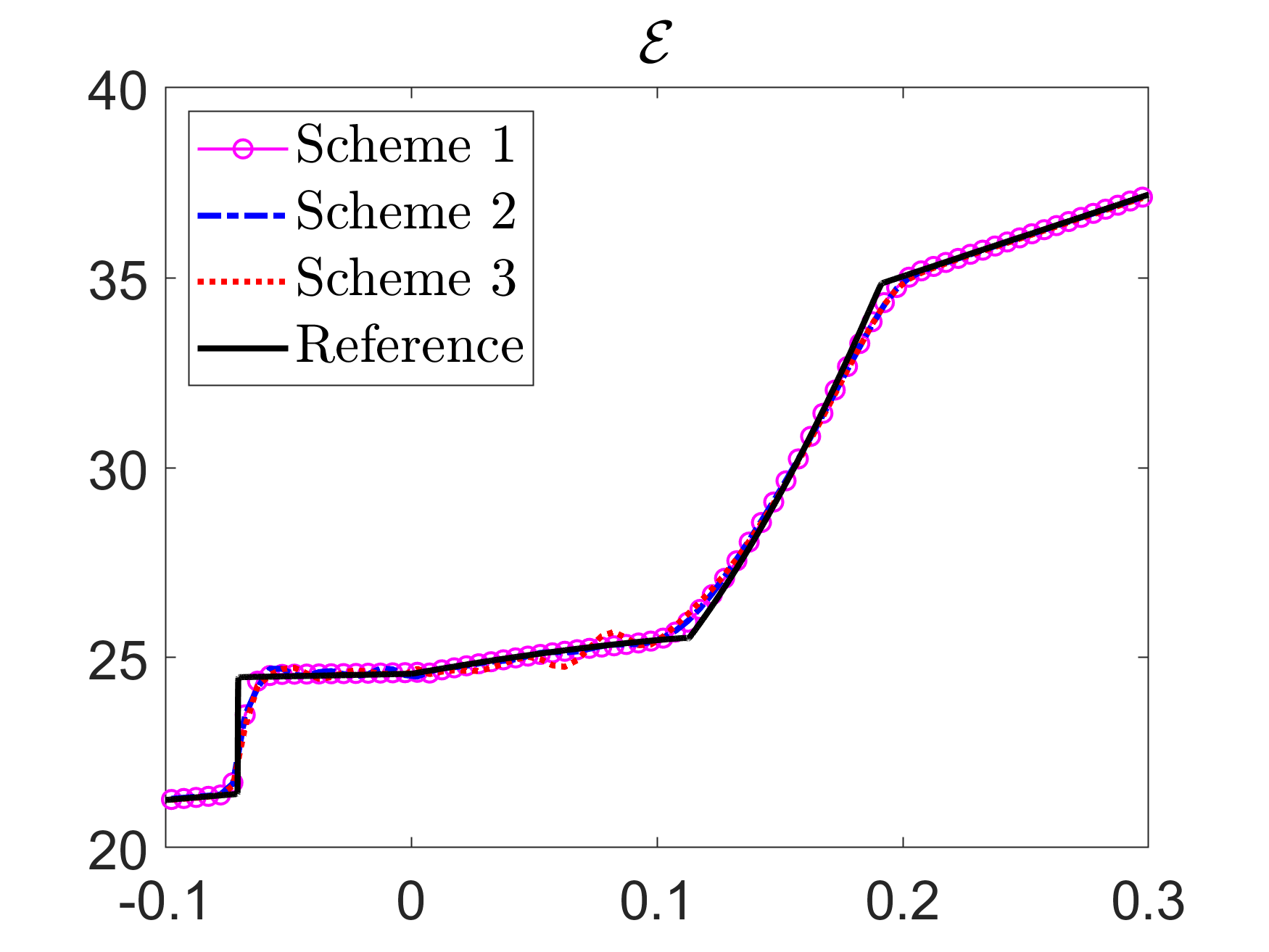}}
\caption{\sf Example 3: Water depth $h$, discharge $q$, and energy {$\cal E$} computed by Schemes 1--3 using $\dx=1/200$.\label{fig211}}
\end{figure}

\begin{figure}[ht!]
\centerline{\includegraphics[trim=0.4cm 0.3cm 0.9cm 0.1cm, clip, width=5.7cm]{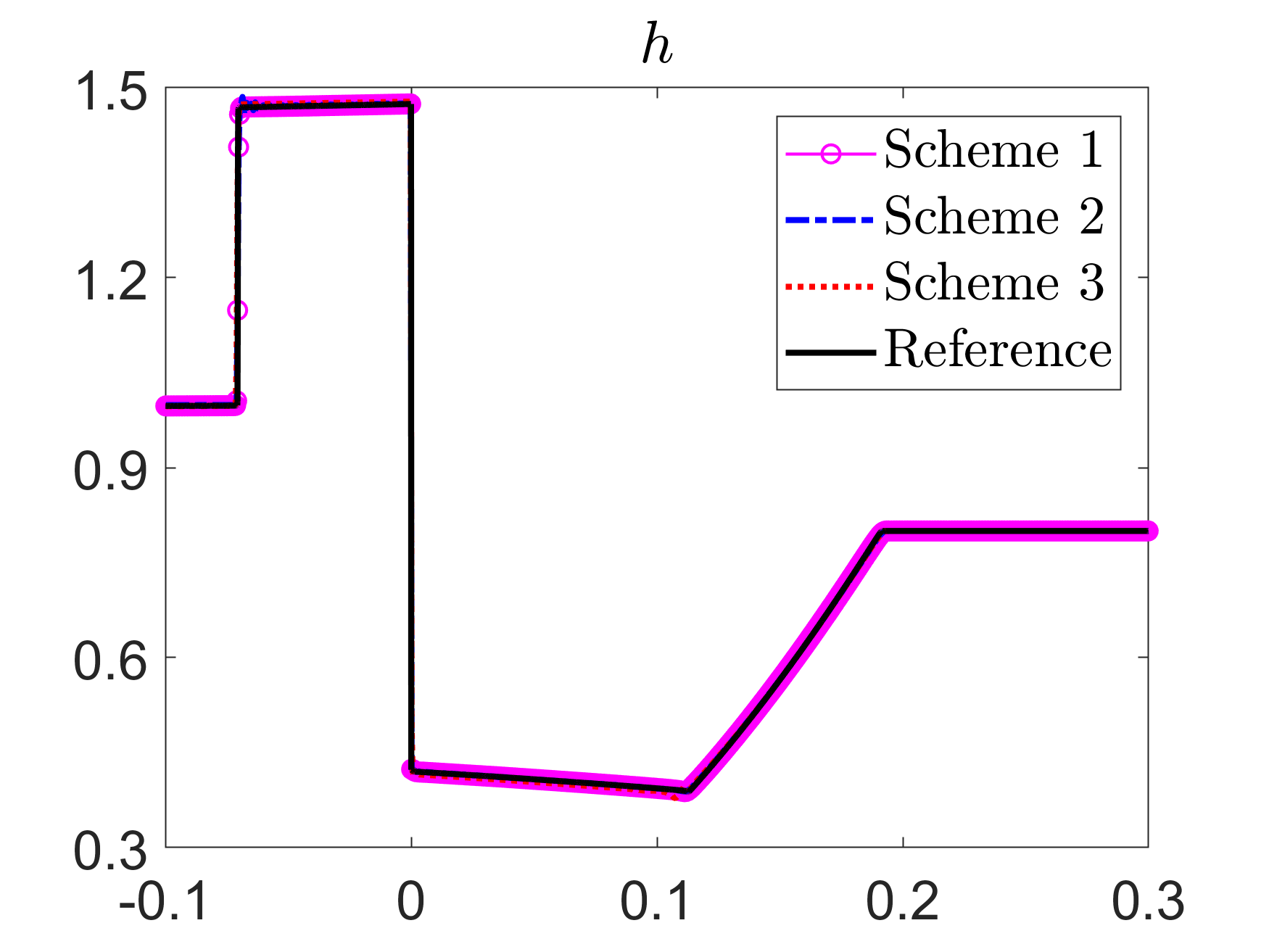}\hspace*{0.3cm}
	    \includegraphics[trim=0.4cm 0.3cm 0.9cm 0.1cm, clip, width=5.7cm]{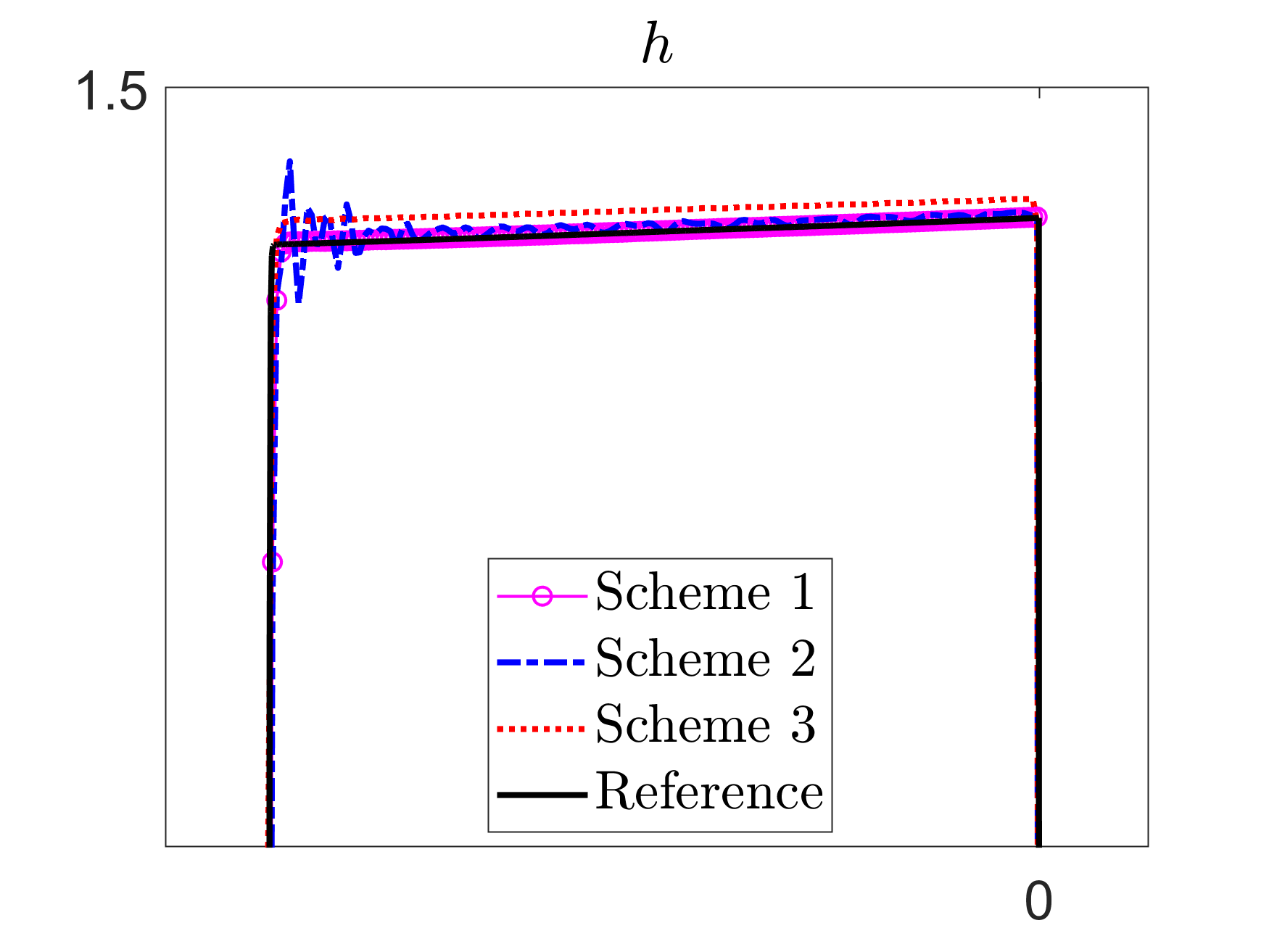}}
\vskip7pt
\centerline{\includegraphics[trim=0.4cm 0.3cm 0.9cm 0.1cm, clip, width=5.7cm]{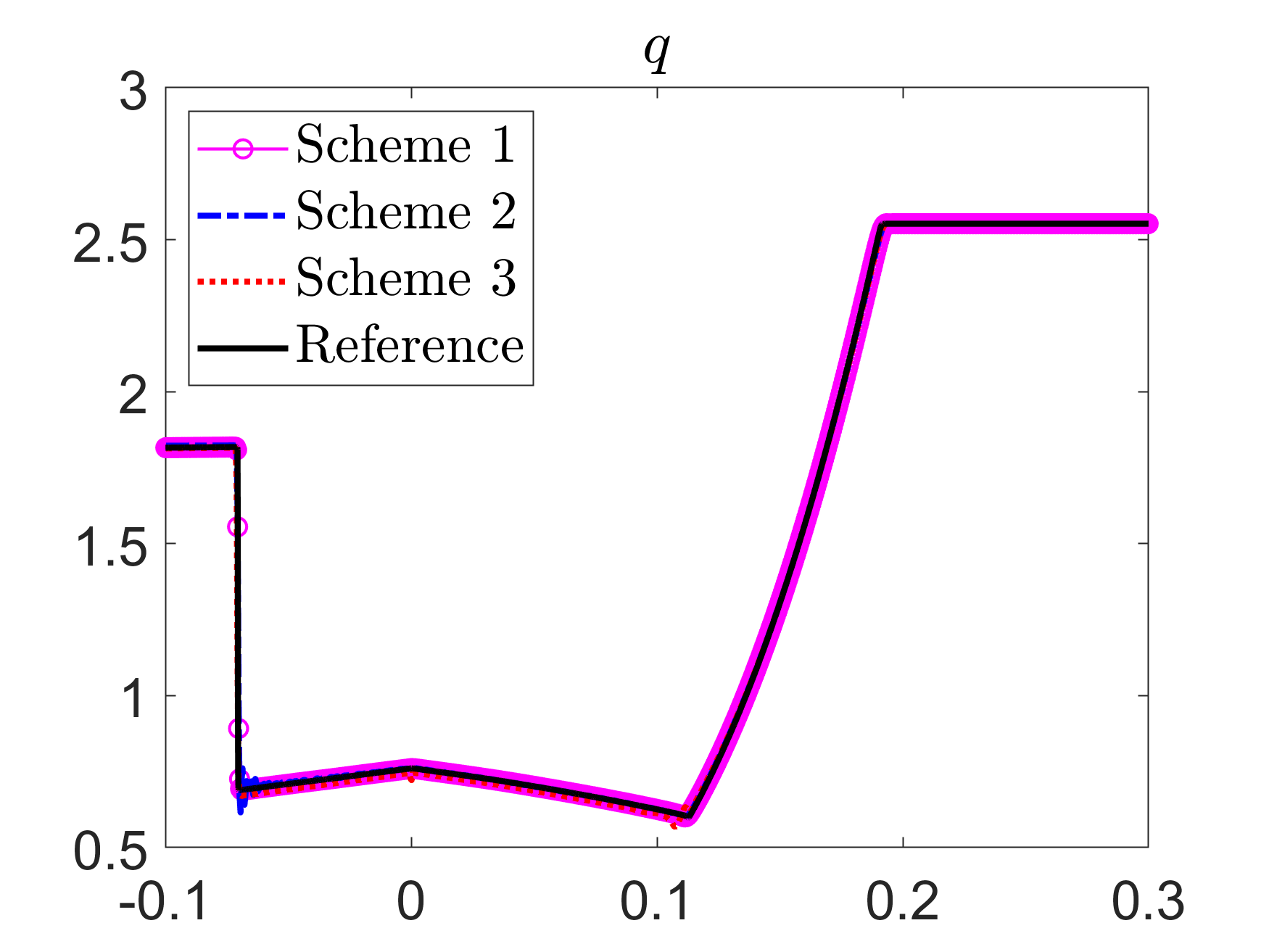}\hspace*{0.3cm}
            \includegraphics[trim=0.4cm 0.3cm 0.9cm 0.1cm, clip, width=5.7cm]{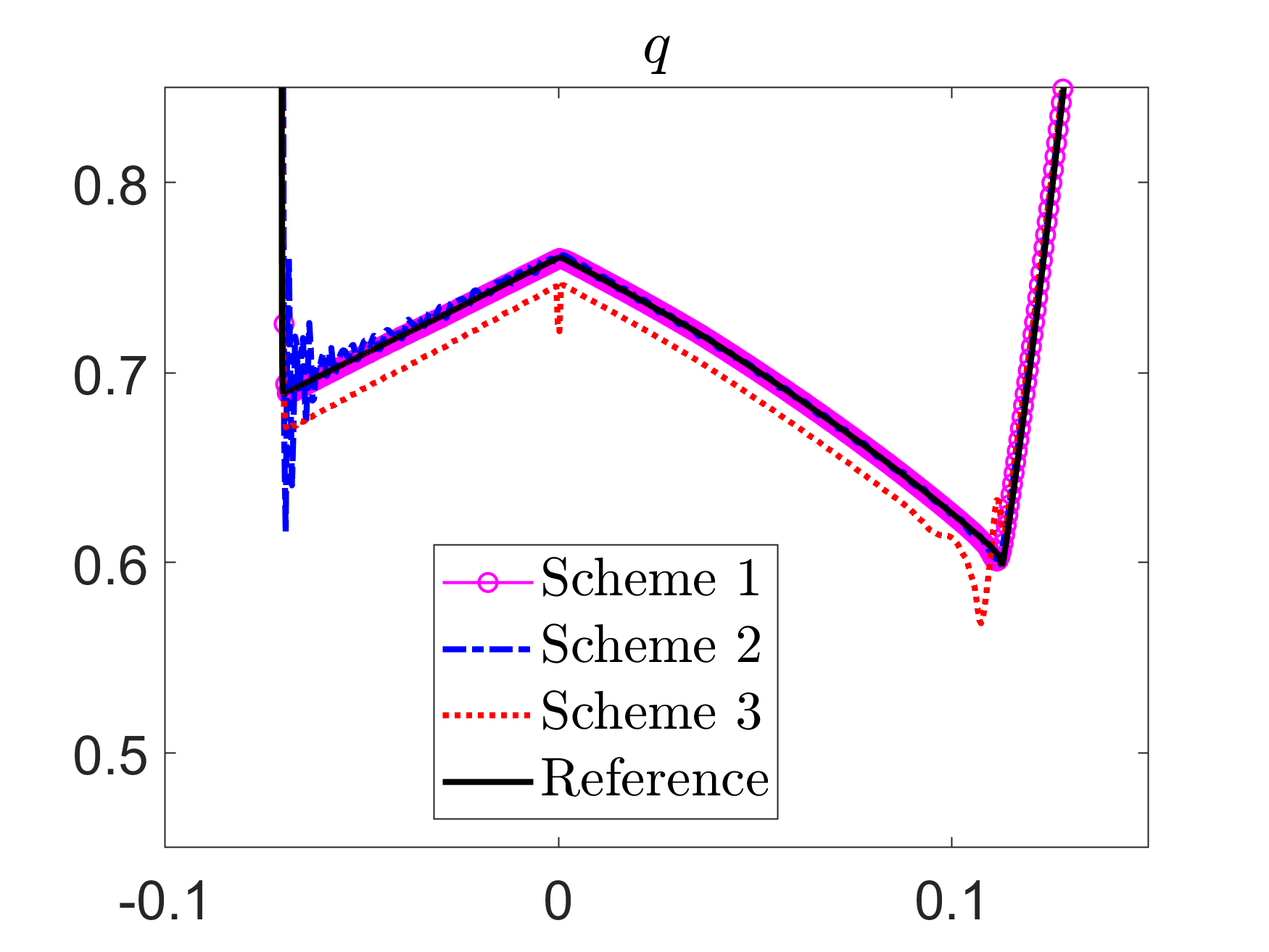}}
\vskip7pt
\centerline{\includegraphics[trim=0.4cm 0.3cm 0.9cm 0.1cm, clip, width=5.7cm]{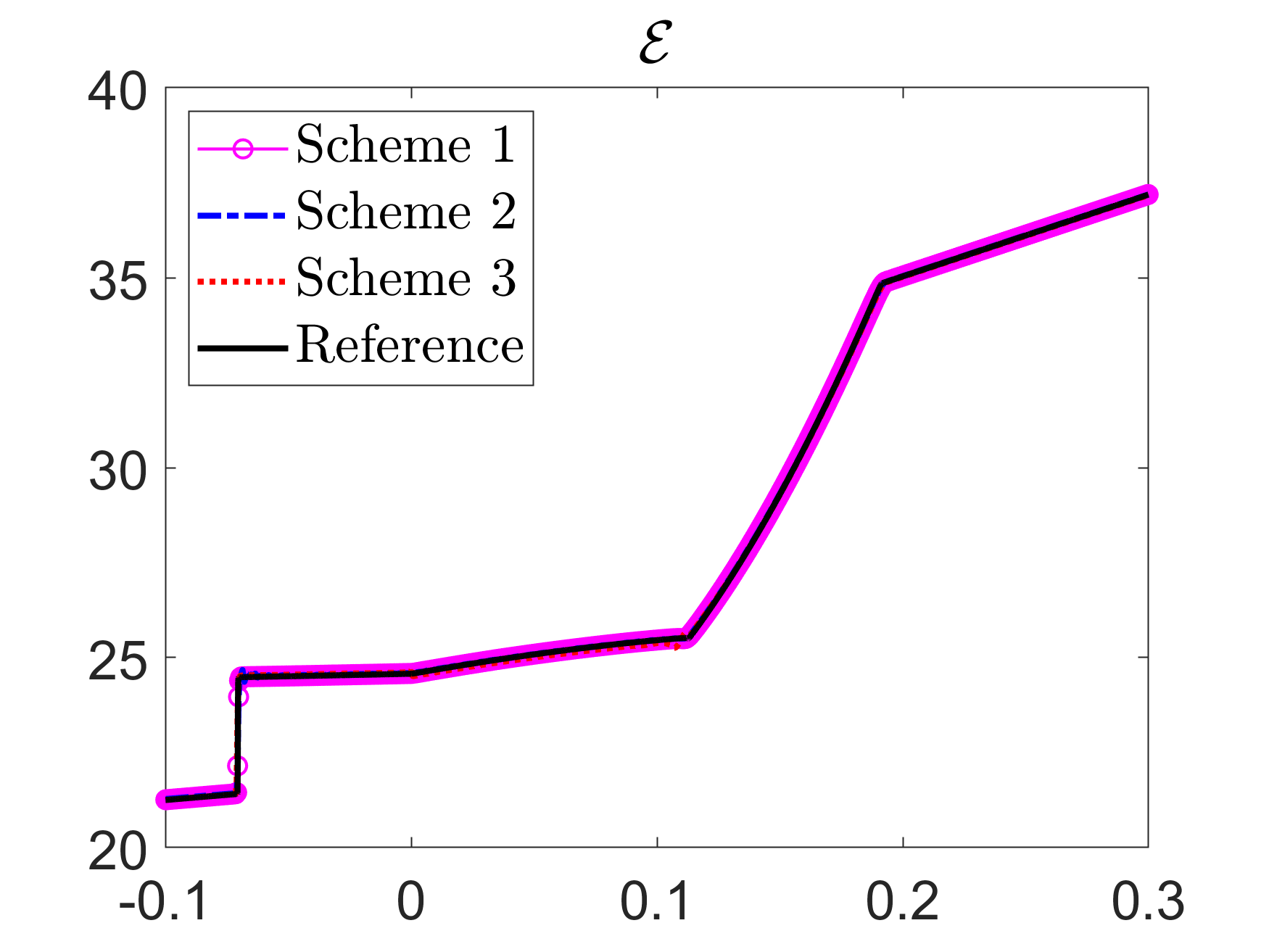}\hspace*{0.3cm}
            \includegraphics[trim=0.4cm 0.3cm 0.9cm 0.1cm, clip, width=5.7cm]{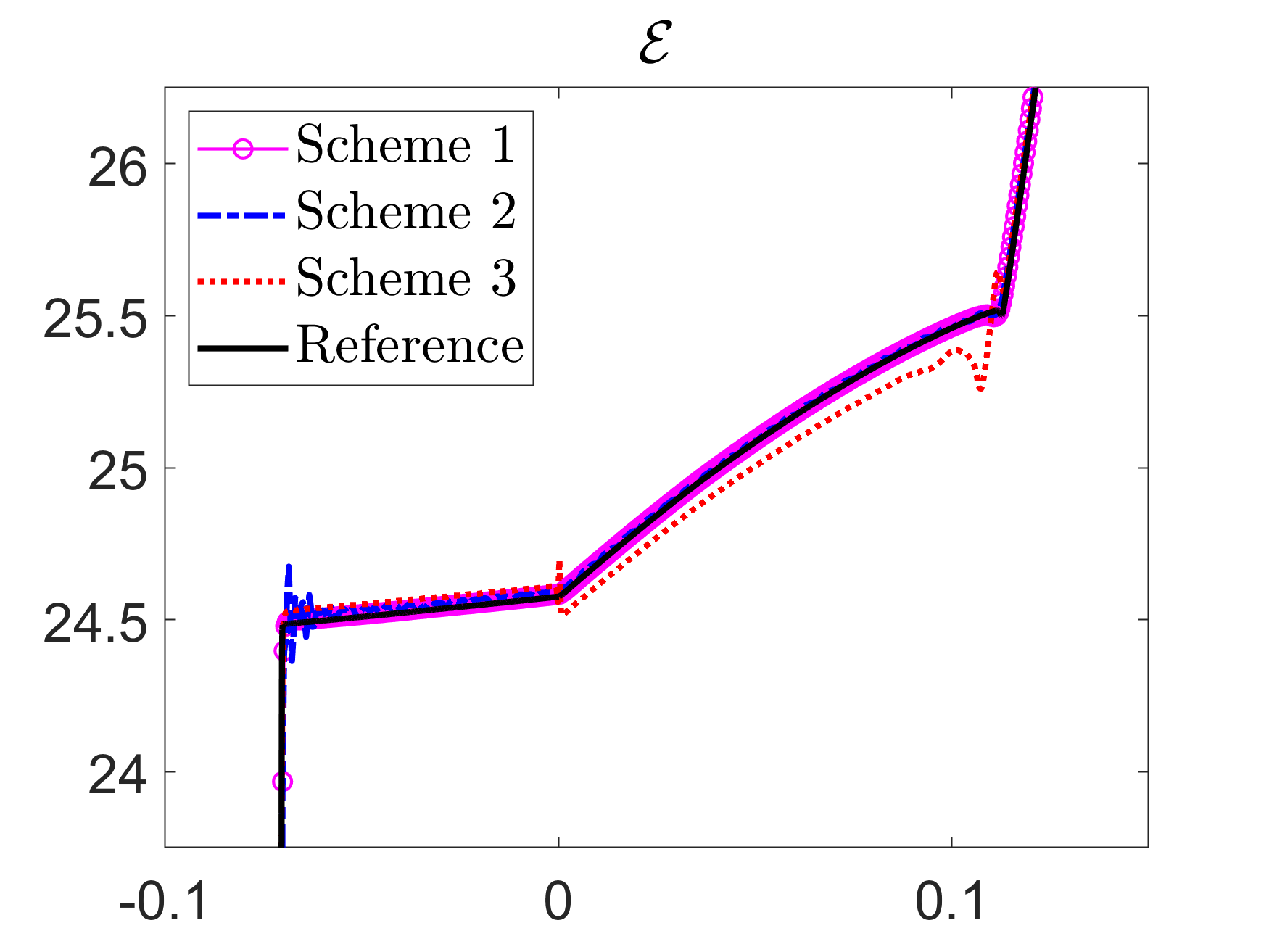}}
\caption{\sf Example 3: Water depth $h$, discharge $q$, and energy ${\cal E}$ computed by Schemes 1--3 using $\dx=1/2500$ (left column), and
zoom at the areas containing oscillations in Scheme 2 solution (right column).\label{fig212}}
\end{figure}

We compute the numerical solutions by the three studied schemes until the final time $t=0.03$ on a uniform mesh with
$\dx=1/200$ together with the reference solution computed by Scheme 1 on a much finer mesh with $\dx=1/10000$. The obtained results ($h$,
$q$, and ${\cal E}$) are shown in Figure \ref{fig211}, where one can observe that the solutions computed by both Schemes 2 and 3 are
oscillatory, whereas Scheme 1 solution is oscillation-free. In order to further study the performance of Schemes 1--3, we refine the mesh to
$\dx=1/2500$ and compare the high-resolution solutions of the three studied schemes together with the reference solution computed by Scheme
1 on a much finer mesh with $\dx=1/10000$; see Figure \ref{fig212}. One can observe that the high-frequency oscillations produced by Scheme
3 on a coarse mesh were almost suppressed when the mesh was refined, whereas Scheme 2, which does not employ any LCD, still produces
oscillations, which can be clearly seen in the zoom views shown in the right column.

\subsubsection*{Example 4---Convergence to a Steady State ($n=0.15$)}
In this example, we study the convergence of the solutions computed by Schemes 1--3 towards the steady flow over a hump. We consider the
continuous bottom topography given by
\begin{equation*}
Z(x)=\left\{\begin{aligned}
&0.2&&\mbox{if}~8\le x\le12,\\
&0&&\mbox{otherwise},
\end{aligned}\right.
\end{equation*}
and the initial and boundary data that correspond to a subcritical flow:
\begin{equation*}
h(x,0)\equiv2-Z(x),\quad q(x,0)\equiv0,\quad q(0,t)=4.42,\quad h(25,t)=2,
\end{equation*}
with the boundary conditions for $h$ at $x=0$ and $q$ at $x=25$ set to be homogeneous Neumann.

\begin{figure}[ht!]
\centerline{\includegraphics[trim=0.1cm 0.3cm 0.9cm 0.1cm, clip, width=5.7cm]{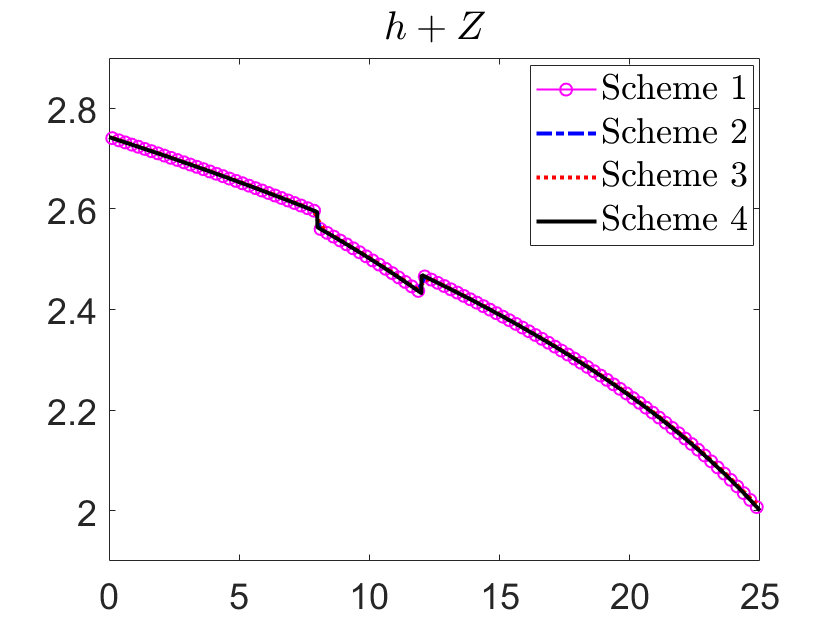}\hspace*{0.0cm}
            \includegraphics[trim=0.1cm 0.3cm 0.9cm 0.1cm, clip, width=5.7cm]{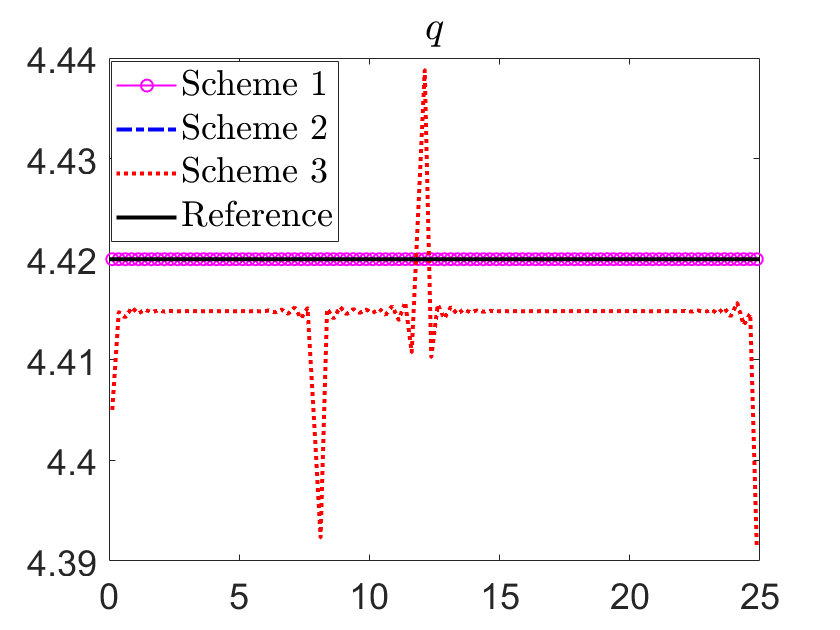}\hspace*{0.0cm}
            \includegraphics[trim=0.1cm 0.3cm 0.9cm 0.1cm, clip, width=5.7cm]{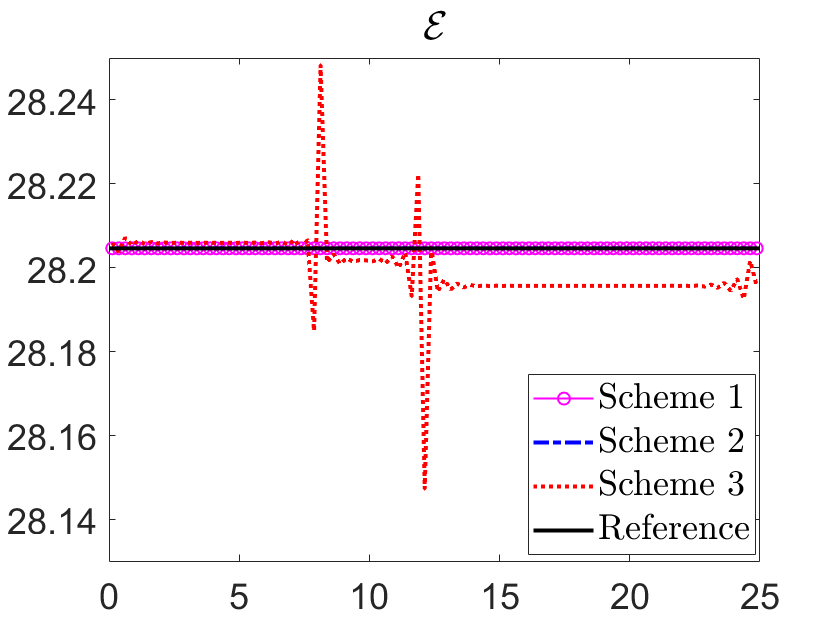}}
\caption{\sf Example 4 (moving water steady state): $h+Z$, $q$, and ${\cal E}$ computed by Schemes 1--3.\label{fig21}}
\end{figure}

We compute the numerical solutions until the final time $t=500$ on the computational domain $[0,25]$ covered by a uniform mesh with
$\dx=1/4$ together with the reference solution computed by Scheme 1 on a much finer mesh with $\dx=1/80$. The obtained numerical solutions
$h+Z$, $q$, and ${\cal E}-I_\hf$ (note that we subtract $I_\hf$ from $\cal E$ to obtain the quantity, which is supposed to be constant
disregarding what values of $\widehat x$ were used in the evaluation of ${\cal E}_j$) are plotted in Figure \ref{fig21}. One can clearly see
that the WB Schemes 1 and 2 converge to the constant $q$ and ${\cal E}$, whereas the non-WB Scheme 3 generates spurious oscillations.

\begin{figure}[ht!]
\centerline{\includegraphics[trim=0.4cm 0.3cm 0.9cm 0.1cm, clip, width=5.7cm]{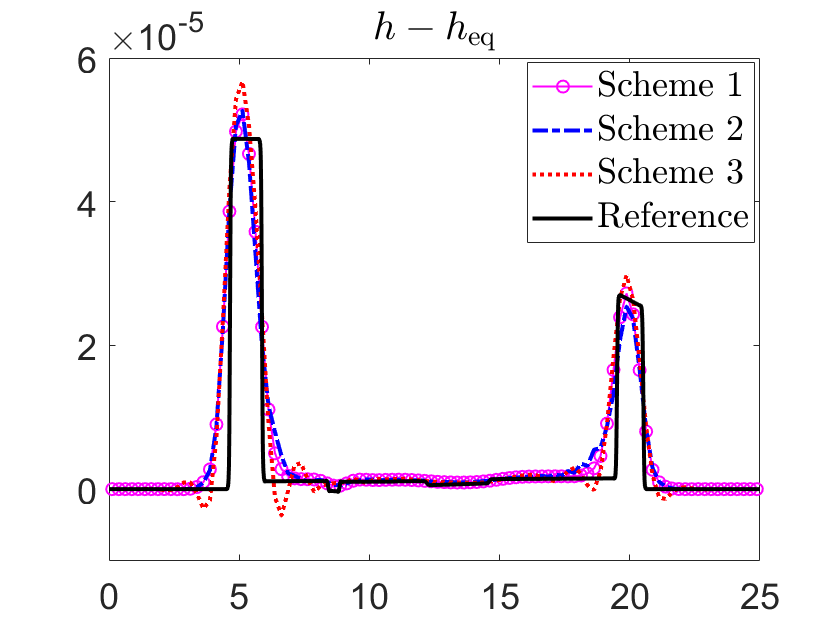}\hspace*{0.3cm}
            \includegraphics[trim=0.4cm 0.3cm 0.9cm 0.1cm, clip, width=5.7cm]{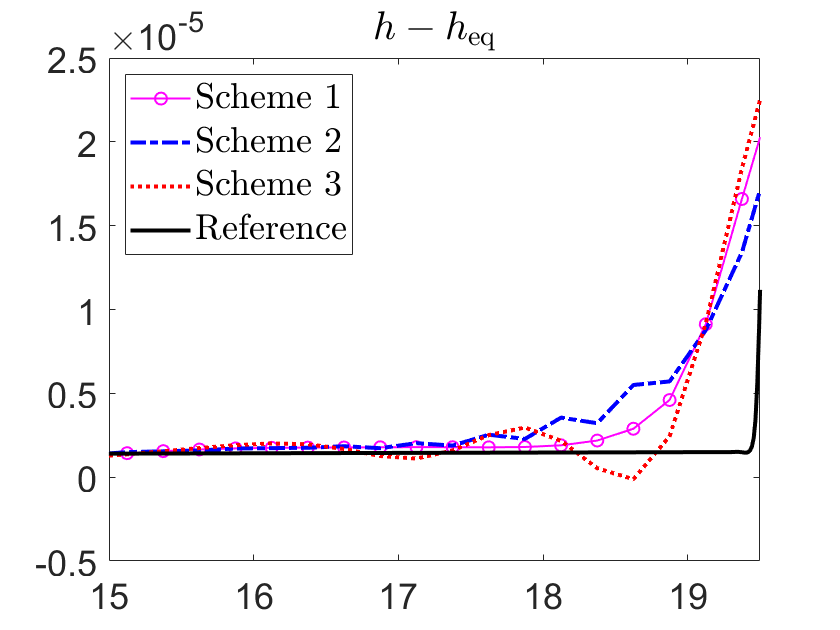}}
\caption{\sf Example 4: The differences $h(x,1.5)-h_{\rm eq}(x)$ computed by Schemes 1--3 (left) and zoom at $x\in[15,19.5]$ (right).\label{fig22}}
\end{figure}

We then test the ability of the studied schemes to capture the propagation of a small perturbation of the obtained moving-water equilibria.
To this end, we denote the obtained steady states by $h_{\rm eq}(x)$ and $q_{\rm eq}(x)$ (notice that each scheme has its own discrete
equilibrium), and then consider the following initial data:
\begin{equation*}
h(x,0)=h_{\rm eq}(x)+\left\{\begin{aligned}&10^{-4},&&9.5\le x\le10.5,\\&0,&&\mbox{otherwise},\end{aligned}\right.\qquad
q(x,0)=q_{\rm eq}(x).
\end{equation*}
We compute the solutions by the three studied schemes until the final time $t=1.5$ on the same uniform mesh with $\dx=1/4$. The obtained
differences $h(x,1.5)-h_{\rm eq}(x)$ are plotted in Figure \ref{fig22}. One can clearly see that both the WB Schemes 1 and 2 can capture the
time evolution of the perturbation quite accurately, whereas the non-WB Scheme 3 generates spurious oscillations as it can preserve
still-water equilibria only. Even though Scheme 2 is WB, it still generates some oscillations, and thus Scheme 1 clearly outperforms both of
its counterparts in this example.

\subsubsection*{Example 5---Efficiency vs. Accuracy ($n=0$)}
In this example, we assess both the WB property and the computational efficiency of Schemes 1--3 on a nontrivial moving-water steady state.
The bottom topography is
\begin{equation*}
Z(x)=-0.2e^{-40(x-10)^2},
\end{equation*}
the equilibrium state is
\begin{equation*}
E_{\rm eq}(x)\equiv32,\quad q_{\rm eq}(x)\equiv2,
\end{equation*}
and the discrete values of the corresponding equilibrium water depth $h_{\rm eq}(x_j)$ are obtained by solving the nonlinear equations
\begin{equation*}
\frac{4}{2h_{\rm eq}^2(x_j)}+g\left(h_{\rm eq}(x_j)+Z(x_j)\right)=32.
\end{equation*}
These data are prescribed in the computational domain $[0,25]$ subject to the homogeneous Neumann boundary conditions.

To test the ability of the three schemes to preserve this steady state and capture small perturbations, we take the following initial data:
\begin{equation*}
h(x,0)=h_{\rm eq}(x)+10^{-4}e^{-4(x-12)^2},\quad q(x,0)=q_{\rm eq}(x)\equiv2,
\end{equation*}
in which a small Gaussian-shaped perturbation is added to the equilibrium state.

\begin{figure}[ht!]
\centerline{\includegraphics[trim=0.4cm 0.3cm 0.9cm 0.1cm, clip, width=5.7cm]{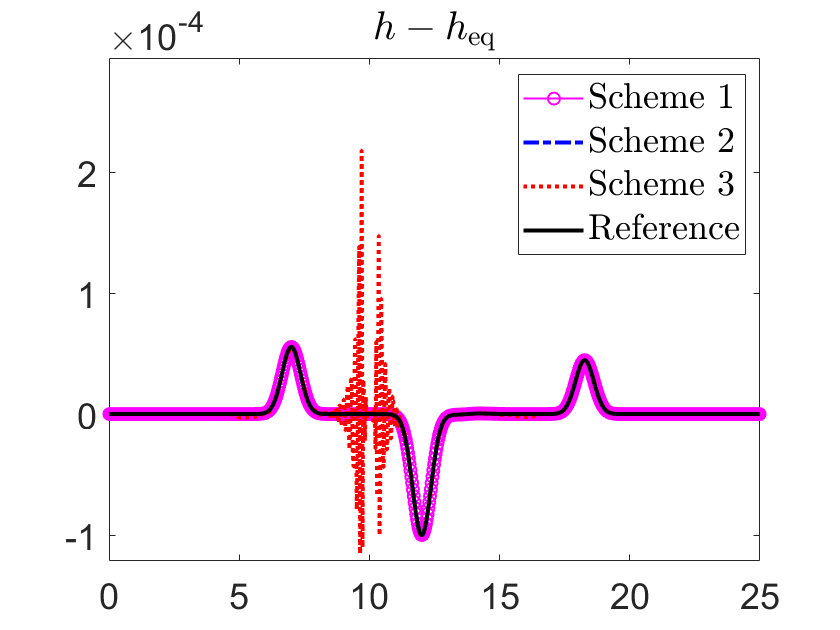}\hspace*{0.3cm}
            \includegraphics[trim=0.4cm 0.3cm 0.9cm 0.1cm, clip, width=5.7cm]{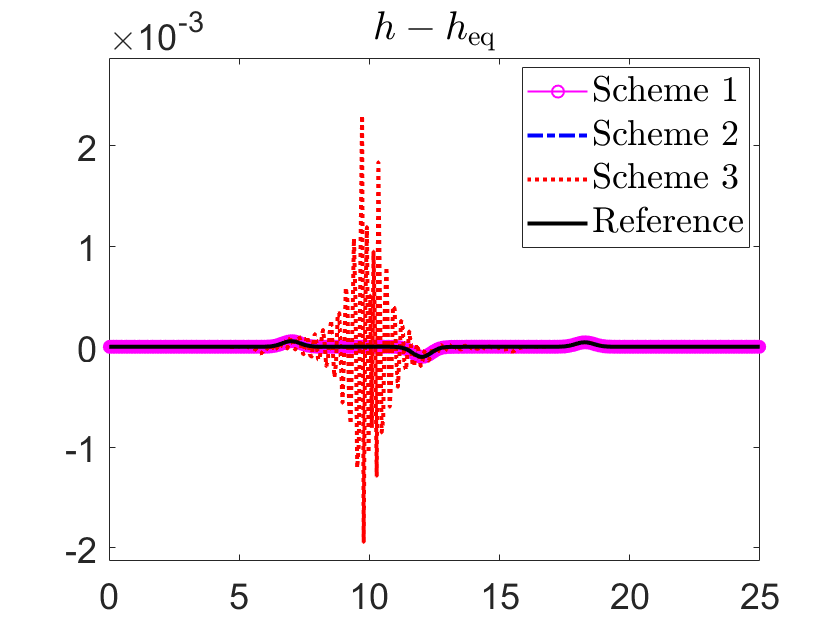}}
\caption{\sf Example 5: The differences $h(x,1)-h_{\rm eq}(x)$ computed by Schemes 1--3 for $\dx=1/16$ (left) and $1/32$ (right).
\label{fig23}}
\end{figure}

We compute the numerical solutions up to the final time $t=1$ on a sequence of uniform meshes with $\dx=1/2$, $1/4$, $1/8$, $1/16$, $1/32$,
and $1/64$. We measure the $L^1$-errors in $h$ using the Runge formula, which is based on the solutions computed on the three consecutive
uniform grids with the mesh sizes $\dx$, $2\dx$, and $4\dx$ and denoted by $(\cdot)^{\dx}$, $(\cdot)^{2\dx}$, and $(\cdot)^{4\dx}$,
respectively:
\begin{equation}
{\rm Error}(\dx)\approx\frac{\delta_{12}^2}{|\delta_{12}-\delta_{24}|},
\label{runge}
\end{equation}
where $\delta_{12}:=\|(\cdot)^{\dx}-(\cdot)^{2\dx}\|_{L^1}$ and $\delta_{24}:=\|(\cdot)^{2\dx}-(\cdot)^{4\dx}\|_{L^1}$. We report the
$L^1$-errors along the CPU times (in seconds) in Table \ref{tab41f}, where one can see that Scheme 3 is the least computationally expensive,
but it produces errors several orders of magnitude larger than those produced by Schemes 1 and 2. Concerning the WB property, it is also
evident that Schemes 1 and 2 capture the small perturbation of the steady state very accurately, while Scheme 3 fails to suppress spurious
oscillations, resulting in substantially larger errors. The latter can be clearly seen in Figure \ref{fig23}, where we plot
$h(x,1)-h_{eq}(x)$ computed by the three studied schemes on coarse meshes with $\dx=1/16$ and $1/32$ (here, the finest-mesh solution
computed by Scheme 1 is being considered as the reference solution).
\begin{table}[ht!]
\centering
\begin{tabular}{ccccccc}
\toprule
\multirow{2}{*}{$\dx$}&\multicolumn{2}{c}{Scheme 1}&\multicolumn{2}{c}{Scheme 2}&\multicolumn{2}{c}{Scheme 3}\\
\cline{2-7}
      &Error                &CPU   &Error                &CPU   &Error               &CPU\\
\hline
$1/16$&$9.57\mathrm{e}{-08}$&$1.66$&$8.69\mathrm{e}{-07}$&$1.50$&$4.93\mathrm{e}{-02}$&$0.78$\\
$1/32$&$5.57\mathrm{e}{-09}$&$6.56$&$1.40\mathrm{e}{-09}$&$5.98$&$3.28\mathrm{e}{-03}$&$2.95$\\
$1/64$&$6.48\mathrm{e}{-10}$&$25.7$&$6.67\mathrm{e}{-10}$&$23.5$&$2.54\mathrm{e}{-05}$&$11.6$\\
\bottomrule 
\end{tabular}
\caption{$L^1$-errors in $h$ and CPU times (in seconds).}
\label{tab41f}
\end{table}

\subsection{Two-Layer Shallow Water System}
\subsubsection*{Example 6---Experimental Order of Accuracy}
In this example taken from \cite{KPmultil}, we consider the following initial data:
\begin{equation*}
h_1(x,0)=5+e^{\cos(2\pi x)},\quad h_2(x,0)=5-e^{\cos(2\pi x)}-\sin^2(\pi x),\quad q_1(x,0)=q_2(x,0)\equiv0,
\end{equation*}
and a continuous bottom topography
\begin{equation*}
Z(x)=\sin^2(\pi x)-10,
\end{equation*}
both prescribed in the computational domain $[0,1]$ subject to the periodic boundary conditions.

We compute the numerical solution until the final time $t=0.1$ by the studied Schemes 1--3 on a sequence of uniform meshes with $\dx=1/40$,
$1/80$, $1/160$, $1/320$, and $1/640$. We measure the $L^1$-errors in $h_1$ using the Runge formula \eref{runge} and estimate the
experimental convergence rates as follows:
$$
{\rm Rate}(\dx)\approx\log_2\left(\frac{\delta_{24}}{\delta_{12}}\right).
$$
The obtained results are reported in Table \ref{tab71}, where one can clearly see that the expected orders of accuracy are achieved by all
the three studied schemes. Note that in order to achieve the fifth order of accuracy, we have used smaller time steps with
$\dt\sim(\dx)^\frac{5}{3}$.
\begin{table}[ht!]
\centering
\begin{tabular}{ccccccccc}
\toprule
\multirow{2}{1em}{$\dx$}&\multicolumn{2}{c}{Scheme 1}&\multicolumn{2}{c}{Scheme 2}&\multicolumn{2}{c}{Scheme 3}\\
\cline{2-7}&Error&Rate&Error&Rate&Error&Rate\\
\hline
$1/160$&$1.90\mathrm{e}{-07}$&$4.83$&$1.79\mathrm{e}{-07}$&$5.05$&$6.14\mathrm{e}{-08}$&$4.78$\\
$1/320$&$5.10\mathrm{e}{-09}$&$5.02$&$5.01\mathrm{e}{-09}$&$5.10$&$1.67\mathrm{e}{-09}$&$4.99$\\
$1/640$&$1.57\mathrm{e}{-10}$&$5.02$&$1.58\mathrm{e}{-10}$&$5.04$&$5.13\mathrm{e}{-11}$&$5.01$\\
\bottomrule
\end{tabular}
\caption{\sf Example 6: The $L^1$-errors and experimental convergence rates for $h_1$.\label{tab71}}
\end{table}

\subsubsection*{Example 7---Small Perturbation of a Discontinuous Steady State}
In this example taken from \cite{KLX_21,CKX_24WB}, we consider a discontinuous steady state given by
\begin{equation*}
\begin{aligned}
(h_1)_{\rm eq}(x):&=\begin{cases}1.22373355048230,&x<0,\\1.44970064153589,&x>0,\end{cases}\qquad(q_1)_{\rm eq}(x)\equiv12,\\
(h_2)_{\rm eq}(x):&=\begin{cases}0.968329515483846,&x<0,\\1.12439026921484,&x>0,\end{cases}\qquad(q_2)_{\rm eq}(x)\equiv10,\\
\end{aligned}
\end{equation*}
and a discontinuous bottom topography
\begin{equation*}
Z(x)=\begin{cases}-2,&x<0,\\-1,&x>0.\end{cases}
\end{equation*}
In order to test the ability of the studied schemes to capture quasi-steady solutions, we add a small perturbation to the upper layer depth
and take the following initial data:
\begin{equation*}
\begin{aligned}
h_1(x,0)&=(h_1)_{\rm eq}(x)+\begin{cases}0.04,&x\in[-0.9,-0.8],\\0,&\mbox{otherwise},\end{cases}\\
h_2(x,0)&=(h_2)_{\rm eq}(x),\quad q_1(x,0)=(q_1)_{\rm eq}(x),\quad q_2(x,0)=(q_2)_{\rm eq}(x),
\end{aligned}
\end{equation*}
prescribed in the computational domain $[-1,1]$ subject to the homogeneous Neumann boundary conditions.

\begin{figure}[ht!]
\centerline{\includegraphics[trim=0.0cm 0.3cm 0.7cm 0.2cm, clip, width=5.7cm]{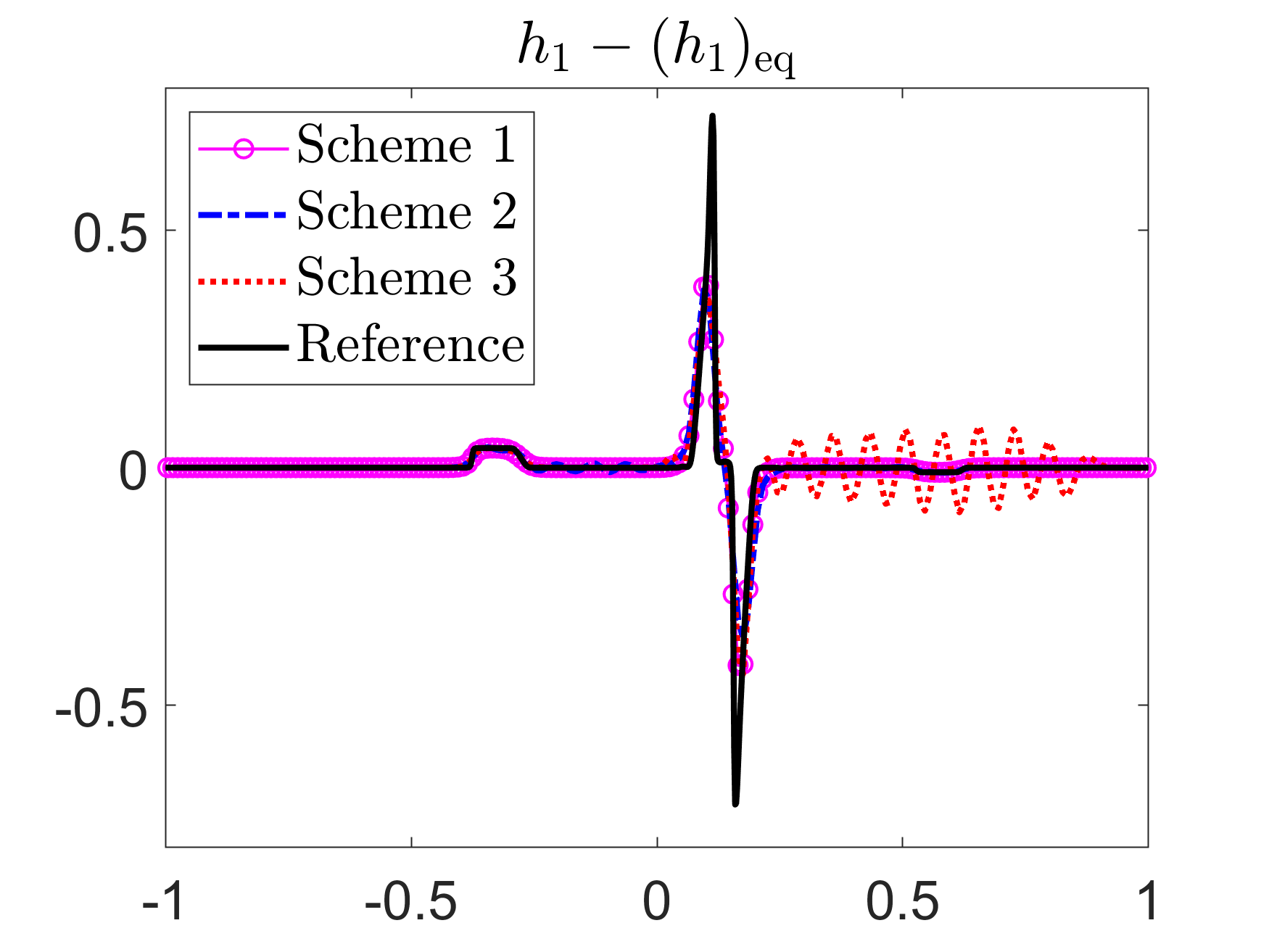}\hspace{0.3cm}
            \includegraphics[trim=0.0cm 0.3cm 0.7cm 0.2cm, clip, width=5.7cm]{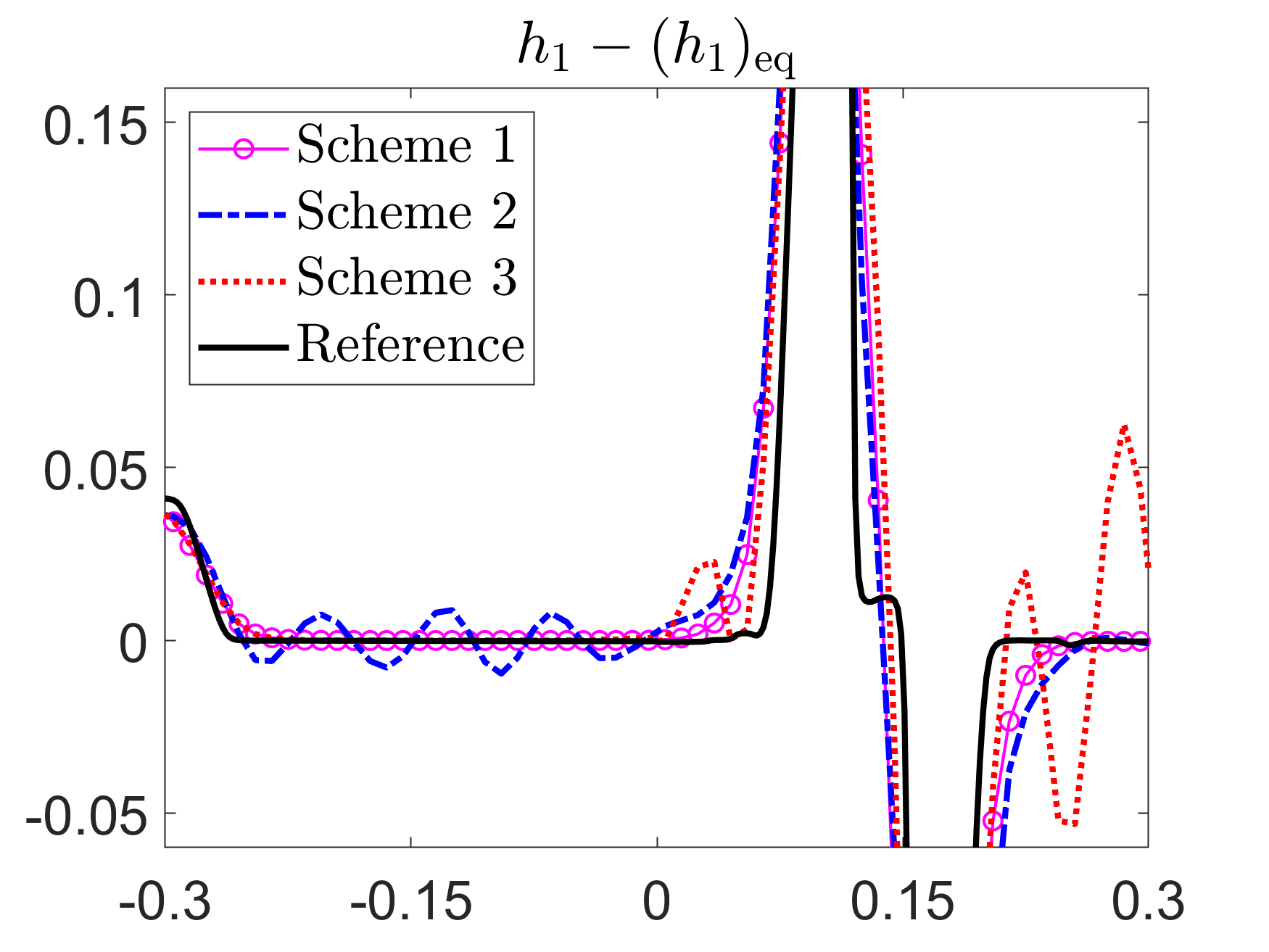}}
\vskip7pt
\centerline{\includegraphics[trim=0.0cm 0.3cm 0.7cm 0.2cm, clip, width=5.7cm]{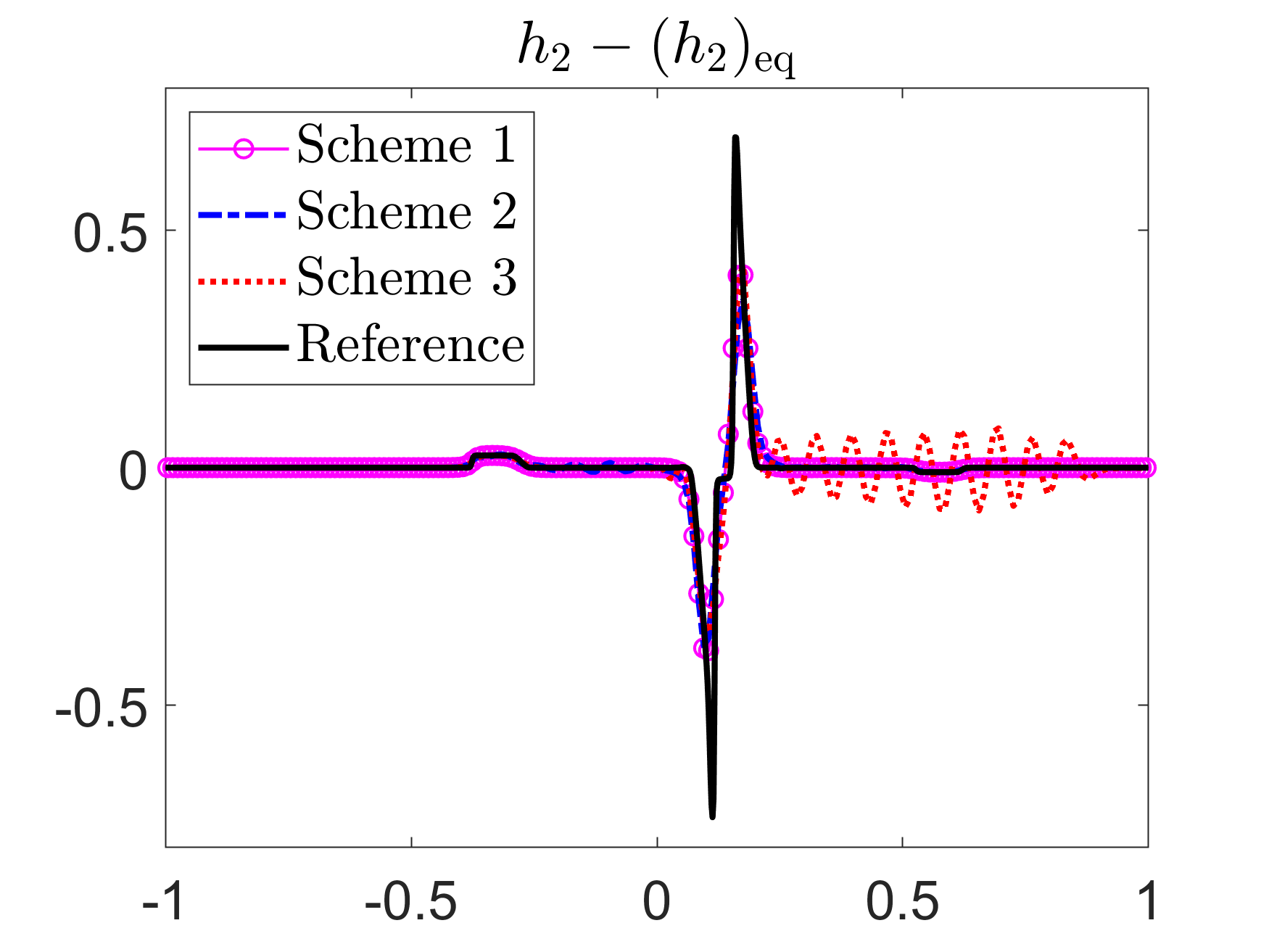}\hspace{0.3cm}
            \includegraphics[trim=0.0cm 0.3cm 0.7cm 0.2cm, clip, width=5.7cm]{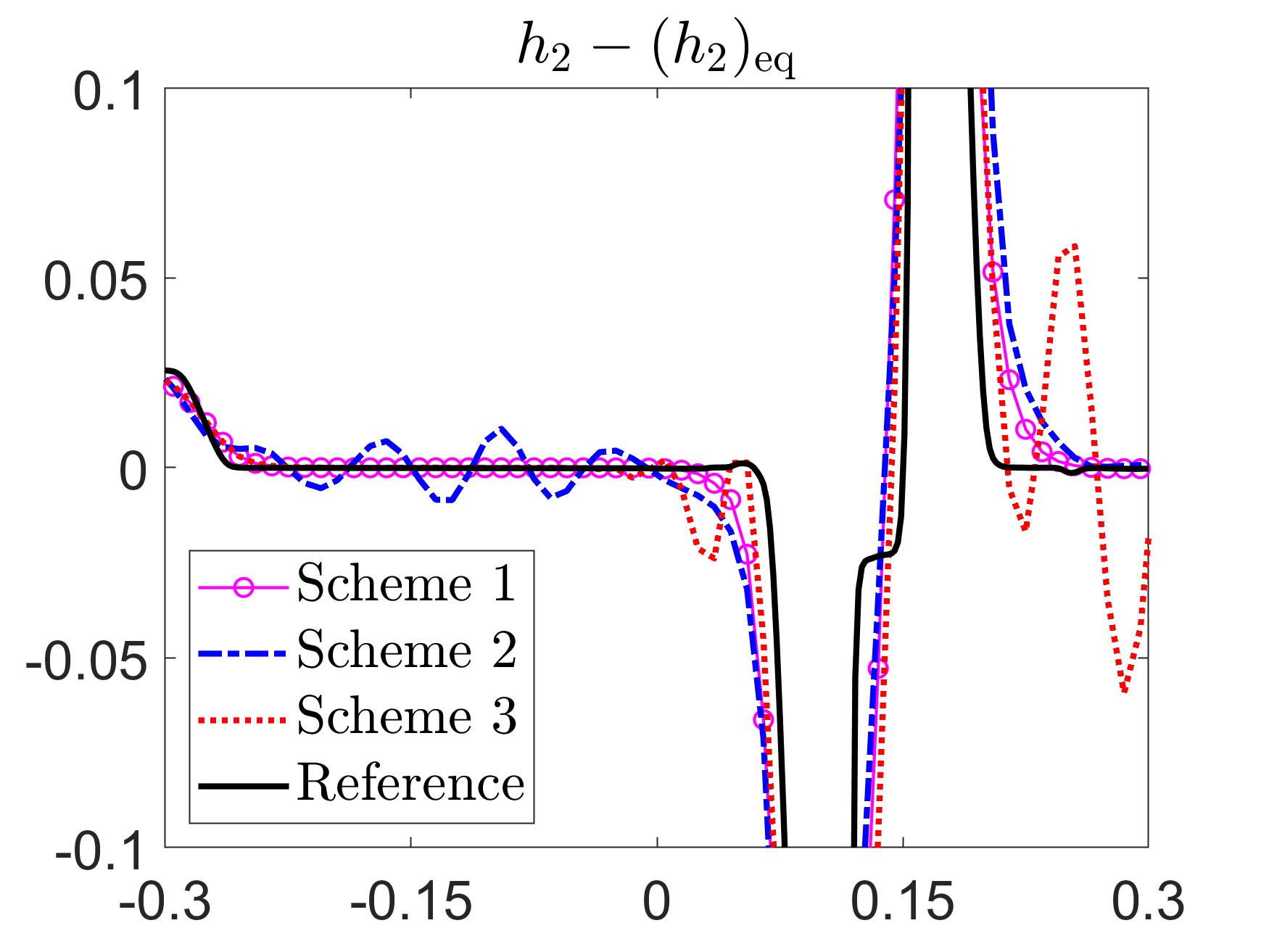}}
\caption{\sf Example 7: The differences $h_1(x,0.1)-(h_1)_{\rm eq}(x)$ (top row) and $h_2(x,0.1)-(h_2)_{\rm eq}(x)$ (bottom row), and zoom
at $x\in[-0.3,0.3]$ (right column).\label{fig7}}
\end{figure}

We compute the numerical solutions until the final time $t=0.1$ by Schemes 1--3 on a uniform mesh with $\dx=1/100$ and obtain the reference
solution using Scheme 1 on a much finer mesh with $\dx=1/1000$. The differences $h_1(x,0.1)-(h_1)_{\rm eq}(x)$ and
$h_2(x,0.1)-(h_2)_{\rm eq}(x)$ are plotted in Figure \ref{fig7}, where one can clearly see that, unlike the proposed Scheme 1, Schemes 2 and
3 produce oscillatory numerical results.

\subsubsection*{Example 8---Riemann Problem}
In this example taken from \cite{KLX_21,CKX_24WB}, we numerically solve a Riemann problem with the following initial data:
\begin{equation*}
(h_1,q_1,h_2,q_2)(x,0)=\begin{cases}(1,1.5,1,1),&x<0,\\(0.8,1.2,1.2,1.8),&\mbox{otherwise},\end{cases}
\end{equation*}
and discontinuous bottom topography:
\begin{equation*}
Z(x)=\begin{cases}-2,&x<0,\\-1.5,&\mbox{otherwise},\end{cases}
\end{equation*}
prescribed in the computational domain $[-1,1]$ subject to the homogeneous Neumann boundary conditions.

\begin{figure}[ht!]
\centerline{\includegraphics[trim=0.0cm 0.3cm 0.7cm 0.2cm, clip, width=5.7cm]{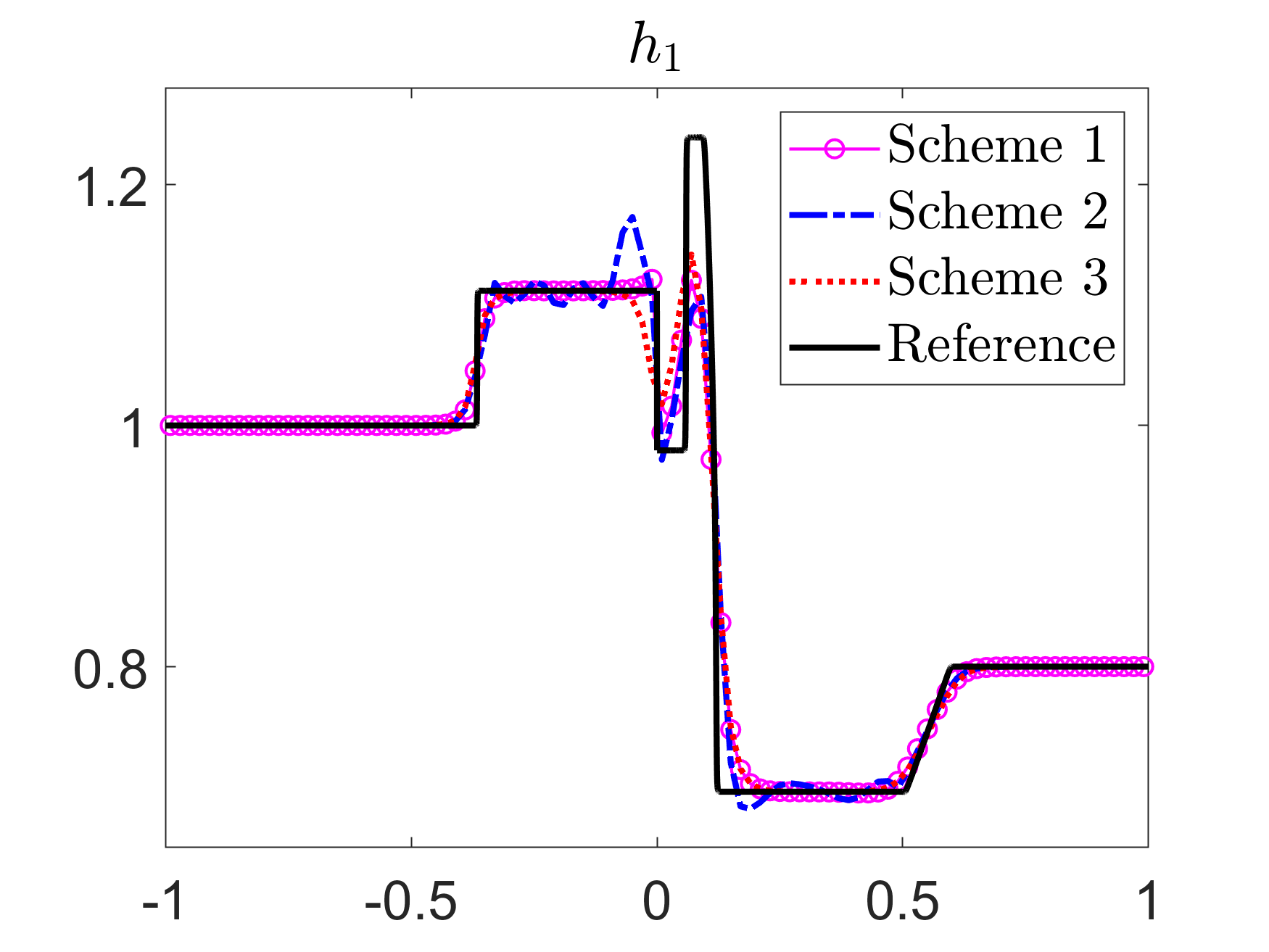}\hspace{0.3cm}
            \includegraphics[trim=0.0cm 0.3cm 0.7cm 0.2cm, clip, width=5.7cm]{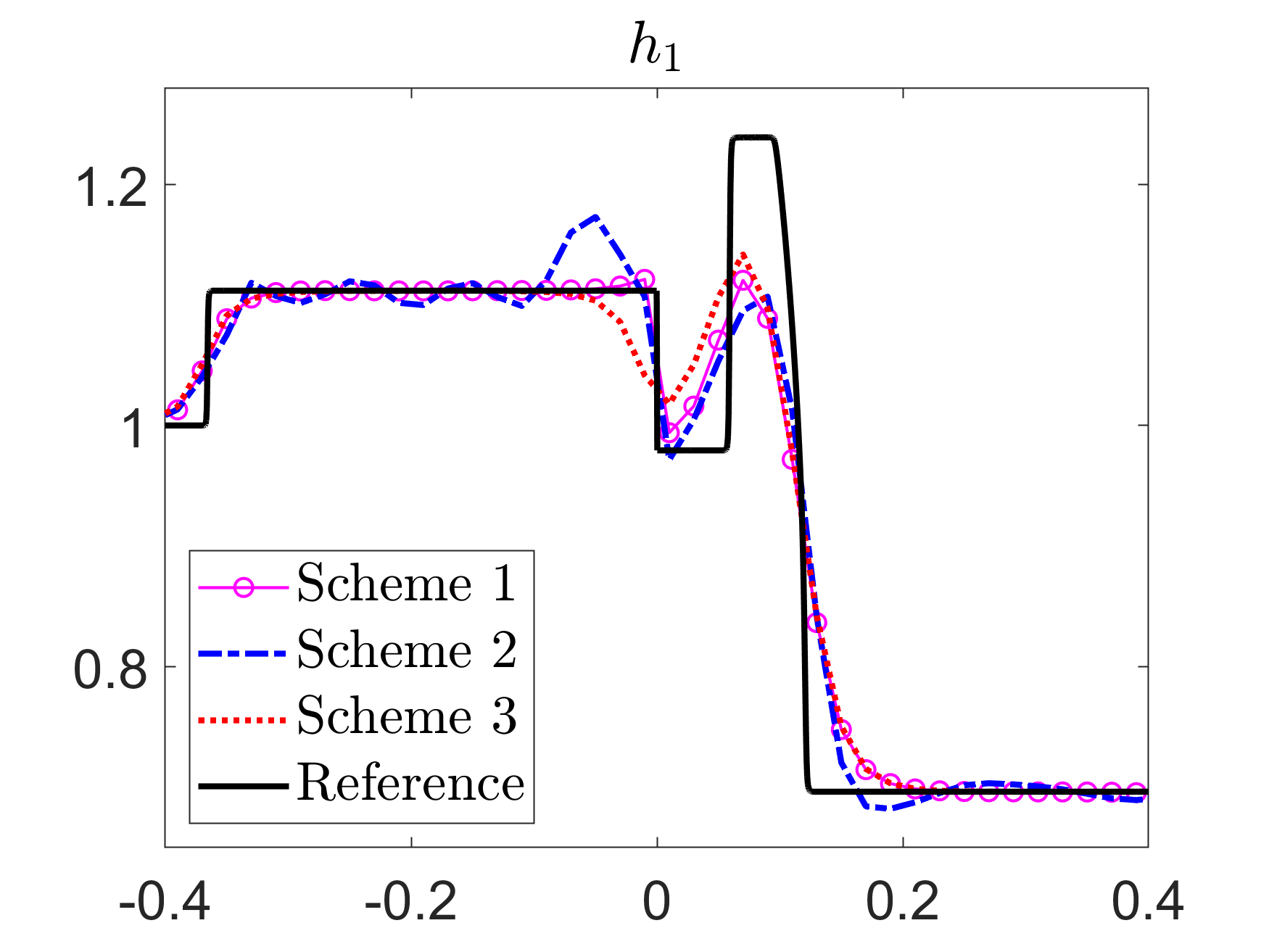}}
\vskip7pt
\centerline{\includegraphics[trim=0.0cm 0.3cm 0.7cm 0.2cm, clip, width=5.7cm]{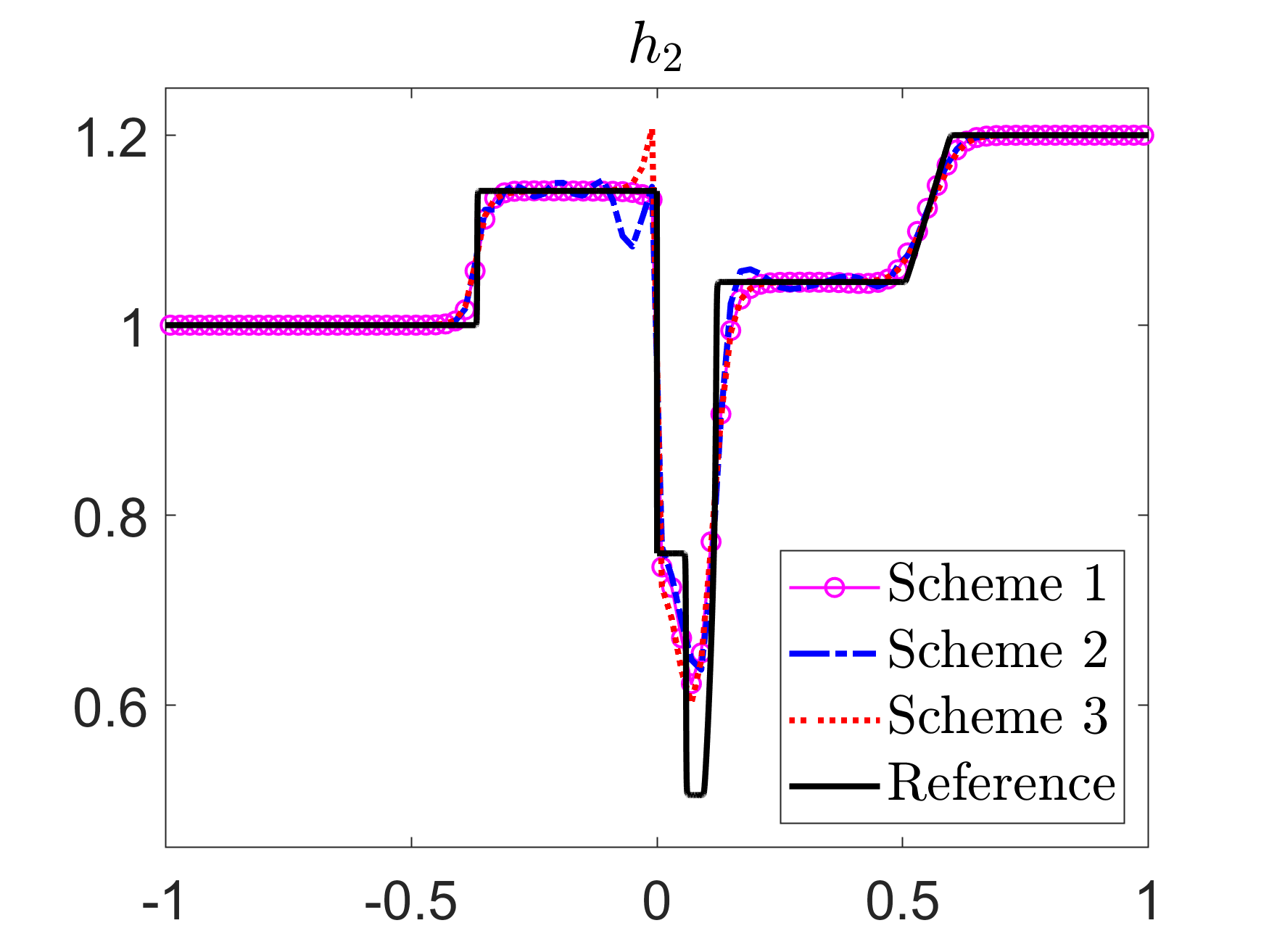}\hspace{0.3cm}
            \includegraphics[trim=0.0cm 0.3cm 0.7cm 0.2cm, clip, width=5.7cm]{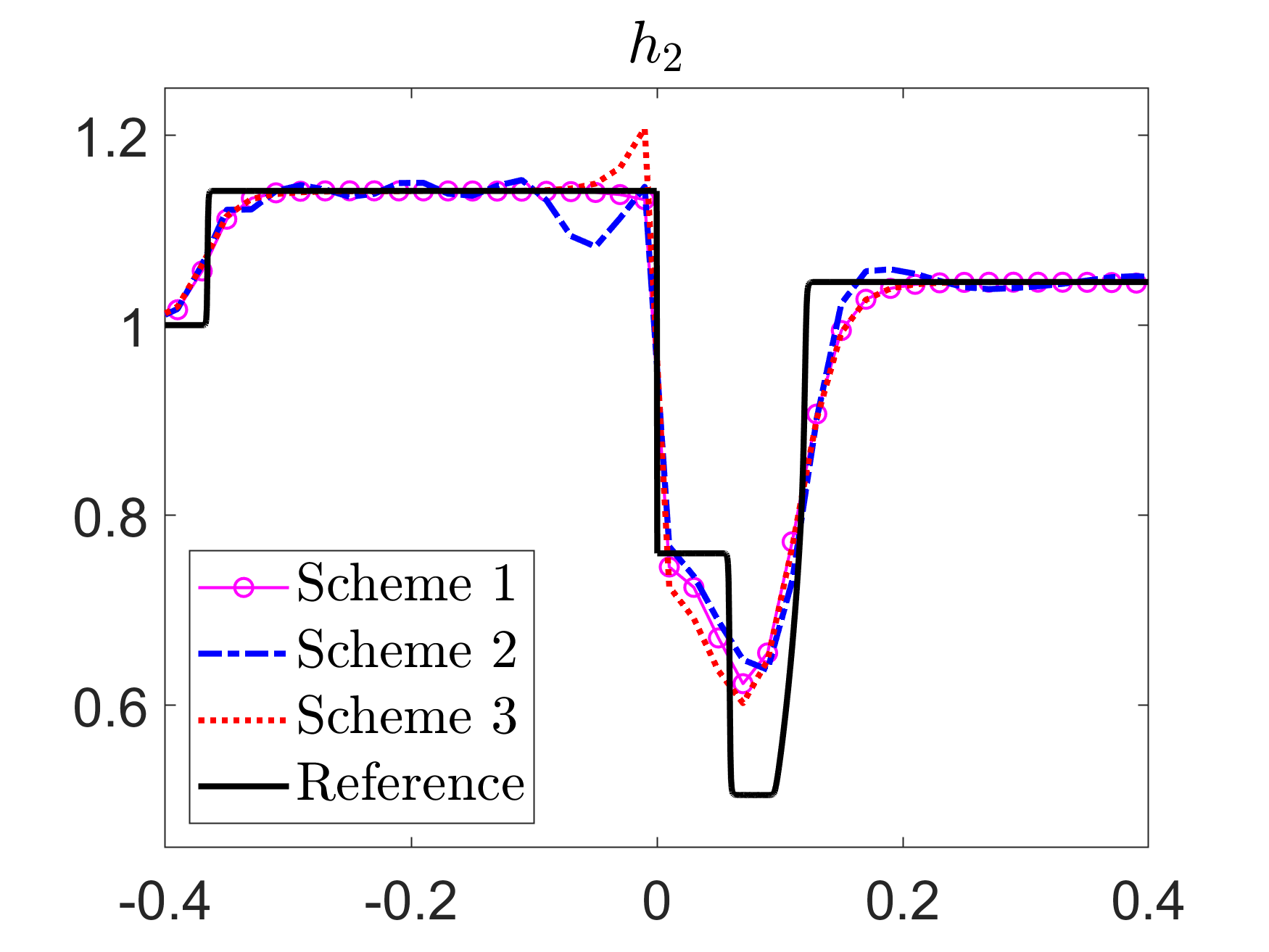}}
\caption{\sf Example 8: Upper layer depth $h_1$ (top row) and lower layer depth $h_2$ (bottom row), and zoom at $x\in[-0.4,0.4]$ (right
column).\label{fig9}}
\end{figure}

We compute the numerical solutions until the final time $t=0.1$ by Schemes 1--3 on a uniform mesh with $\dx=1/50$. The obtained upper layer
depth $h_1$ and lower layer depth $h_2$ are plotted in Figure \ref{fig9} together with the reference solution computed by Scheme 1 on a much
finer mesh with $\dx=1/2000$. As one can see, Scheme 1 clearly outperforms Schemes 2 and 3 as there are no oscillations in the results
computed by Scheme 1.

\subsection{1-D Euler Equations with Gravitation}
\subsubsection*{Example 9---Shock Tube Problem}
In this example, which is a modification of the example studied in \cite{Xing13,Luo11}, we set a nonlinear gravitational potential
$\phi(x)=\frac{1}{10x+1}$ in the computational domain $[0,1]$, where we prescribe the following initial data:
\begin{equation*}
(\rho(x,0),u(x,0),p(x,0))=\left\{\begin{aligned}&(1,0,1),&&x\le0.5,\\&(0.125,0,0.1),&&x>0.5,\end{aligned}\right.
\end{equation*}
and the reflecting (solid wall) boundary conditions, which are implemented using the ghost cell technique: we set the same values of $\rho$
and $p$ in the ghost cells while for $u$ the sign is negated.

\begin{figure}[ht!]
\centerline{\includegraphics[trim=0.6cm 0.3cm 1.3cm 0.1cm, clip, width=5.7cm]{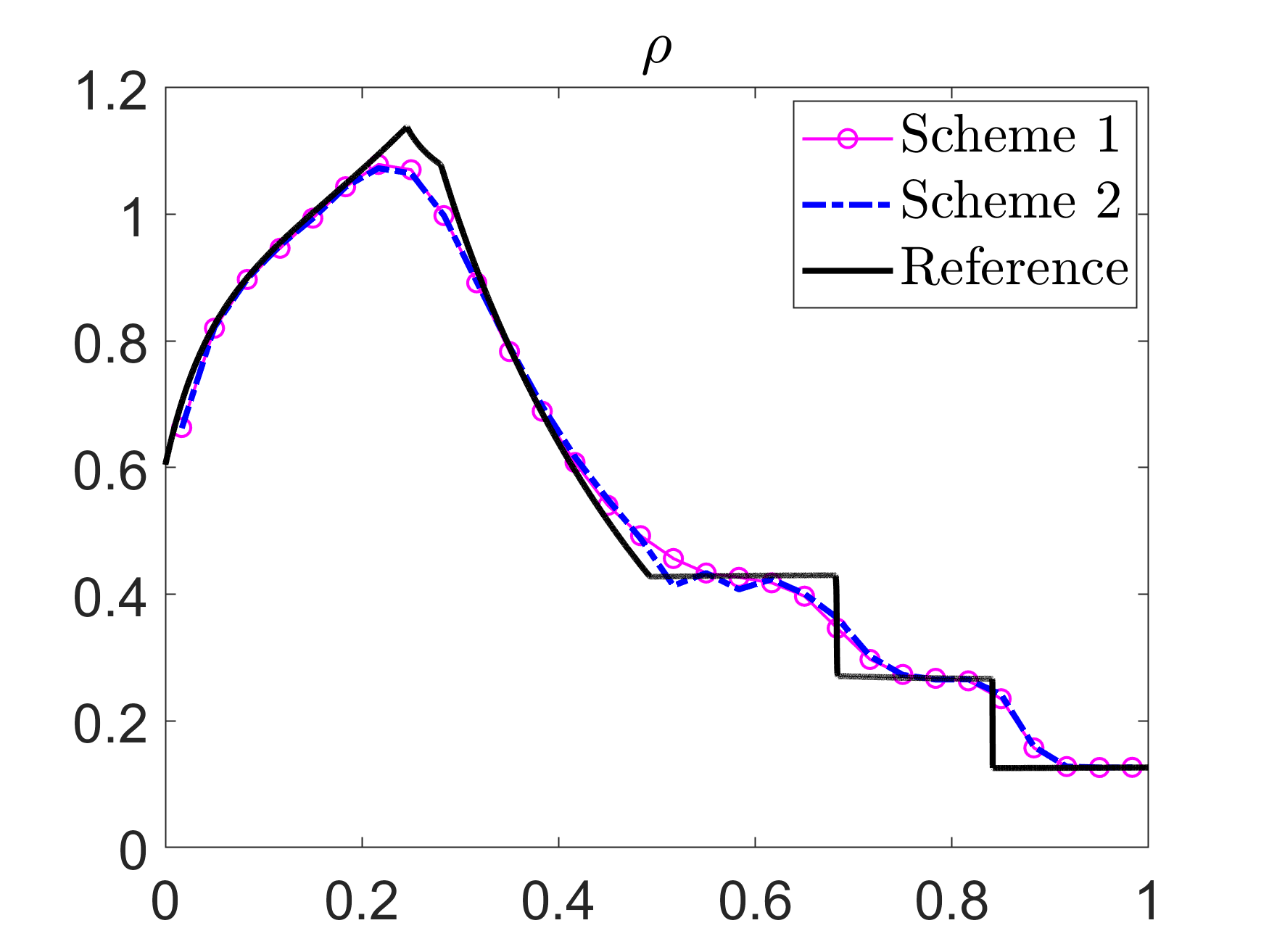}\hspace*{0.0cm}
            \includegraphics[trim=0.6cm 0.3cm 1.3cm 0.1cm, clip, width=5.7cm]{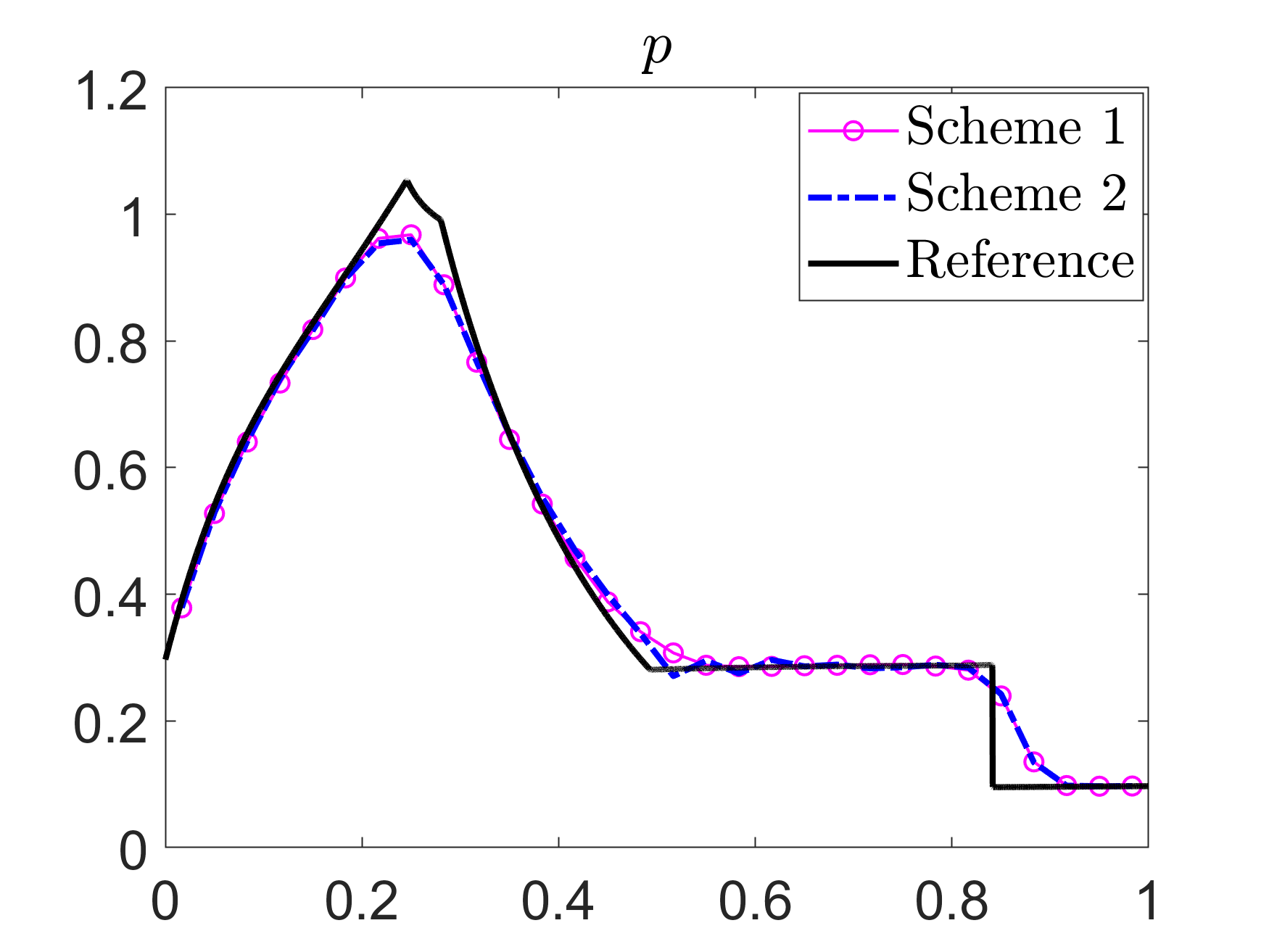}\hspace*{0.2cm}
	    \includegraphics[trim=0.6cm 0.3cm 1.3cm 0.1cm, clip, width=5.7cm]{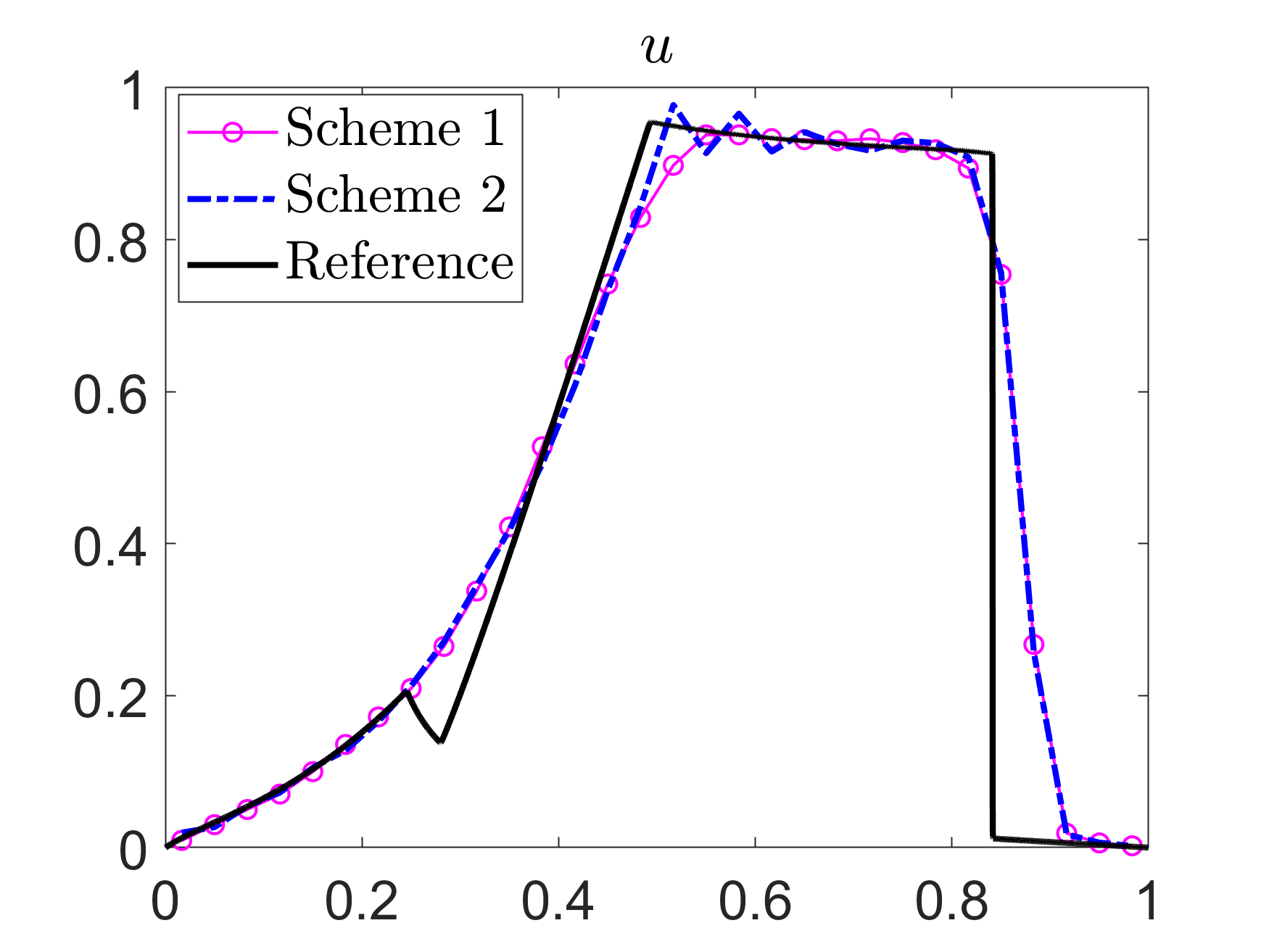}}
\caption{\sf Example 9: $\rho$, $p$, and $u$ computed by Schemes 1 and 2.\label{fig31}}
\end{figure}

We compute the numerical solutions until the final time $t=0.2$ by Schemes 1 and 2 on a uniform grid with $\dx=1/30$ together with the
reference solution computed by Scheme 1 on a much finer mesh with $\dx=1/8000$. The numerical results are shown in Figure \ref{fig31}, where
one can clearly see that Scheme 1 outperforms Scheme 2 as there are no oscillations in the results computed by Scheme 1.

\subsection{2-D Euler Equations with Gravitation}
\subsubsection*{Example 10---Discontinuous Perturbation of a Steady State}
In this example, which is a modification of the example studied in \cite{CCKOT,Xing13}, we consider the following hydrostatic equilibrium:
\begin{equation*}
\rho_{\rm eq}(x,y)=1.21e^{-1.21\phi(x,y)},\quad u_{\rm eq}(x,y)=v_{\rm eq}(x,y)\equiv0,\quad p_{\rm eq}(x,y)=e^{-1.21\phi(x,y)},
\end{equation*}
with the gravitational potential $\phi(x,y)=x+y$. We take the computational domain $[0,1]\times[0,1]$ and impose the homogeneous Neumann
boundary conditions.

We construct the discrete steady state as it was described in \cite[Example 7]{Kurganov25Na}, but with the integrals in $K^x$ and $K^y$
evaluated within the fifth order of accuracy.

Equipped with the discrete steady state, we first numerically verify that Schemes 1 and 2 can preserve it within the machine accuracy, and
then introduce a discontinuous pressure perturbation and consider the perturbed initial data:
\begin{equation*}
(\rho,u,v)(x,y,0)=(\rho_{\rm eq},u_{\rm eq},v_{\rm eq})(x,y),\quad p(x,y,0)=p_{\rm eq}(x,y)+
\left\{\begin{aligned}&0.5,&&(x-0.5)^2+(y-0.5)^2\le0.2^2,\\&0,&&\mbox{otherwise}.\end{aligned}\right.
\end{equation*}

\begin{figure}[ht!]
\centerline{\includegraphics[trim=0.6cm 0.6cm 1.2cm 0.3cm, clip, width=5.3cm]{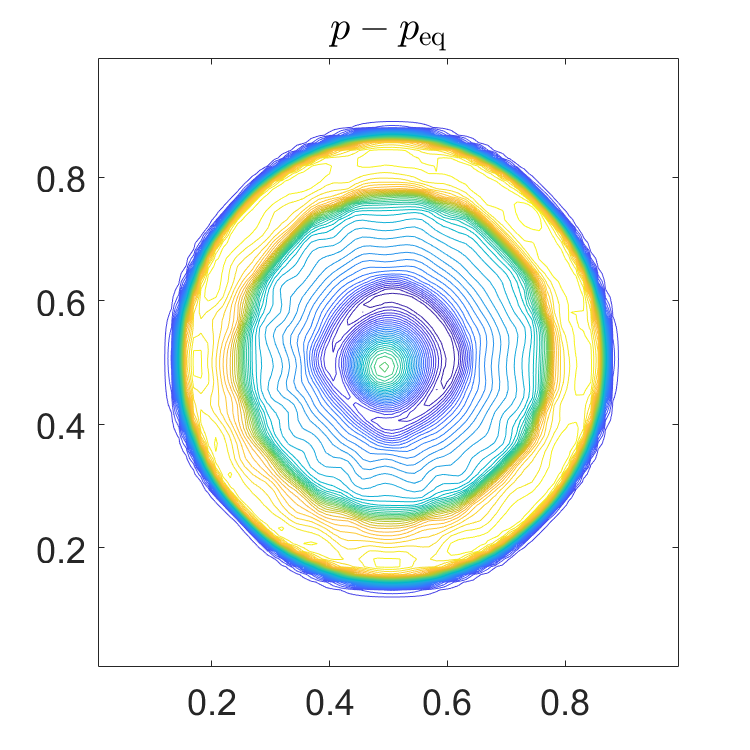}\hspace*{0.3cm}
            \includegraphics[trim=0.6cm 0.6cm 1.2cm 0.3cm, clip, width=5.3cm]{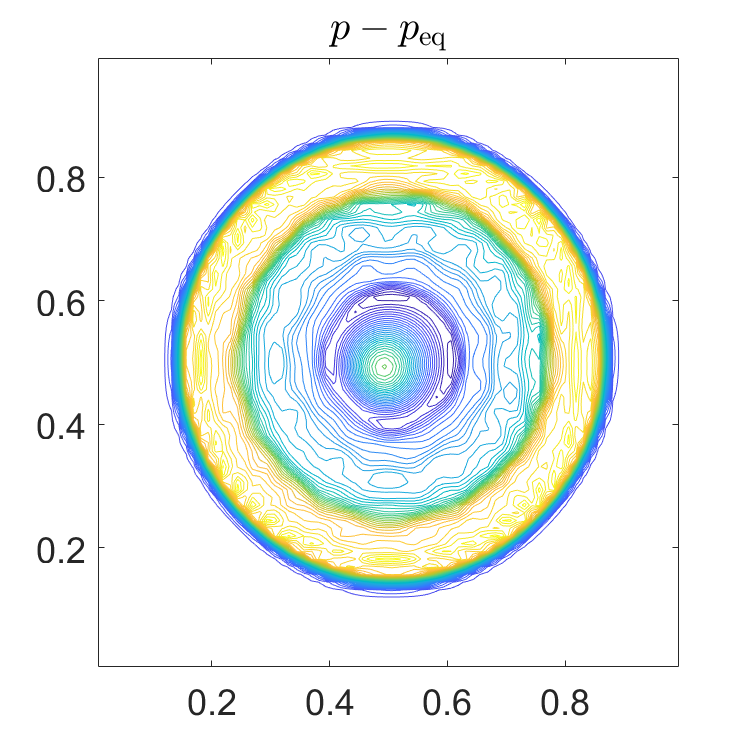}\hspace*{0.5cm}
            \includegraphics[trim=0.1cm 0.4cm 1.2cm 0.3cm, clip, width=5.5cm]{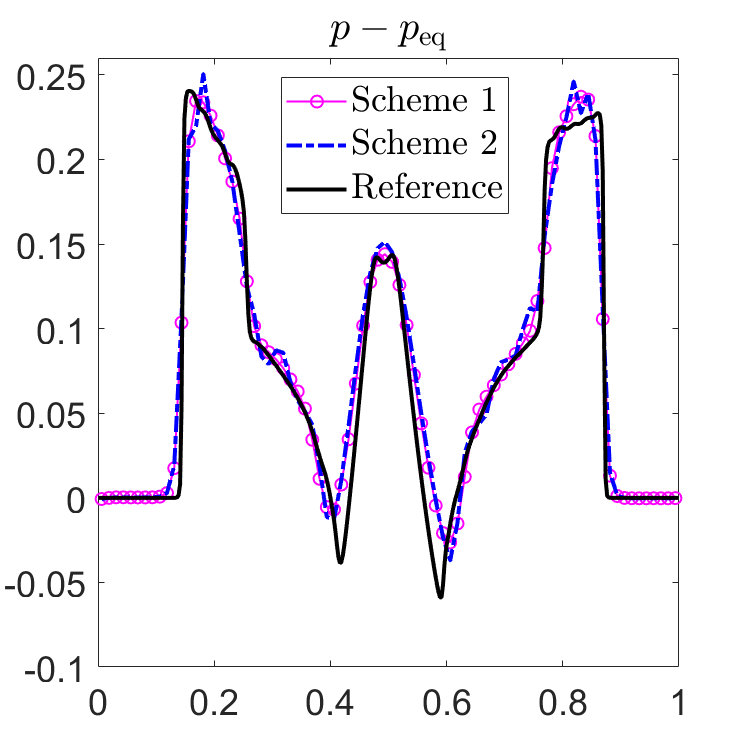}}
\caption{\sf Example 10: Pressure perturbation ($p(x,y,0.12)-p_{\rm eq}(x,y)$) captured by Schemes 1 (left) and 2 (middle) and their 1-D
slices along $y=0.5$ (right).\label{fig38}}
\end{figure}

We compute the numerical solutions by Schemes 1 and 2 on a uniform grid with $\dx=\dy=1/80$ together with the reference solution computed by
Scheme 1 on a much finer mesh with $\dx=\dy=1/400$ at time $t=0.12$ and plot the obtained results in Figure \ref{fig38}. As one can see,
Scheme 1 clearly outperforms Scheme 2 as there are almost no oscillations in the results computed by Scheme 1.
\begin{rmk}
In the above examples, we have shown that reconstructing equilibrium variables through the LCD is advantageous as it helps to remove
WENO-type oscillations while keeping the scheme WB. In the 1-D case, the proposed LCD of the equilibrium variables is based on the
projections onto the eigenvectors of the matrices $C_j=C(\mU_j)$ as explained in \S\ref{sec3}. Instead, one may try to base the LCD of
equilibrium variables on the eigenvectors of the matrices ${\cal A}_j={\cal A}(\mU_j):={\partial\mF}/{\partial\mU}(\mU_j)-B(\mU_j)$. To this
end, one needs to compute the matrices $\widehat Q_j$ and $\widehat Q_j^{-1}$ such that $\widehat Q_j^{-1}{\cal A}_j\widehat Q_j$ is a
diagonal matrix, and then to apply the Ai-WENO-Z interpolation to a different set of the local characteristic variables
$\widehat{\bm\Gamma}_\ell$, which are defined by
\begin{equation*}
\widehat{\bm\Gamma}_\ell=\widehat Q_j^{-1}\mE_\ell,\quad \ell=j\pm2,j\pm1,j.
\end{equation*}
The Ai-WENO-Z reconstruction then gives $\widehat{\bm\Gamma}^\pm_{j\mp\hf}$ and hence
$\mE^\pm_{j\mp\hf}=\widehat Q_j\bm\Gamma^\pm_{j\mp\hf}$, which are different from the values $\mE^\pm_{j\mp\hf}$ obtained in
\eref{3.3}--\eref{3.4}.

However, this kind of LCD does not lead to very good results. To demonstrate this, we recompute the numerical solutions in Examples 2, 4,
and 7 using Scheme 1, but with the aforementioned alternative LCD based on ${\cal A}(\mU)$ rather than $C(\mU)$. The obtained results shown
in Figures \ref{fig3a}--\ref{fig7a} confirm the advantages of the LCD based on $C(\mU)$.
\begin{figure}[ht!]
\centerline{\includegraphics[trim=0.0cm 0.3cm 0.7cm 0.2cm, clip, width=5.7cm]{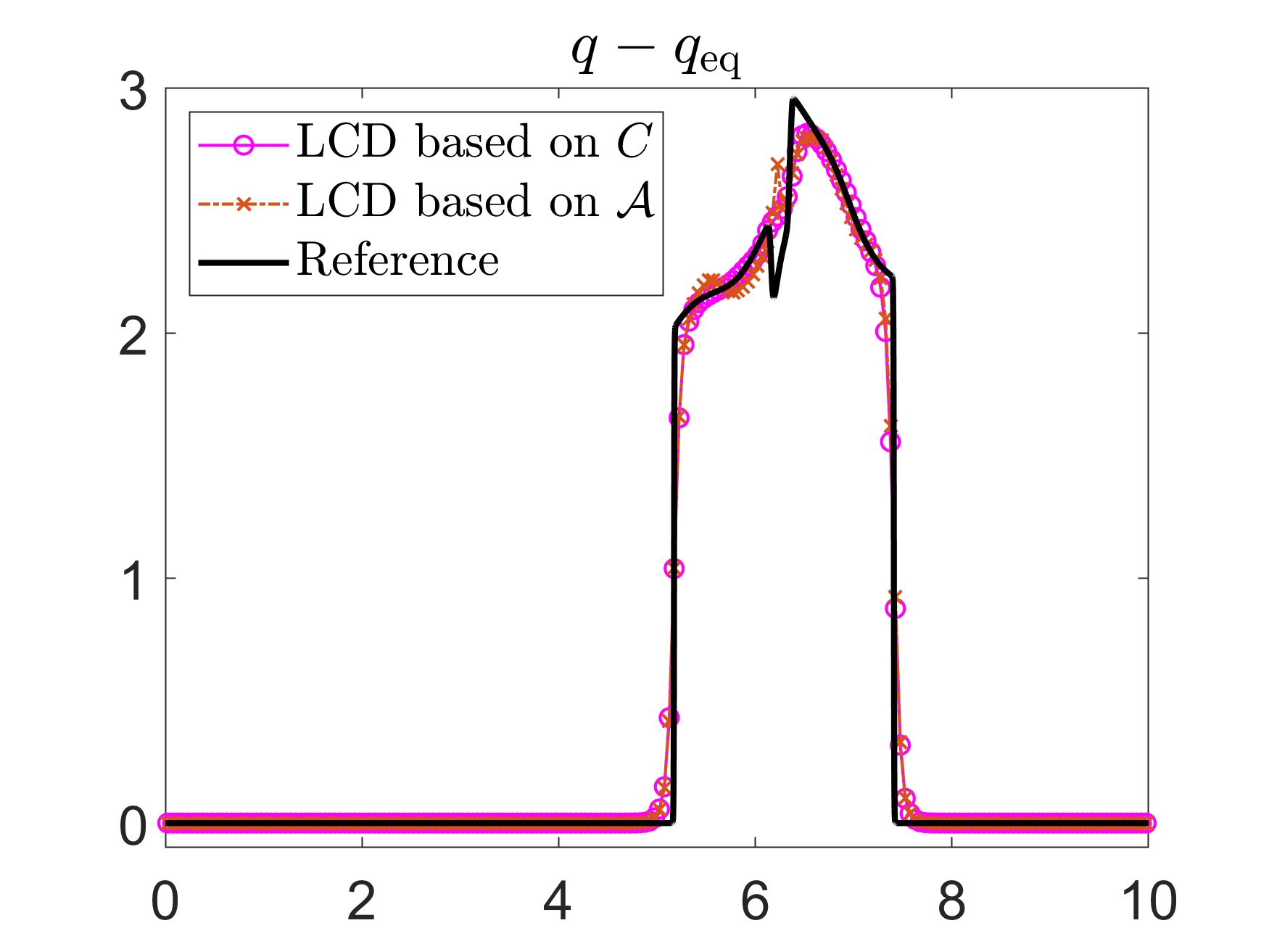}\hspace*{0.3cm}
            \includegraphics[trim=0.0cm 0.3cm 0.7cm 0.2cm, clip, width=5.7cm]{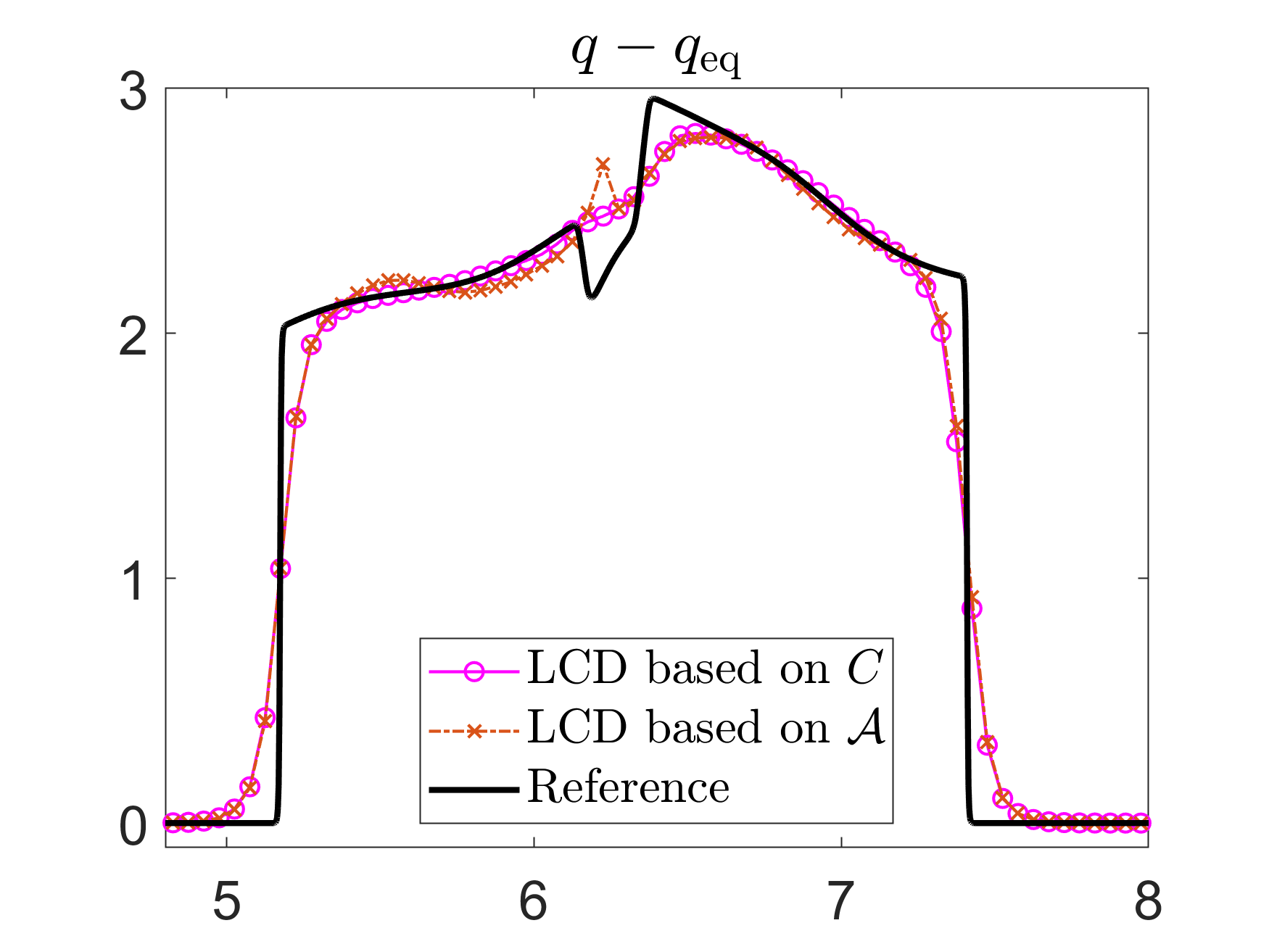}}
\caption{\sf Example 2: The differences $q(x,0.5)-q_{\rm eq}(x)$ computed by Scheme 1 with two different LCDs of the equilibrium variables
(left) and zoom at $x\in[4.8,8]$ (right). \label{fig3a}}
\end{figure}
\begin{figure}[ht!]
\centerline{\includegraphics[trim=0.0cm 0.3cm 0.9cm 0.1cm, clip, width=5.7cm]{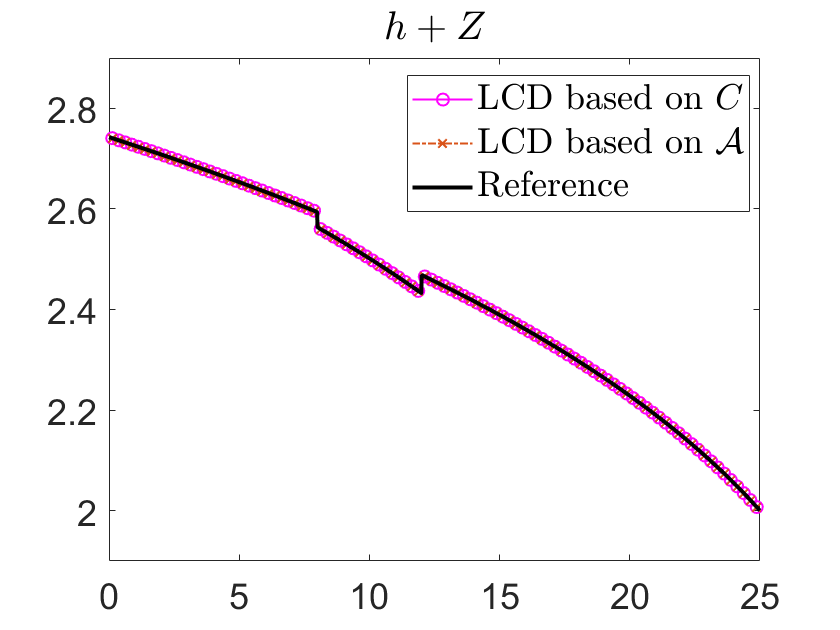}\hspace*{0.0cm}
            \includegraphics[trim=0.0cm 0.3cm 0.9cm 0.1cm, clip, width=5.7cm]{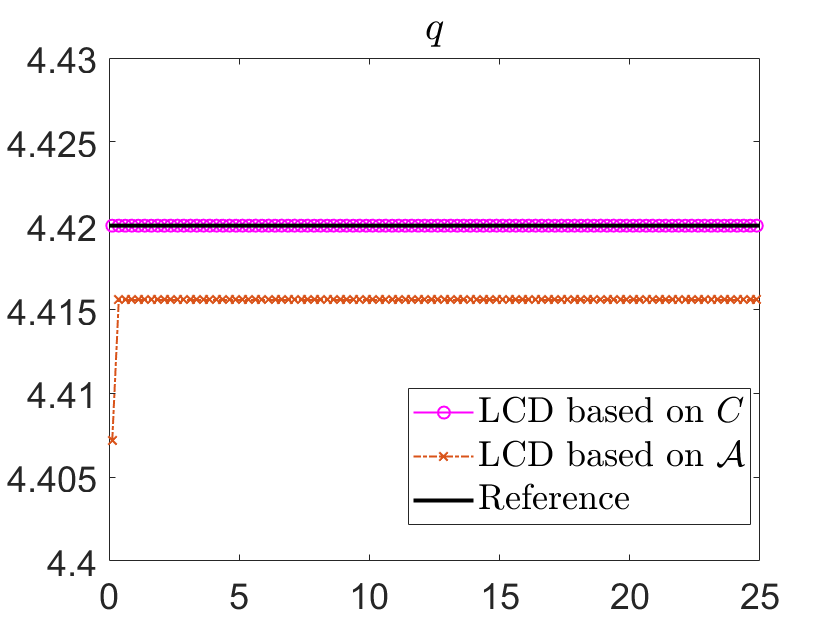}\hspace*{0.0cm}
            \includegraphics[trim=0.0cm 0.3cm 0.9cm 0.1cm, clip, width=5.7cm]{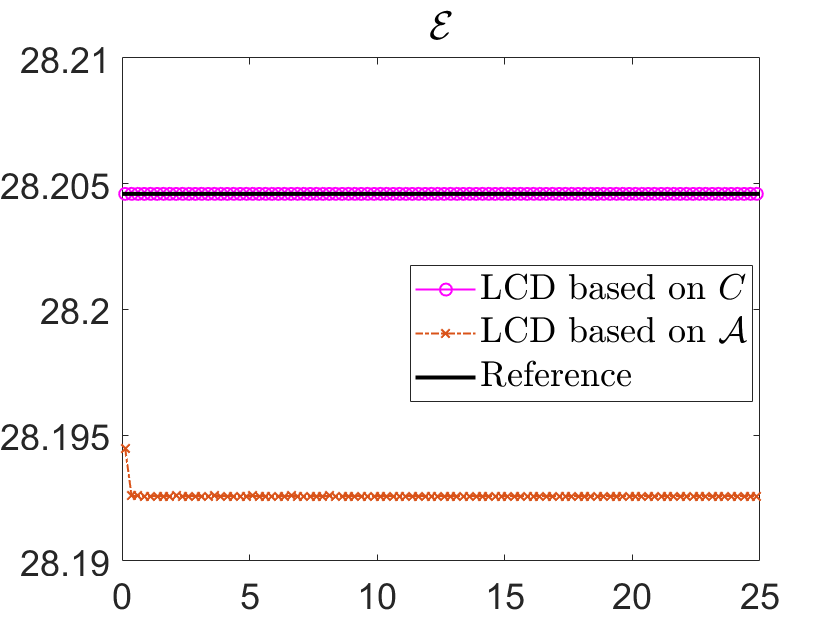}}
\caption{\sf Example 4 (moving water steady state): $h+Z$, $q$, and ${\cal E}$ computed by Scheme 1  with two different LCDs of the
equilibrium variables.\label{fig21a}}
\end{figure}
\begin{figure}[ht!]
\centerline{\includegraphics[trim=0.4cm 0.3cm 0.9cm 0.1cm, clip, width=5.7cm]{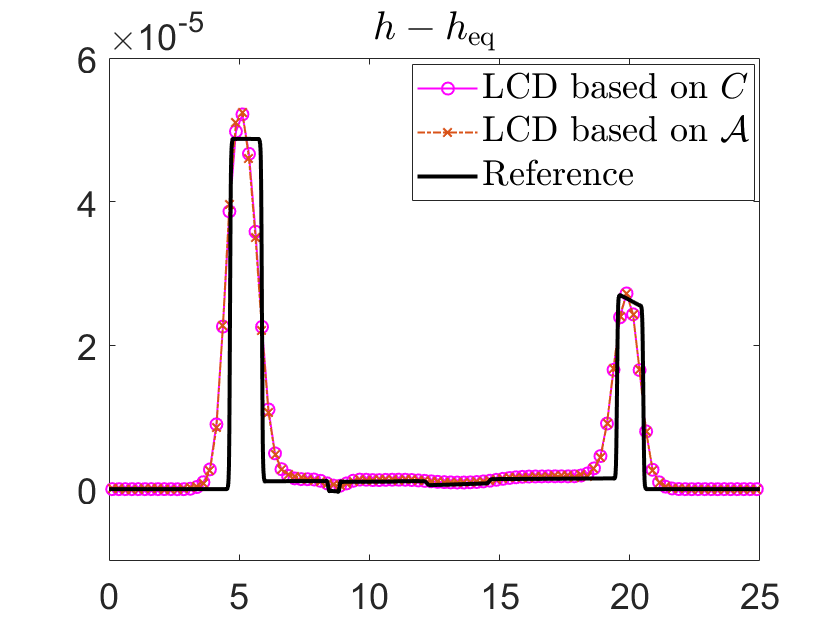}\hspace*{0.3cm}
            \includegraphics[trim=0.4cm 0.3cm 0.9cm 0.1cm, clip, width=5.7cm]{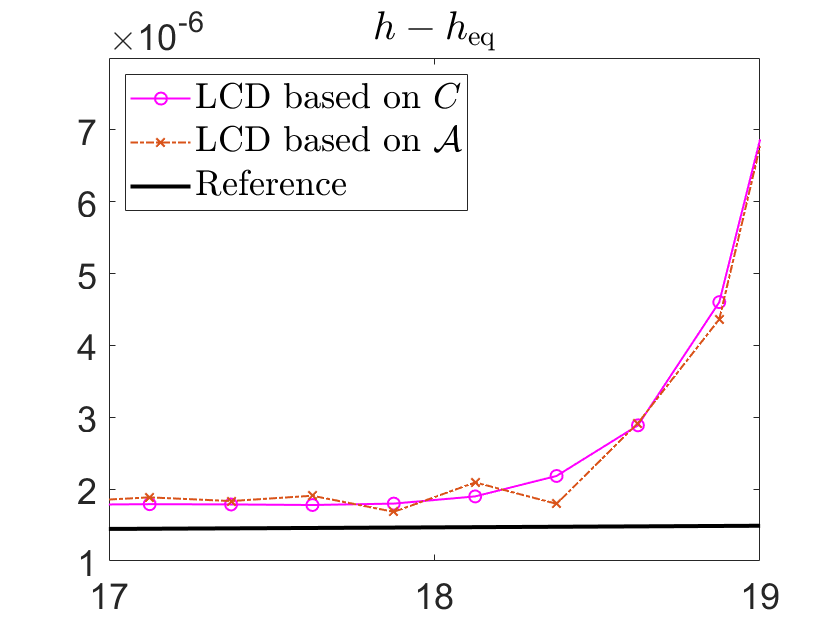}}
\caption{\sf Example 4: The differences $h(x,1.5)-h_{\rm eq}(x)$ computed by Scheme 1  with two different LCDs of the equilibrium variables
(left) and zoom at $x\in[17,19]$ (right).\label{fig22a}}
\end{figure}
\begin{figure}[ht!]
\centerline{\includegraphics[trim=0.0cm 0.3cm 0.7cm 0.2cm, clip, width=5.7cm]{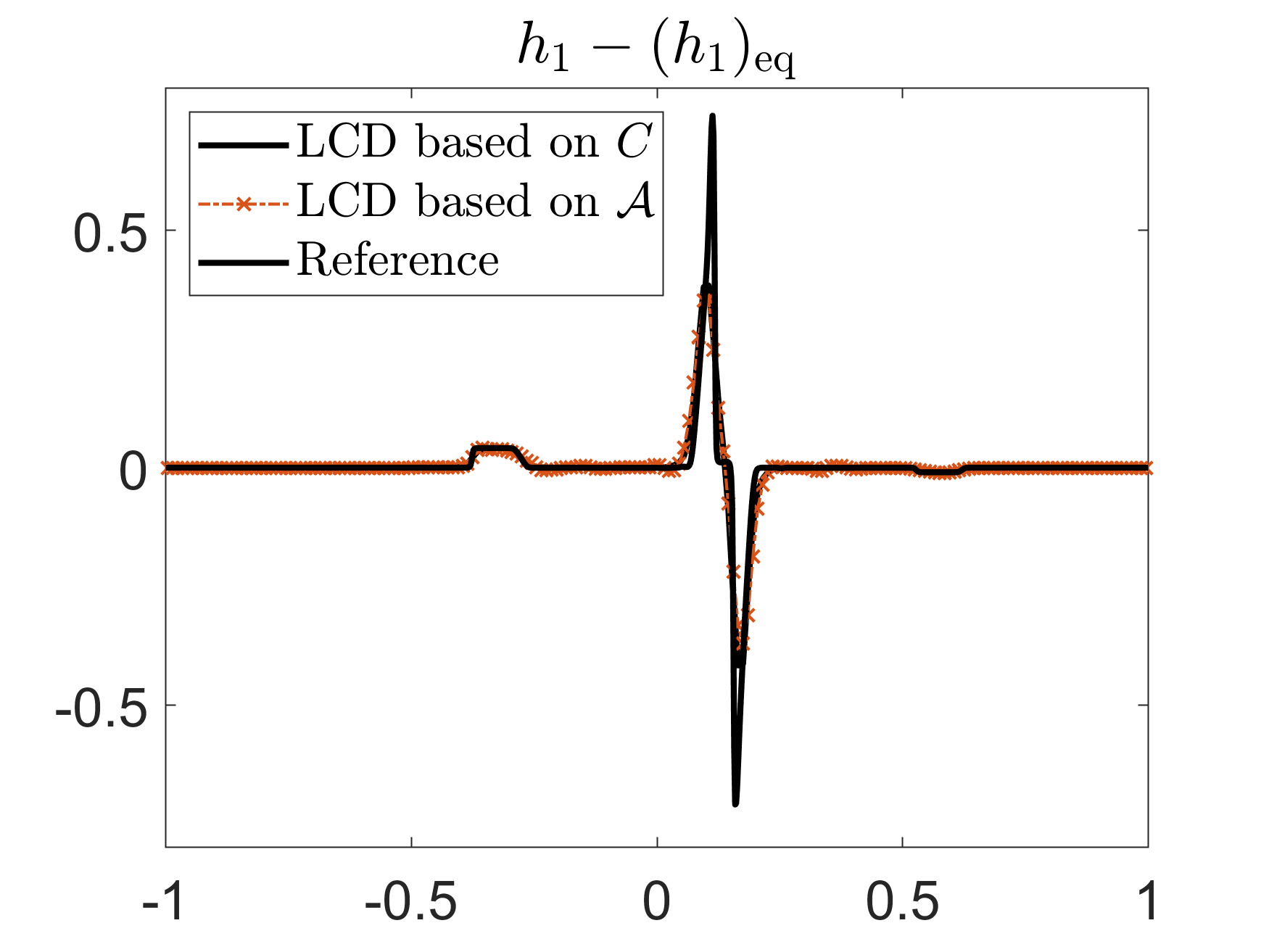}\hspace*{0.3cm}
            \includegraphics[trim=0.0cm 0.3cm 0.7cm 0.2cm, clip, width=5.7cm]{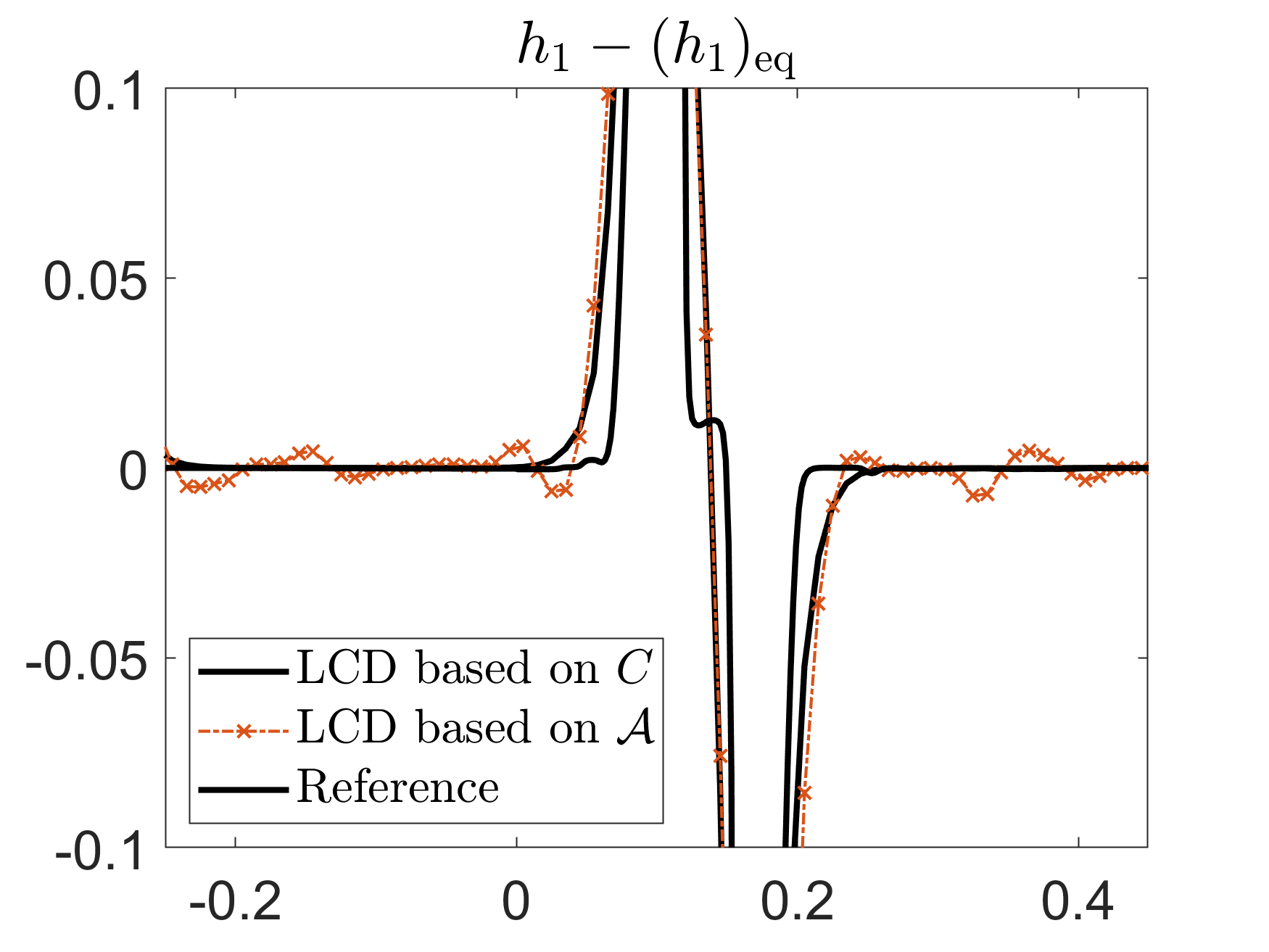}}
\caption{\sf Example 7: The differences $h_1(x,0.1)-(h_1)_{\rm eq}$ computed by Scheme 1 with two different LCDs of the equilibrium
variables (left) and zoom at $x\in[-0.25,0.45]$ (right).\label{fig7a}}
\end{figure}
\end{rmk}

\subsection{2-D Ripa System}
\subsubsection*{Example 11---Quasi 1-D Moving-Water Equilibrium and its Circular Perturbation}
In this example, we take a continuous bottom topography, which varies in the $x$-direction only:
\begin{equation*}
Z(x)=\left\{\begin{aligned}
&0.2-0.05(x-10)^2,&&8\le x\le12,\\
&0,&&\mbox{otherwise},
\end{aligned}\right.
\end{equation*}
and the following subcritical initial and boundary conditions:
\begin{equation*}
\begin{aligned}
&(\mathcal{E}^x,q^x,q^y,\theta)\Big|_{(x,y,0)}=
\big(\mathcal{E}^x_{\rm eq},q^x_{\rm eq},q^y_{\rm eq},\theta_{\rm eq}\big)\Big|_{(x,y)}\equiv(110.33025,4.42\sqrt{5},0,49.06),\\
&h(25,y,t)=2,~~q^x(0,y,t)=4.42\sqrt{5},
\end{aligned}
\end{equation*}
prescribed in the computational domain $[0,25]\times[0,10]$. The rest of the boundary conditions in the $x$-direction are homogeneous
Neumann, and the reflecting (solid wall) boundary conditions are imposed in the $y$-direction. The latter are implemented using the ghost
cell technique: we set the same values of $h$, $u$, and $\theta$ in the ghost cells while for $v$ the sign is negated.

We compute the solution by Schemes 1--3 until the final time $t=20$ on a uniform mesh with $\dx=\dy=1/4$ and present the equilibrium errors
in Table \ref{tab67}. As one can see, the errors of Schemes 1 and 2 are close to the machine errors, while Scheme 3 generates large
$L^1$-errors. We also note that the errors of Scheme 1 are smaller than those of Scheme 2, which implies that Scheme 1 performs slightly
better than Scheme 2 in maintaining the moving-water equilibrium.
\begin{table}[!ht]
\centering
\begin{tabular}{lcccccccccccccccc}
\toprule 
Scheme & ${||h-h_{\rm eq}||}_{L^1}$&${||q^x-q^x_{\rm eq}||}_{L^1}$&${||q^y-q^y_{\rm eq}||}_{L^1}$&${||h\theta-h\theta_{\rm eq}||}_{L^1}$&
${||{\cal E}^x-{\cal E}^x_{\rm eq}||}_{L^1}$\\[0.8ex]
\hline
Scheme 1&1.83e-15&1.36e-14&0.00&1.26e-13&8.60e-14\\
Scheme 2&2.17e-15&2.54e-14&0.00&3.49e-13&2.51e-13\\
Scheme 3&2.51e-03&1.09e-03&0.00&1.23e-01&8.66e-02\\
\bottomrule 
\end{tabular}
\caption{Example 11: $L^1$-errors in $h$, $q^x$, $q^y$, $h\theta$, and ${\cal E}^x$ at $t=20$.}\label{tab67}
\end{table}

We then add a small genuinely 2-D circular perturbation to the water depth without changing the initial $q^x$, $q^y$, and $h\theta$. The
initially perturbed $h$ is
\begin{equation*}
h(x,y,0)=h_{\rm eq}(x,y)+
\left\{\begin{aligned}
&0.01,&&(x-6)^2+(y-5)^2<0.25,\\
&0,&&\mbox{otherwise}.
\end{aligned}\right.
\end{equation*}

\begin{figure}[ht!]
\centerline{\includegraphics[trim=0.0cm 0.3cm 0.7cm 0.2cm, clip, width=5.7cm]{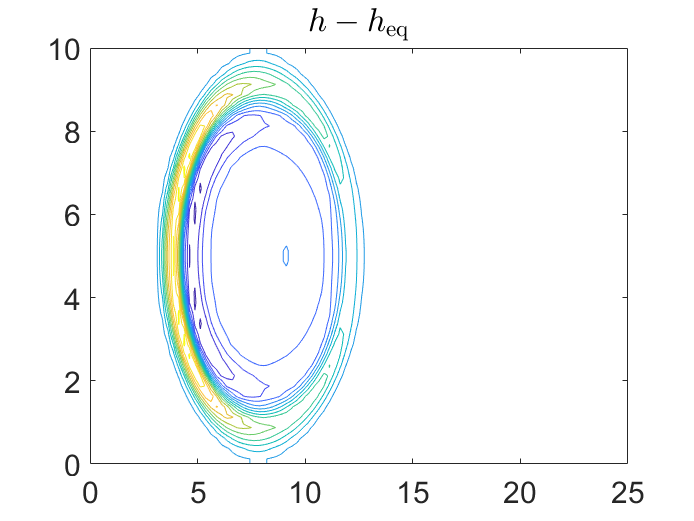}\hspace{0.3cm}
            \includegraphics[trim=0.0cm 0.3cm 0.7cm 0.2cm, clip, width=5.7cm]{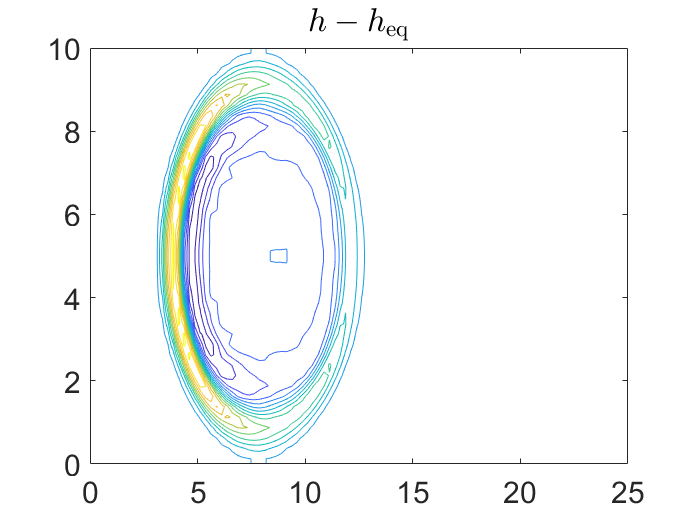}}
\vskip7pt
\centerline{\includegraphics[trim=0.0cm 0.3cm 0.7cm 0.2cm, clip, width=5.7cm]{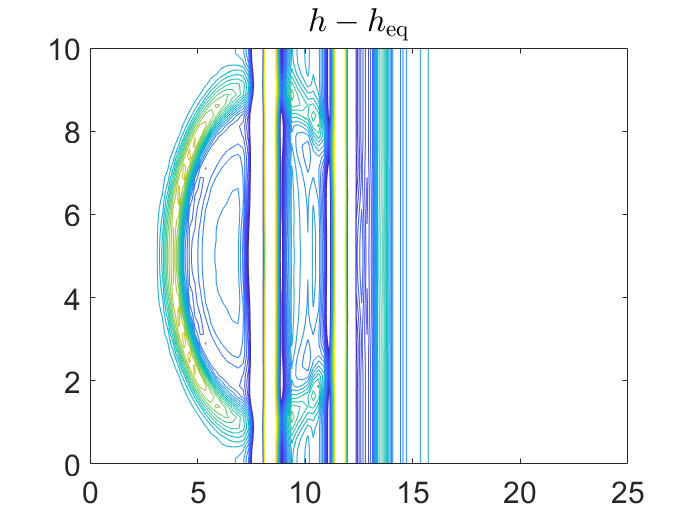}\hspace{0.3cm}
            \includegraphics[trim=0.0cm 0.3cm 0.7cm 0.2cm, clip, width=5.7cm]{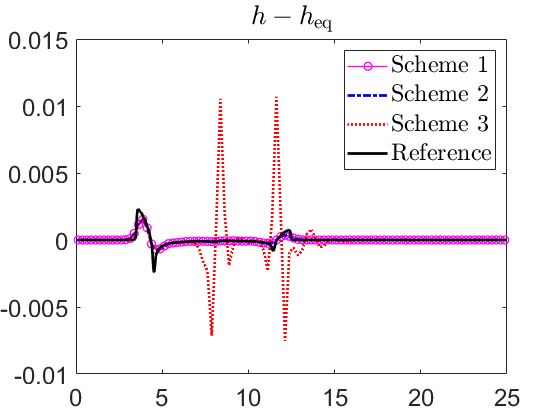}}
\caption{\sf Example 11 (circular perturbation): The differences $h(x,0.4)-h_{\rm eq}$ computed by Schemes 1 (top left), 2 (top right), and
3 (bottom left), and their 1-D slices along the line $y=5.125$ (bottom right).\label{figRipa}}
\end{figure}

We compute the numerical solutions by Schemes 1--3 until the final time $t=0.4$ on a uniform mesh with $\dx=\dy=1/4$ together with the
reference solution computed by Scheme 1 on a much finer mesh with $\dx=\dy=1/40$ and plot the obtained results in Figure \ref{figRipa}. As
one can see, Scheme 3, which cannot preserve the moving-water equilibrium, produces large spurious oscillations near $x=8$ and $x=12$ due to
the effect of the bottom topography. On the other hand, the solutions computed by Schemes 1 and 2 are almost oscillation-free.

\subsubsection*{Example 12---2-D Dam-Break Problem}
In the final example, we consider the initial data
\begin{equation*}
(h,u,v,\theta)\Big|_{(x,y,0)}=
\left\{\begin{aligned}
&(5,0.5,0,9.812),&&\max\{|x|,|y|\}\le0.5,\\
&(3,2.75,0,15.2086),&&\quad~~~\mbox{otherwise},
\end{aligned}\right.
\end{equation*}
prescribed in the computational domain $[-1,1]\times[-1,1]$ subject to the homogeneous Neumann boundary conditions. The bottom topography
is flat ($Z(x,y)\equiv0$).

\begin{figure}[ht!]
\centerline{\includegraphics[trim=0.0cm 0.3cm 0.7cm 0.15cm, clip, width=5.7cm]{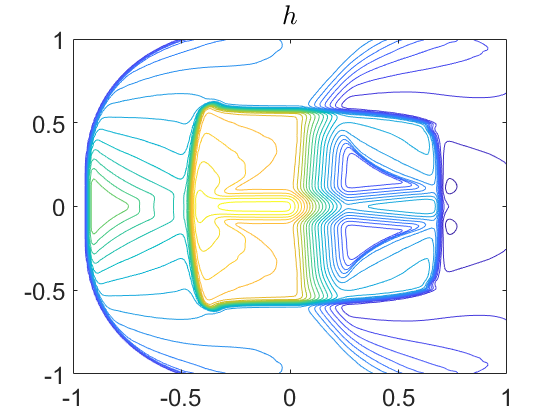}\hspace{0.2cm}
            \includegraphics[trim=0.0cm 0.3cm 0.7cm 0.15cm, clip, width=5.7cm]{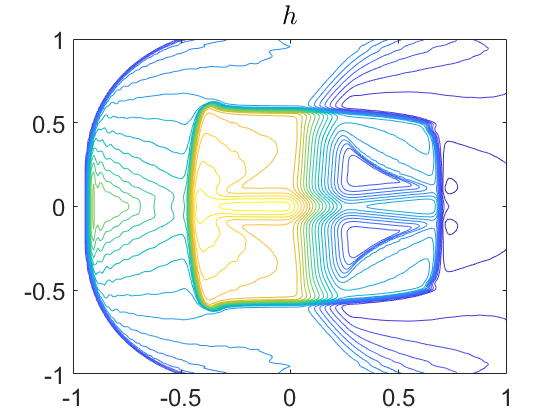}\hspace{0.2cm}
            \includegraphics[trim=0.0cm 0.3cm 0.7cm 0.15cm, clip, width=5.7cm]{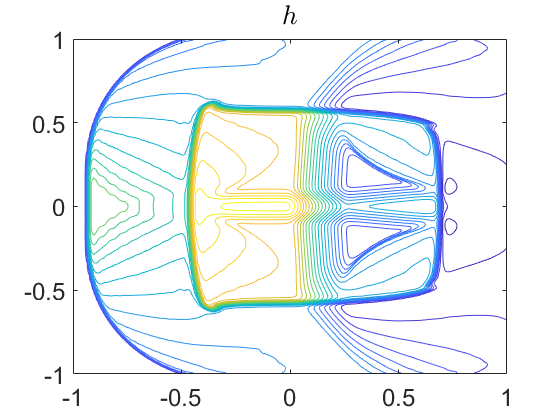}}
\vskip7pt
\centerline{\includegraphics[trim=0.0cm 0.3cm 0.7cm 0.15cm, clip, width=5.7cm]{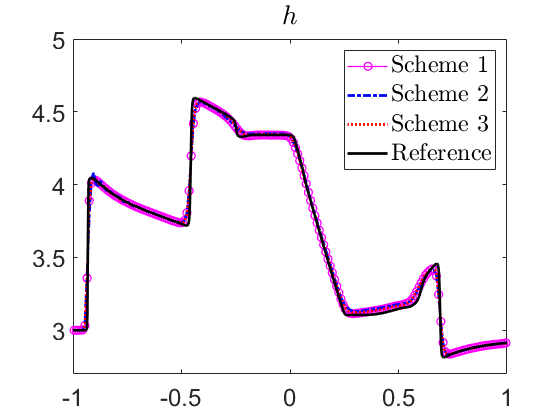}\hspace{0.3cm}
            \includegraphics[trim=0.0cm 0.3cm 0.7cm 0.15cm, clip, width=5.7cm]{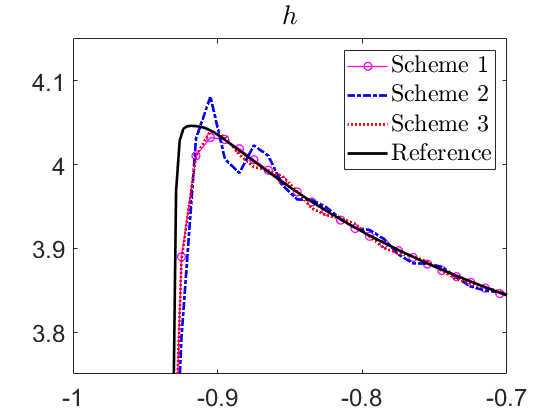}}
\caption{\sf Example 12: Top row: water depth $h$ computed by Schemes 1 (left), 2 (middle), and 3 (right). Bottom row: 1-D slices of the
computed $h$ along the line $y=0.125$ (left) and zoom at $x\in[-1,-0.7]$ (right).\label{figDam-Break}}
\end{figure}

We compute the solutions by Schemes 1--3 until the final time $t=0.075$ on a uniform mesh with $\dx=\dy=1/100$ and plot the obtained results
in Figure \ref{figDam-Break}. The reference solution is computed by Scheme 1 on a much finer mesh with $\dx=\dy=1/400$. As one can see,
Scheme 2, which does not use the LCD, generates spurious oscillations, which can be clearly seen in the 1-D slices along the line $y=0.125$.

\section*{Acknowledgment}
The work of S. Chu was funded by the Deutsche Forschungsgemeinschaft (DFG, German Research Foundation) - SPP 2410 Hyperbolic Balance Laws in
Fluid Mechanics: Complexity, Scales, Randomness (CoScaRa) within the Project(s) HE5386/26-1 (Numerische Verfahren f\"ur gekoppelte
Mehrskalenprobleme,525842915) and (Zuf\"allige kompressible Euler Gleichungen: Numerik und ihre Analysis, 525853336) HE5386/27-1, and the
Deutsche Forschungsgemeinschaft (DFG, German Research Foundation) - SPP 2183: Eigenschaftsgeregelte Umformprozesse with the
Project(s) HE5386/19-2,19-3 Entwicklung eines flexiblen isothermen Reckschmiedeprozesses f\"ur die eigenschaftsgeregelte Herstellung von
Turbinenschaufeln aus Hochtemperaturwerkstoffen (424334423). The work of A. Kurganov was supported in part by NSFC grants 12171226 and
W2431004. The work of B.-S. Wang was supported in part by NSFC grant 12301530 and the startup funding provided by the Ocean University of
China.


\begin{thebibliography}{10}

\bibitem{CKL23}
{\sc Y.~Cao, A.~Kurganov, and Y.~Liu}, Flux globalization based well-balanced
  path-conservative central-upwind scheme for the thermal rotating shallow
  water equations, Commun. Comput. Phys., 34 (2023), pp.~993--1042.

\bibitem{CKLX_22}
{\sc Y.~Cao, A.~Kurganov, Y.~Liu, and R.~Xin}, Flux globalization based
  well-balanced path-conservative central-upwind schemes for shallow water
  models, J. Sci. Comput., 92 (2022).
\newblock Paper No. 69.

\bibitem{CCK23_Adaptive}
{\sc A.~Chertock, S.~Chu, and A.~Kurganov}, Adaptive high-order {A-WENO}
  schemes based on a new local smoothness indicator, E. Asian. J. Appl. Math.,
  13 (2023), pp.~576--609.

\bibitem{CCKOT}
{\sc A.~Chertock, S.~Cui, A.~Kurganov, {\c{S}}.~N. \"{O}zcan, and E.~Tadmor},
  Well-balanced schemes for the {E}uler equations with gravitation:
  {C}onservative formulation using global fluxes, J. Comput. Phys., 358 (2018),
  pp.~36--52.

\bibitem{Chu21}
{\sc S.~Chu, A.~Kurganov, and M.~Na}, Fifth-order {A}-{WENO} schemes based on
  the path-conservative central-upwind method, J. Comput. Phys., 469 (2022).
\newblock Paper No. 111508.

\bibitem{CKX23}
{\sc S.~Chu, A.~Kurganov, and R.~Xin}, New {M}ore {E}fficient {A}-{WENO}
  {S}chemes, J. Sci. Comput., 104 (2025).
\newblock Paper No. 53.

\bibitem{CKX_24WB}
{\sc S.~Chu, A.~Kurganov, and R.~Xin}, A well-balanced fifth-order {A-WENO}
  scheme based on flux globalization, Beijing J. Pure Appl. Math., 2 (2025),
  pp.~87--113.

\bibitem{DLGW}
{\sc W.~S. Don, D.-M. Li, Z.~Gao, and B.-S. Wang}, A characteristic-wise
  alternative {WENO}-{Z} finite difference scheme for solving the compressible
  multicomponent non-reactive flows in the overestimated quasi-conservative
  form, J. Sci. Comput., 82 (2020).
\newblock Paper No. 27.

\bibitem{DLWW22}
{\sc W.~S. Don, R.~Li, B.-S. Wang, and Y.~H. Wang}, A novel and robust
  scale-invariant {WENO} scheme for hyperbolic conservation laws, J. Comput.
  Phys., 448 (2022).
\newblock Paper No. 110724.

\bibitem{Gottlieb11}
{\sc S.~Gottlieb, D.~Ketcheson, and C.-W. Shu}, Strong stability preserving
  {R}unge-{K}utta and multistep time discretizations, World Scientific
  Publishing Co. Pte. Ltd., Hackensack, NJ, 2011.

\bibitem{Gottlieb12}
{\sc S.~Gottlieb, C.-W. Shu, and E.~Tadmor}, Strong stability-preserving
  high-order time discretization methods, SIAM Rev., 43 (2001), pp.~89--112.

\bibitem{GrosheintzLavalKaeppeli2019}
{\sc L.~Grosheintz-Laval and R.~K{\"a}ppeli}, High-order well-balanced finite
  volume schemes for the {Euler} equations with gravitation, J. Comput. Phys.,
  378 (2019), pp.~324--343.

\bibitem{JSZ}
{\sc Y.~Jiang, C.-W. Shu, and M.~Zhang}, An alternative formulation of finite
  difference weighted {ENO} schemes with {L}ax-{W}endroff time discretization
  for conservation laws, SIAM J. Sci. Comput., 35 (2013), pp.~A1137--A1160.

\bibitem{Joh}
{\sc E.~Johnsen}, On the treatment of contact discontinuities using {WENO}
  schemes, J. Comput. Phys., 230 (2011), pp.~8665--8668.

\bibitem{KaeppeliMishra2014}
{\sc R.~K{\"a}ppeli and S.~Mishra}, Well-balanced schemes for the {Euler}
  equations with gravitation, J. Comput. Phys., 259 (2014), pp.~199--219.

\bibitem{KlingenbergPuppoSemplice2019}
{\sc C.~Klingenberg, G.~Puppo, and M.~Semplice}, Arbitrary order finite volume
  well-balanced schemes for the {Euler} equations with gravity, SIAM J. Sci.
  Comput., 41 (2019), pp.~A695--A721.

\bibitem{KLX_21}
{\sc A.~Kurganov, Y.~Liu, and R.~Xin}, Well-balanced path-conservative
  central-upwind schemes based on flux globalization, J. Comput. Phys., 474
  (2023).
\newblock Paper No. 111773.

\bibitem{Kurganov25Na}
{\sc A.~Kurganov and M.~Na}, Flux globalization based well-balanced
  central-upwind schemes for the {E}uler equations with gravitation, Computers
  \& Fluids, 300 (2025).
\newblock Paper No. 106713.

\bibitem{KPmultil}
{\sc A.~Kurganov and G.~Petrova}, Central-upwind schemes for two-layer shallow
  water equations, SIAM J. Sci. Comput., 31 (2009), pp.~1742--1773.

\bibitem{LiXing2016}
{\sc G.~Li and Y.~Xing}, High order finite volume {WENO} schemes for the
  {Euler} equations under gravitational fields, J. Comput. Phys., 316 (2016),
  pp.~145--163.

\bibitem{LiXing2018}
{\sc G.~Li and Y.~Xing}, Well-balanced finite difference weighted essentially
  non-oscillatory schemes for the {Euler} equations with static gravitational
  fields, Comput. Math. Appl., 75 (2018), pp.~2071--2085.

\bibitem{LLWDW23}
{\sc P.~Li, T.~T. Li, W.~S. Don, and B.-S. Wang}, Scale-invariant
  multi-resolution alternative {WENO} scheme for the {E}uler equations, J. Sci.
  Comput., 94 (2023).
\newblock Paper No. 15.

\bibitem{LiWangDon2021}
{\sc P.~Li, B.-S. Wang, and W.-S. Don}, Sensitivity parameter-independent
  characteristic-wise well-balanced finite volume {WENO} scheme for the {Euler}
  equations under gravitational fields, J. Sci. Comput., 88 (2021).
\newblock Paper No. 47.

\bibitem{Liu17}
{\sc H.~Liu}, A numerical study of the performance of alternative weighted
  {ENO} methods based on various numerical fluxes for conservation law, Appl.
  Math. Comput., 296 (2017), pp.~182--197.

\bibitem{Luo11}
{\sc J.~Luo, K.~Xu, and N.~Liu}, A well-balanced symplecticity-preserving
  gas-kinetic scheme for hydrodynamic equations under gravitational field, SIAM
  J. Sci. Comput., 33 (2011), pp.~2356--2381.

\bibitem{Man91}
{\sc R.~Manning}, On the flow of water in open channels and pipes, Trans. Inst.
  Civ. Eng. Irel., 20 (1891), pp.~161--207.

\bibitem{Nonomura20}
{\sc T.~Nonomura and K.~Fujii}, Characteristic finite-difference {WENO} scheme
  for multicomponent compressible fluid analysis: overestimated
  quasi-conservative formulation maintaining equilibriums of velocity,
  pressure, and temperature, J. Comput. Phys., 340 (2017), pp.~358--388.

\bibitem{Qiu02}
{\sc J.~Qiu and C.-W. Shu}, On the construction, comparison, and local
  characteristic decomposition for high-order central {WENO} schemes, J.
  Comput. Phys., 183 (2002), pp.~187--209.

\bibitem{QGKWW}
{\sc Y.~Qiu, Z.~Gao, A.~Kurganov, B.-S. Wang, and X.~Wen}, Fifth-order
  well-balanced path-conservative {A-WENO} scheme for the {R}ipa model.
\newblock Submitted; arXiv:2607.09293.

\bibitem{Rip93}
{\sc P.~Ripa}, Conservation laws for primitive equations models with
  inhomogeneous layers, Geophys. Astrophys. Fluid Dynam., 70 (1993),
  pp.~85--111.

\bibitem{Rip95}
{\sc P.~Ripa}, On improving a one-layer ocean model with thermodynamics, J.
  Fluid Mech., 303 (1995), pp.~169--201.

\bibitem{Shu20}
{\sc C.-W. Shu}, Essentially non-oscillatory and weighted essentially
  non-oscillatory schemes, Acta Numer., 29 (2020), pp.~701--762.

\bibitem{WD22}
{\sc B.-S. Wang and W.~S. Don}, Affine-invariant {WENO} weights and operator,
  Appl. Numer. Math., 181 (2022), pp.~630--646.

\bibitem{WDGK_20}
{\sc B.-S. Wang, W.~S. Don, N.~K. Garg, and A.~Kurganov}, Fifth-order {A-WENO}
  finite-difference schemes based on a new adaptive diffusion central numerical
  flux, SIAM J. Sci. Comput., 42 (2020), pp.~A3932--A3956.

\bibitem{WDKL}
{\sc B.-S. Wang, W.~S. Don, A.~Kurganov, and Y.~Liu}, Fifth-order {A-WENO}
  schemes based on the adaptive diffusion central-upwind {R}ankine-{H}ugoniot
  fluxes, Commun. Appl. Math. Comput., 5 (2023), pp.~295--314.

\bibitem{WLGD_18}
{\sc B.-S. Wang, P.~Li, Z.~Gao, and W.~S. Don}, An improved fifth order
  alternative {WENO-Z} finite difference scheme for hyperbolic conservation
  laws, J. Comput. Phys., 374 (2018), pp.~469--477.

\bibitem{Xing13}
{\sc Y.~Xing and C.-W. Shu}, High order well-balanced {WENO} scheme for the gas
  dynamics equations under gravitational fields, J. Sci. Comput., 54 (2013),
  pp.~645--662.

\bibitem{XuShu2024}
{\sc Z.~Xu and C.-W. Shu}, Local characteristic decomposition-free high-order
  finite difference {WENO} schemes for hyperbolic systems endowed with a
  coordinate system of {Riemann} invariants, SIAM J. Sci. Comput., 46 (2024),
  pp.~A1352--A1372.

\end{thebibliography}

\end{document}